\documentclass[10pt]{article}

\usepackage[english]{babel}
\usepackage{AVDoc_OM}

\title{ Adaptation of a population to a changing environment\\
under the light of quasi-stationarity
 }
\author{Aurélien Velleret \footnote{Institut für Mathematik, Goethe Universität, Fachbereich 12, 60054 Frankfurt am Main, Germany,
		email: velleret@math.uni-frankfurt.de}}
\begin{document}

\maketitle


\section*{Abstract}
\setcounter{eq}{0}

We analyze the long-term stability 
of a stochastic model
designed to illustrate 
the adaptation of a population to variation in its environment.
A piecewise-deterministic process modeling adaptation
is coupled to a Feller logistic diffusion 
modeling  population size.
As the individual features in the population 
become further away from the optimal ones,
the growth rate declines,
making population extinction  more likely.
Assuming that the environment 
changes  deterministically and steadily 
in a constant direction,
we obtain the existence and uniqueness 
of the quasi-stationary distribution,
the associated survival capacity and the Q-process.
Our approach also provides
several exponential convergence results 
(in total variation for the measures).
From this synthetic information,
we can characterize the efficiency 
of internal adaptation
(i.e. population turnover from mutant invasions).
When the latter is lacking, 
there is still stability,
but because of the high level of population extinction.
Therefore, such a characterization 
must be based on specific features of this quasi-ergodic regime.
\\

Keywords: mobile optimum, quasi-stationary distribution, evolution, ecology, jump processes, Markov process in continuous time and continuous space

\section{Introduction}
\setcounter{eq}{0}

\subsection{Eco-evolutionary motivations}

Our objective is to study the relative contribution 
of mutations with various strong effects
 to the adaptation of a population.
Our goal is therefore to analyze a model as simple as possible 
 in which these mutations are filtered 
according to the advantage they provide.
 This advantage can be immediately significant
 (better growth rate of the mutant subpopulation)
 or play a role in the future adaptation
 (the population is doomed without mutants). 
 The stochastic model considered
  takes into account these two aspects.
It extends the one introduced by \cite{KH09}
and described  more formally 
in \cite{NP17} and \cite{KNP18}. 

Similarly,
we assume that the population is described by a certain value $\hat{x} \in \bR^d$,
hereafter referred to as its trait.
For the sake of a simple theoretical model, 
spatial dispersion as well as phenotypic heterogeneity 
(at least for the individual features of interest)
are neglected.
We therefore that the population is monomorphic at all times
 and that $\hat{x}$ represents
 the phenotype of the individuals in the population.
Nonetheless, we allow for variations of this trait $\hat{x}$
due to stochastic events, 
namely when a subpopulation issued from a mutant with trait $\hat{x} + w$ 
manages to persist and invade the "resident" population.
In the model, such events are assumed
 to occur instantaneously.

The main novelty of our approach
is that we couple this "adaptive" process
with a Feller diffusion process $N$
with a logistic drift.
This diffusion describes the dynamics of the population size 
in a limit where it is large.
We mean here that individual birth and death events 
have a negligible impact,
but that the accumulation of these events
has a visible and stochastic effect.
 In particular, the introduction of the "size" in the model
allows us to easily translate
the notion of maladaptation, 
in the form of a poor growth rate.

For the long-time dynamics,
we are mainly interested in considering 
only surviving populations,
that is conditioning the process
upon the fact that the population size has not decreased to 0.
The implication of taking size into account 
is twofold.
On the one hand,
extinction occurs much  more rapidly 
when adaptation is poor.
Indeed, the population size is then very rapidly declining.
So a natural selection effect 
can be observed at the population level. 
On the other hand,
the better the adaptation,
the larger the population size can be
and the more frequent the birth of new mutants
 in the population.
In our simple model, 
a mutant trait that is better suited 
for the survival of the population as a whole
is also characterized by a greater probability
that the resident population gets invaded,
once a single mutant is introduced.
Compared to the case of a fixed size 
as in \cite{NP17}  and \cite{KNP18},
this second implication means a stabilizing effect for the phenotype
when the population size is large enough;
but also a destabilizing effect when the population size decreases.
This is in contrast to natural selection 
at the individual level 
(which is the main effect detailed in \cite{KH09}).
Indeed, when adaptation is already nearly optimal,
very few among the mutants that appear in the population
can successfully maintain themselves 
and eventually invade the resident trait.

Let us assume here that mutations 
can allow the individuals to survive in these new environments.
In this context, 
how resilient is the population 
to environmental changes?
Is there a clear threshold to the rate of change
that such a population can handle?
How can we describe the interplay between the  above properties?

To begin to answer these questions,
and like \cite{KH09},
we assume for simplicity
 that the environmental change 
is given by a constant speed translation 
of the profile of fitness,
with $v$ this speed and 
$\mathbf{e_1}$ the direction of this change. 
In practice, this means that the growth rate 
of the population at time $t$
is expressed as a function of $x:= \hat{x} - v\,t\, \mathbf{e_1}$,
for a monomorphic population with trait $\hat{x}$ at time $t$.
Naturally, 
the phenotypic lag $x$ 
becomes the main quantity of interest for varying $t$.
Likewise, 
we can express as a function of $x$ and $w$
the probability that a mutant individual,
with mutation $w$,
will lead to the invasion of a resident population
 with trait $\hat{x}$ at time $t$.
This probability should indeed be related 
to the difference between the growth rate at $x$ and at $x+w$,
although we will not require any precise relationship 
in the results of this paper.
Furthermore, we assume
that the distribution of the additive effect for the new mutations 
is constant over time 
and independent of the trait $\hat{x}$ of the population 
before the mutation
(thus independent of $x$ in the moving frame of reference).
\\

In this context, 
we can exploit the notion of quasi-stationary distribution (QSD, cf Remark 2.2.3)
to characterize what would be an equilibrium for these dynamics 
prior to extinction.
The main contribution of the current paper 
is to ensure that this notion is  unambiguously defined here.
To the best of our knowledge, 
this is the first time 
that the existence and uniqueness of the QSD is proved
for a piecewise deterministic process 
coupled to a diffusion.

By our proof, 
we also provide a justification 
of the notion of typical relaxation time and extinction time.
The quasi-stationary description is well suited
provided the latter is much longer than the former.
As can be verified by simulations, 
typical convergence to the QSD is exponential in such cases.
However, 
the marginal starting from certain initial conditions 
may take long before it approximates the QSD,
mainly in cases where extinction is initially very likely.

In the following subsections of the introduction,
we present the stochastic process under consideration
then some elementary notations.
The main results are described in Section \ref{O_Res},
starting with our hypothesis in Subsection \ref{O_hyp2} 
and the theorem statement 
in Section \ref{O_Thm}.
In Subsection \ref{O_sec_ecoevo},
we discuss its interpretation in terms of ecology and evolution.
Its connection with related adaptation models  
is given in Subsection \ref{O_QErel},
and with the classical techniques of quasi-stationarity in Subsection \ref{O_MQS}.
The rest of the paper is devoted to  proofs.
We prove the existence and uniqueness of the process in Section \ref{O_exU},
and introduce in the next Section \ref{O_str}
the main theorems on which 
our main Theorem \ref{O_ECVGnl} is based.
Two alternative hypotheses are considered,
 with some variations in the proofs.
We choose to group the theorems in the three following sections
according to the property of the process they imply 
for the various sets of assumptions.
The definition of a specific sigma-field
and its property are reported in the Appendix,
as well as some illustrations of the asymptotic profiles
given by simulations.

\subsection{The stochastic model}
\label{O_model}

Following \cite{KH09} as explained in the introduction
for the definition of the adaptive component,
the system 
that describes the combined evolution 
of the population size and of its phenotypic lag
is then given by:
\begin{equation*}
(S_0)\left\{
\begin{aligned}
X_t 
&= x - v\, t\, \mathbf{e_1}
+ \int_{[0,t] \times \bR^d \times (\bR_+)^2 } 
w \; \varphi_0 \left( X_{s^-},\,N_s,\,w,\,u_f,\, u_g \right) \Md\\ 
N_t 
&= n + \int_{0}^{t}  \left( r(X_s) \; N_s - \gamma_0 \ltm (N_s)^2 \right) ds 
+ \sigma \int_{0}^{t} \sqrt{N_s} \; dB_s,
\end{aligned}
\right.
\end{equation*}
where $N_t$ describes the size of the population and $X_t$ the phenotypic lag of this population.

Here, $v>0$ is the speed of environmental change (in direction $\mathbf{e_1}$), 
$B_t$ is a standard $\mathcal{F}_t$ Brownian motion 
and $M$ is a Poisson Random Measure (PRaMe) on $\bR_+ \times \bR^d \times \bR_+ $,
also adapted to $\cF_t$,
with intensity:
\begin{align*}
\pi (ds,\, dw,\, du_f,\, du_g ) = ds \; \nu(dw) \; du_f\; du_g,
\end{align*}
where $\nu(dw)$ is a measure describing the distribution  of new mutations,
and:
\begin{align*}
\varphi_0 ( x,\, n,\, w,\, u_f,\, u_g) 
= \idc{ u_f\le f_0(n) } \times \idc{ u_g \le g(x,w )  }.  
\end{align*}

The independence between $M$ and $B$ is automatically deduced 
from the following proposition.
\begin{prop}
	\label{O_Qind}
	A Brownian Motion and a PRaMe 
	that are adapted to the same filtration 
	and such that their increase after time $t$ is independent from $\cF_t$
	are necessarily independent.
\end{prop}

\paragraph{Proof of Proposition \ref{O_Qind}.} 
Thanks to Theorem 2.1.8 of \cite{dT13},
if $X_1$, $X_2$ are 
additive functionals and semi-martingales with respect to a common filtration, 
both starting from zero,
and such that their quadratic covariations $[X_1,X_2]$ is a.s. zero,
then the random vector 
$(X_1(t) - X_1(s), X_2(t) - X_2(s))$
is independent of $\cF_s$, for every
$0 \le s \le t$. 
Moreover, the vector $(X_1, X_2)$ of additive processes is independent.

Note $B$ the Brownian Motion 
and $M$ the PRaMe on $\bR_+\times \cX$ 
For any test function $F:\cX \mapsto \bR$,
define 
$Z(t):= \int_{[0,t] \times \cX} F(x)\; M(ds, dx).$
Both $Z$ and $B$
are additive functionals and semi-martingales 
with respect to the filtration $\cF_t$, 
both starting from zero.
$Z$ being a jump process and $B$ continuous, 
their quadratic covariation equals a.s. 0.
Since it applies to any $F$,
exploiting Theorem 2.18 of \cite{dT13} implies that
$B$ and $M$ are independent.
\hfill $\square$
\\

In the model of the moving optimum originally considered in \cite{KH09}, 
$X = 0$ corresponds to the optimal state in terms 
of some reproductive value function $R(x)$,
for $x\in \bR$.
This function $R$ is also assumed to be symmetrical  
and decreasing with $|X|$.
Here we consider a possibly multidimensional state space for $X$
and will usually not require any assumption on the related function $g$.

$X$ is described as the phenotypic lag
because 
$X_t + v\, t\, \mathbf{e_1}$ is the character of the individuals at time $t$ in the population
while in this original model, 
the mobile optimum is located at trait $v\, t\, \mathbf{e_1}$.
These assumptions on the fitness landscape are natural,
and we abide by them in our simulations.
Nonetheless,
they are mainly assumed for simplicity
and we have chosen here to be as general as possible in the definition of $r$.
$X_t$ is thus a lag as compared to the trait $v\, t\, \mathbf{e_1}$
that is merely a reference value.

$g(X_t, w)$ 
is the  mutation kernel
that describes the rate of fixation
at which a mutant subpopulation of trait $X_t +v\,t\, \mathbf{e_1} + w$ 
invades a resident population 
of trait $X_t + v\,t\, \mathbf{e_1}$.
Although the rate at which the mutations occur in one individual 
can reasonably assumed to be symmetrical in $w$,
it is clearly not the case for $g$.
In a large population, 
the filtering of considering only fixing mutations
highly restricts the occurrence of strongly deleterious mutations,
strongly favors strongly advantageous mutations.
For mutations with little effects, 
there is only a slight bias.
To cover both of these situations, 
we consider in our analysis 
both the case where any mutation effect is permitted 
and the case where only advantageous ones are.
Although the latter case will raise more difficulty
in terms of accessibility of the domain, 
the core of the argument is quite the same
and the simulations seem to provide similar results in both cases.

The term $f(N_t)$ is introduced to model the fact 
that for a constant mutation rate by individual,
the mutation rate for the population is all the larger 
than the population size is large. 
$f(N_t):= N_t$ is the first reasonable choice,
but we may also be interested in introducing an effect of the population size 
in the fixation rate.

$N$ follows the equation of a Feller logistic diffusion 
where the growth rate $r$ at time $t$ only depends on $X_t$, 
while the strength of competition $c$ and the coefficient of diffusion $\sigma$ 
are kept constant.
Such a process is the most classical ones for the dynamics of a large population size
in a continuous space setting and such that explosion is prevented.
It is described in \cite{L05} (with fixed growth rate),
notably as a limit of some individual-based model.
$\sigma$ is related to the proximity between to uniformly sampled individuals
in terms of their filiation links:
$1/\sigma^2$ scales as the population size 
and is sometimes describes as the "effective population size".
\\

From a biological perspective, 
$X$ has no reason to explode.
Under our assumption $[H11]$ below,
such explosion is clearly prevented. 
Yet, we won't focus on conditions ensuring non-explosion for $X$.
Indeed, it would mean 
(by assumption $[H8]$ below)
that the growth rate becomes extremely negative.
It appears very natural 
to consider that it would lead
to the extinction of the population.
So, we define the extinction time as:
\begin{align*}
\ext:= \inf\{t \ge 0 \pv N_t = 0\}\wedge \Tsup{k\ge 1} T_X^k,
\where T_X^k:=  \inf\{t \ge 0 \pv \|X_t\| \ge  k\}.
\EQn{O_TNO}{\ext}
\end{align*}

Because it simplifies many of our calculations, 
in the following,
we will 
consider $Y_t:= \frac{2}{\sigma} \sqrt{N_t}$ 
rather than $N_t$. 

\begin{fact}
With the previous notations,
$(X, \,Y)$ satisfies the following SDE:
\begin{align*}
	(S)&\left\{
	\begin{aligned}
		X_t 
		&= x - v\,t\, \mathbf{e_1} 
		+ \int_{[0,t] \times \bR^d \times (\bR_+)^2 } w \; 
		\varphi \left( X_{s^-},\,Y_s,\,w,\,u_f,\, u_g \right)\; \Md,\\ 
		Y_t 
		&= y + \int_{0}^{t} \psi\lp X_s,\, Y_s\rp ds + B_t,
	\end{aligned}
	\right.
	\\&\text{where we define:}
	\hcm{0.5}
	\psi(x, y) =  - \frac{1}{2\, y} + \frac{r(x)\, y}{2} - \gamma\, y^3, \quad
	\text{ with }
	\gamma:= \frac{\gamma_0\, \sigma^2}{8},
	\\ &\hcm{3}
	\varphi(x,\, y,\,w,\,u_f,\, u_g) 
	:= \varphi_0  \lp x,\, \sigma^{2} y^{2}/4 ,\,
	w,\, u_f,\, u_g \rp.
	\\& 
	\text{ Thus with }f(y):= f_0  [ \sigma^{2} y^{2}/4 ],
	\quad
	\varphi ( x,\, y,\, w,\, u_f,\, u_g) 
	= \idc{ u_f\le f(y) } \times \idc{ u_g \le g(x,w )  }.  
\end{align*}
\end{fact}
An elementary application of the Ito formula proves this fact.
\\

The aim of the following theorems 
is to describe the law of the marginal 
of the process $(X, Y)$
at large  time $t$ conditionally 
upon the fact that the extinction has not occurred,
in short  the MCNE at time $t$. 
Considering the conditioning at the current time 
leads to considering properties of quasi-stationarity;
while a conditioning at a much more future time 
leads to a Markov process usually referred to as the Q-process,
in some sense the process conditioned on never going extinct.
The two aspects are clearly complementary 
and our approach will treat both in the same framework,
in the spirit initiated by \cite{ChQSD}.

\subsection{Elementary notations}
\label{O_sec_not}
In the following, the notation $k\ge 1$ is to be understood as $k\in \mathbb{N}$ while $\tp \ge 0$ --resp. $c>0$-- should be understood as $\tp \in \bR_+:= [0, \infty)$ --resp. 
$c\in \bR_+^* $ $:= (0, \infty)$. In this context (with $m\le n$),
we denote classical sets of integers by:
$\quad \mathbb{Z}_+:= \Lbr 0,1,2...\Rbr$,$\; \mathbb{N}:= \Lbr 1,2, 3...\Rbr$,
$\; [\![m, n ]\!]:= \Lbr m,\, m+1, ..., n-1,\, n\Rbr$,
where the notation $:=$ makes explicit that we define some notation by this equality.
For maxima and minima, we usually denote:
$s \vee t:= \max\{s, t\}$,\,  
$s \wedge t:= \min\{s, t\}.$
Accordingly, for a function $\varphi$, $\varphi^{\wedge}$ 
--resp. $\varphi^{\vee}$--  will be the notation
 for a lower-bound --resp. for an upper-bound-- of $\varphi$.
$C^0(X, Y)$ denotes the set of continuous functions from any $X$ to any $Y$.
$\cB(X)$ is the set of bounded functions from any $X$ to $\bR$.
$\cM(X)$ and $\cM_1(X)$  denote the sets of resp. positive measures
and probability measures on any state space $X$.
Numerical indices are rather indicated in superscript, 
while specifying notations are often in subscript.
By notation, $\{y\in \cY\pv A(y)\mVg B(y)\}$ denotes the set
of values $y$ of $\cY$ such that both $A(y)$ and $B(y)$ hold true.
Likewise, for two probabilistic conditions $A$ and $B$ on $\omega\in \Omega$, 
and a r.v. $X$,
we may use $\E(X\pv A\mVg B)$ instead 
of $\E(X\idg{\Gamma})$ 
where $\Gamma := \{\omega \in \Omega\pv A(\omega)\mVg B(\omega)\}.$

\section{Exponential convergence to the $QSD$}
\label{O_Res}
\setcounter{eq}{0}

\subsection{Hypothesis}
\label{O_hyp2}

We will consider two different sets of assumptions,
including or rejecting 
the possibility for deleterious mutations to invade the population.

First, the following set $(H)$ of assumptions 
can always be assumed throughout the paper, although some assumptions may be mentioned as not involved.
	\Hyp{HfPos}		 \Hyp{HgB}			\Hyp{HrN}	
	\Hyp{HdS}	\Hyp{Hac} 	
\begin{enumerate}
\item[\RHyp{HfPos}]
$f\in \mathcal{C}^0\big( \bR_+^*, \,\bR_+\big)$
is positive.

\item[\RHyp{HgB}]
$g\in \mathcal{C}^0\big( \bR^d \times \bR^d, \,\bR_+\big)$ 
and is bounded on any $K\times \bR^d$, 
where $K$ is a compact set of $\bR^d$.

\item[\RHyp{HrN}] 
$r$ is locally  Lipschitz-continuous on $\bR^d$
and $r(x)$ tends to $ -\infty$ as $\|x\|$ tends to $\infty$.

\item[\RHyp{HdS}] 
$\nu(\bR^d) < \infty$.
Moreover, there exist $\theta, \nu_\wedge > 0$
and $\dS \in (0, \theta)$
such that: 
$$\nu(dw) \ge \nu_\wedge \; \idg{B(\theta + \dS)\setminus B(\theta - \dS)}\; dw,$$
where $B(R)$, for $R>0$, denotes the open ball of radius $R$ centered at the origin.

\item[\RHyp{Hac}]
provided $d\ge 2$, 
$\nu(dw)<< dw$
 and the density  $g(x, w)\, \nu(w)$ 
 (for a jump from $x$ to $x+w$),
 of the jump size law w.r.t. Lebesgue's measure,
satisfies:
\begin{align*}
\frl{x_\vee>0}\quad
\sup \Lbr \dfrac{g(x, w)\, \nu(w)}{\int_{\bR^d} g(x, w')\, \nu(w')\, dw'}\pv
\|x\| \le x_\vee,\, w\in \bR^d \Rbr < \infty.
\end{align*}
\end{enumerate}

When we allow deleterious mutations to invade the population,
we actually mean that the rate of invasion is always positive, leading to the following assumption:
\begin{center}
\HgD \qquad
$g$ is positive.	
\end{center}

Otherwise, we consider the case 
where deleterious mutations are forbidden, 
in the sense the rate is zero for mutations 
that would induce an increase of $\|X\|$.
The invasion rate of advantageous mutations is however still assumed to be positive. This is stated in the following assumption \HgA\, as the alternative to \HgD.
\begin{align*}
\HgA \hcm{1}
	\text{ For any }x, w\in \bR^d,\qquad
	&\|x+w\| < \|x\| 
	\text{ implies }
	g(x, w) > 0
	\\
	&\|x+w\| \ge \|x\| 
	\text{ implies }
	g(x, w) = 0.
\end{align*}

\begin{rem}
$\star$ For $d=1$, 
no condition on the density of $g\ltm \nu$ as in \RHyp{Hac} 
is  required.
	
$\star$ 	It is quite natural to assume that $f(0) = 0$ 
	and  $f(y)$ tends to $\infty$ as $y $ tends to $\infty$, but we will not need those assumptions. 
	
$\star$ $1$ is the natural bound 
with the above-mentioned biological interpretation
of \RHyp{HgB}.
Yet an extension can be introduced
	where $g$ is not exactly the fixation probability, 
	cf. Corollary \ref{O_CorMj}.

	$\star$ Under \RHyp{HgB} and \RHyp{HdS}
	(since $\nu(\bR^d) < \infty$),
	over any finite time-interval, 
	only a finite number of mutations can occur. 
We also need lower-bounds on the probability 
	of specific events 
	which roughly prescribe the dynamics of $X$.
	This is where the lower-bound on the density of $\nu$ 
	is exploited as well as the positivity of $g$,
	deduced from either Assumption \HgD\, or \HgA.
	
$\star$ The fact that $r(x)$ tends to $-\infty$ 
as $x$ tends to $\infty$ 
makes it easier to prove 
that the process remains mostly confined,
say in the time-interval $[0, t]$
under the conditioning that $\{t<\ext\}$.
We would be able to state an explicit value $r_\wedge$ 
depending on the other parameters
so that the proof holds while assuming 
that the limsup of $r(x)$ is upper-bounded by $-r_\wedge$ 
when $\|x\|$ tends to infinity
(instead of being necessarily equal to $-\infty$).
\end{rem}

\subsection{Statement of the main theorems}
\label{O_Thm}

First, we need to ensure that the model specified by equation $(S)$
properly defines a unique solution, which is stated in the next Proposition.
\begin{prop}
\label{O_Ex}
Assume that \HDef\, holds.
Then, for any initial condition 
$(x, y) \in \bR^d\times \bR_+^*$,
there is a unique strong solution 
$(X_t, Y_t)_{t\ge 0}$
in the Skorokhod space
satisfying $(S)$
for any $t < \ext$,
and $X_t = Y_t = 0$ for $t\ge \ext$,
where $\qquad \ext:= \Tsup{n\ge 1} T_Y^n
\wedge \Tsup{n\ge 1} T_X^n\mVg $\\
$T_Y^n:= \inf\{t \ge 0, \; Y_t \le 1/n \}
\mVg \qquad
T_X^n:=  \inf\{t \ge 0, \; \|X_t\| \ge  n\}.$
\end{prop}

\begin{rem}
This proposition makes it possible to express $\ext$ 
as $\inf\{t \ge 0, \; Y_t = 0 \}$.
\end{rem}

We exploit the notion of uniform exponential quasi-stationary convergence
as previously introduced in \cite{AV_disc}.
\begin{defnt}
	\label{UEQS}
	For any linear and bounded semi-group $(P_t)_{t\ge 0}$
	acting on a Polish state space $\cZ$, 
	we say that $P$ displays a uniform exponential quasi-stationary convergence
	with characteristics
	$(\alpha, h, \lambda) \in \cM_1(\cZ)\ltm B(\cZ)\ltm \bR$
	if $\LAg \alpha \bv h\RAg = 1$
	and there exists $C, \gamma>0$ such that for any $t>0$
	and for any measure $\mu\in \cM(\cZ)$ with $\NTV{\mu}\le 1$:
	\begin{align*}
		\NTV{e^{\lambda t} \mu P_t(ds) - \LAg \mu\bv h\RAg \alpha(ds)}
		\le C e^{-\gamma t}.
		\EQn{PtCV}{}
	\end{align*}
\end{defnt}

\begin{rem}
	$\star$	As shown in Fact 2.2.2 of \cite{AV_disc}
	it implies that 
	for any $t>0,$ $\alpha P_t(ds) = e^{-\lambda t} \alpha(ds)$.
Any measure satisfying this property is called a quasi-stationary distribution.
	
	It is elementary that $h_t:x\mapsto e^{\lambda t} \LAg \delta_z P_t\bv \idg{} \RAg$ converge in the uniform norm to $h$. 
	We call $h$ the survival capacity because 
	$e^{\lambda t} \LAg \delta_x P_t\bv \idg{} \RAg
	= \PR_x(t<\ext)/\PR_\alpha(t<\ext)$
	enables to compare the likelihood of survival
with respect to the initial conditions.
	
	Since $h_{t+t'} = e^{\lambda t} P_t h_{t'}$, 
	one can then easily deduce that $e^{\lambda t} P_t h = h$.
	It is also obvious that $h$ is necessarily non-negative.

	$\star$	 By ``characteristics", we express that they are uniquely defined.	
\end{rem}

Our main theorem is stated as follows,
with $\cZ:= \bR^d\ltm \bR_+^*$:
\begin{theo}
\label{O_ECVGnl}
Assume that Assumption $(H)$ holds.
Suppose that either \HgD\, or \HgA\, holds.
If $d\ge2$, assume finally that \RHyp{Hac} holds.
Then, the  semi-group $P$ 
associated to the process $Z:= (X, Y)$
and extinction at time $\ext$
displays a uniform exponential quasi-stationary convergence
with some characteristics
$(\alpha, h, \lambda) \in \cM_1(\cZ)\ltm B(\cZ)\ltm \bR_+$.
Moreover, $h$ is positive.
\end{theo}

\begin{rem}
We refer to Corollary 2.2.1 in \cite{AV_disc},
for the implied convergence result of the renormalized semi-group
to $\alpha$. 
The fact that $h$ is positive implies
 that there is no other QSD in $\cM_1(\cZ)$.
\end{rem}

In \cite{AV_disc} 
is also provided an analysis of the so-called Q-process,
whose properties are as follow:
\begin{theo}
	\label{O_thQ}
Under the same assumptions as in Theorem \ref{O_ECVGnl},
 with $(\alpha, h, \lambda)$ 
the characteristics of exponential convergence of $P$, 
the following properties hold:

(i) \textbf{Existence of the $Q$-process:} \\
There exists a family $(\mathbb{Q}_{(x, y)})_{(x, y)\in \cZ}$ of probability
measures on $\Omega$ defined by:
\begin{align*}
	\limInf{\tp }
	\PR_{(x, y)}(\Lambda_\spr \bv \tp  < \ext) 
	= \mathbb{Q}_{(x, y)}(\Lambda_\spr),
	\EQn{Qxdef}{\mathbb{Q}_{(x, y)}}
\end{align*}
for all $\cF_\spr$-measurable set $\Lambda_\spr$. 
The process $(\Omega;(\cF_\tp )_{\tp \ge 0};
(X_\tp, Y_\tp)_{\tp \ge 0};(\mathbb{Q}_{(x, y)})_{(x, y)\in \cZ)}$
is an $\cZ$-valued homogeneous strong Markov process.

\textit{\textbf{(ii) Weighted exponential ergodicity of the Q-process:} \\
The measure $\quad \beta(dx, dy):= \heig(x, y)\, \alpha(dx, dy)$ 
is the unique invariant probability measure under $\mathbb{Q}$.
\\
Moreover, for any $\mu \in \cM_1(\cZ)$
satisfying $\LAg \mu\bv 1/\heig\RAg < \infty$ 
and $\tp \ge 0$:}
\begin{align*}
& \NTV{ \mathbb{Q}_{\mu} \lc\, (X_{\tp}, Y_\tp) \in (dx, dy)\rc 
	- \beta(dx, dy) }
\le C \,\|\mu -\LAg \mu\bv 1/h\RAg \, \beta\|_{1/\heig}\; e^{-\gamma \; \tp},
\EQn{D:ECvBeta}{\textup{ECv:}\,\beta}
\\ &\text{where }
\mathbb{Q}_\mu (dw):= \textstyle{\int_\cZ}\mu(dx, dy) \, \mathbb{Q}_{(x, y)}) (dw),
\quad \|\mu\|_{1/\heig}:= \|\dfrac{\mu(dx, dy)}{\heig(x, y)}\|_{TV}.
\end{align*}
\end{theo}

\begin{rem}
	$\bullet$ For the total variation norm, 
	considering $(X, Y)$ or $(X, N)$ is equivalent.
	
	$\bullet$ The constant $\LAg \mu\bv 1/h\RAg $ in \Req{D:ECvBeta}
	is optimal up to a factor 2, in the sense that for any $u>0$:
	$\|\mu - u\,\alpha\|_{1/h} 
	\ge \|\mu - \LAg\mu\bv 1/h \RAg \beta\|_{1/h} / 2$
	(cf Remark 2.2.5 of \cite{AV_disc}).
	
	$\bullet$	Since $r$ tends to  $-\infty$ as $\|x\|$ tends to  infinity,
	it is natural to assume that mutations leading $X$
	to be large have a very small probability of fixation.
	Notably,
	it means that we highly expect
	the upper-bound of $g$ in \RHyp{HgB}, uniform over $w$.
	
	$\bullet$ Under hypothesis \HgA,
	one may expect the real probability of fixation $g(x, w)$
	to be at most of order $O(\|w\|)$ for small values of $w$
	(and locally in $x$). 
	In such a case, we can allow $\nu$ to satisfy a smaller integrability condition
	than \RHyp{HdS} while forbidding observable accumulation of mutations.
\end{rem}

\begin{cor}
	\label{O_CorMj}
	Assume that \HDef\, and \HgA\, hold.
	Suppose that $\int_{\bR} (|w|\wedge 1)\; \nu(dw)  <\infty$
	while $\wtd{g}:(x, w) \mapsto g(x, w)/(|w|\wedge 1)$ 
	is bounded on any $K\ltm \bR^d$ for $K$ a compact set of $\bR^d$.
	If $d\ge2$, assume additionally that \RHyp{Hac} holds.
	Then, the conclusions of Theorem \ref{O_ECVGnl} 
	and Theorem \ref{O_thQ} hold true.
\end{cor}

\paragraph{Proof of Corollary \ref{O_CorMj}.}
%
$(X, Y)$ is solution of $(S)$ iff it is solution of:
\begin{equation*}
(S)\left\{
\begin{aligned}
X_t 
&= x - v\,t\, \mathbf{e_1}
+ \int_{[0,t] \times \bR^d \times \bR_+ } w \; 
\wtd{\varphi} \left( X_{s^-},\,Y_s,\,w,\,u_f,\, \wtd{u_g} \right)\; \wtd{M}(ds, dw, du_f, du_g),\\ 
Y_t 
&= y + \int_{0}^{t} \psi\lp X_s,\, Y_s\rp ds + B_t,
\end{aligned}
\right.
\end{equation*}
where $\wtd{M}$ is a PRaMe of intensity $ds \; \wtd{\nu}(dw) \; du_f\; d\wtd{u_g}$\mVg
\begin{align*}
&\wtd{\nu}(dw) 
:= \nu(dw)/(\|w\|\wedge 1)\mVg 
\quad
\wtd{\varphi}( x,\, y,\, w,\, u_f,\, \wtd{u_g}) 
= \varphi ( x,\, y,\, w,\, u_f,\, \wtd{u_g}\ltm(\|w\|\wedge 1))\mVg
\end{align*}
where $\wtd{\varphi}$ is defined as $\varphi$
with $g$ replaced  by $\wtd{g}$.

Thanks to the condition on $\nu$, \RHyp{HdS} holds with $\wtd{\nu}$
instead of $\nu$.
Thanks to the condition on $g$, \RHyp{HgB} still holds with $\wtd{g}$ instead of $g$.
Conditions \HgA\, and \RHyp{Hac} are equivalent for the systems $(g,\nu)$ and $(\wtd{g}, \wtd{\nu})$.
Consequently, if we prove Theorem \ref{O_ECVGnl} 
and Theorem \ref{O_thQ}
with \RHyp{HgB} and \RHyp{HdS}, 
the results follow under the assumptions of Corollary \ref{O_CorMj}.
\hfill $\square$

\subsection{Eco-evolutionary implications of these results}
\label{O_sec_ecoevo}

One of the major motivation 
for the present analysis 
is to make a distinction, as rigorous as possible,
between an environmental change to which the population can spontaneously adapt to
and a change that imposes too much a pressure. 
We recall that in \cite{NP17},
the authors obtain a clear and explicit threshold 
on the speed of this environmental change.
Namely, above this speed, 
the Markov process that they consider is transient, 
whereas it is recurrent below this critical speed.
Thus, it might seem a bit frustrating that such a distinction 
(depending on the speed value $v$)
cannot be observed in the previous theorems.
At least,
these results prove that the distinction 
is not based on the existence nor the uniqueness of the QSD,
and even not on the exponential convergence 
per se.

In fact, 
the reason why
this threshold is so distinct in \cite{NP17}
is that their model is based
on the following underlying assumption:
The poorer the current adaptation is, 
the more efficiently 
mutations are able to fix,
provided that they are then beneficial.
In our case,
a population that is too poorly adapted
is almost doomed to a rapid extinction, 
because the population size
cannot be maintained at large values.
Instead, 
long-term survival is triggered 
by dynamics that maintain 
the population adapted.
Looking back at the history of surviving populations,
it means that we are likely to observe that 
the process has mostly remained confined 
outside of deadly areas.

In order to establish this distinction
between environmental changes 
that are sustainable 
and those that endanger the population,
we need a criterion 
that quantifies the stability 
of such core regions.
Our results provide two exponential rates
whose comparison is enlightening:
if the extinction rate is of the same order 
as the convergence rate 
or larger, 
it means that the dynamics 
is strongly dependent
upon the initial condition.
If the convergence is much faster,
the dynamics shall rapidly 
become similar regardless the initial condition.
This is at least the case
for  initial conditions  that are not too risky
(i.e. where $\heig$ is not too small).
This criterion takes into account 
the intrinsic sustainability 
of the mechanisms involved 
in the adaptation to the current environmental change,
but does not involve the specific initial state of adaptation.

Looking at the simulation results,
the convergence in total variation
indeed appears to happen 
at some exponential rate,
provided that extinction 
does not abruptly wipe out
a large part of the distribution
 at a given time.
To obtain a generic estimate 
of the exponential rate 
at which the effect of the initial condition is lost,
the decay in total variation appears however
computationally expensive  and not very meaningful.
Although they are not as clearly justified,
it seems more practical to exploit
 the decay in time of the correlations of $X$ and/or $N$ 
starting from the QSD profile.
It does not seem very difficult to compare 
the extinction rate
from this estimate.
This is especially true 
in the case where $\cX$ is of dimension one, 
as one can directly estimate the dynamics of the density 
and thus the extinction rate.
Furthermore, it is 
quite reassuring 
to see that the inclusion or not
deleterious mutations 
(for which the invasion probability 
is expected be positive but very small)
is not crucial in the present proof.
We do not see much difference 
by looking at the simulations.

Much more can be said by looking at the simulation estimates
of the quasi-stationary (QSD),
the quasi-ergodic distribution (QED), and the survival capacity.
It is planned to detail these simulation results 
in a later article, 
but let us already give some  insights 
into the comparison between the QSDs and the QEDs
provided in Appendix B.
We see that although the QSDs look very different
at the three different values of mutation rates, 
the QEDs
are in fact very similar.
When extinction plays a notable role,
a tail appears on the QSD
 from the area of concentration 
of the QED
to an area where the population size is close to zero.
The shape of the tail 
and the fact that it does not appear 
on the QED nor for larger mutation rate
suggests that it corresponds 
in some sense to a path towards a rapid extinction.
These regions are clearly more unstable 
than the core areas 
where the Q-process remains confined.
This is probably due to this decline
in population size
when the level of maladaptation becomes more pronounced.
This confinement by the conditioning upon survival
only weakens in the recent past.
Conditions 
most likely leading to extinction 
are allowed
provided the delay is sufficiently large before extinction actually occurs.

\subsection{Quasi-ergodicity of related models}
\label{O_QErel}

The current paper completes the illustrations 
given in Subsection 4.2 of \cite{AV_QSD}
 and Sections 4-5 of \cite{AV_disc}.
If the model of the current paper 
was in fact the original motivation 
for the techniques presented in these two papers,
we can focus more closely 
on each of the difficulties
thanks to these various illustrations.
In any of them, 
the adaptation of the population to its environment
is described by some process $X$ solution to some SDE of the form:
\begin{align*}
X_t = x - \int_0^t V_s\, ds + \int_0^t \Sigma_s \cdot dB_s
+ \int_{[0,t] \times \bR^d \times \bR_+ } 
w \; \idc{u \le U_s(w)}\, M(ds, dw, du),
\end{align*}
where $B$ is a $\cF_t$-adapted Brownian Motion 
and $M$ a $\cF_t$-adapted PRaMe.
$V_s$ and $\Sigma_s$ a priori depend on $X_s$, 
$U_s$ on $X_{s-}$ 
and possibly on a coupled process $N_t$
 describing the population size.
Like the product $f(Y_t)\, g(X_{t-}, w)$ in equation $(S)$,
one specifies in $U_t(w)$ 
the rate at which a mutations of effect $w$ 
invades the population at time $s$.
 $V_s$ both relates to the speed of the environmental change
 and to the mean effects of the mutations 
 invading the population at time $s$
 in a limit of very frequent mutations of very small effects.
  $\Sigma_s$ then relates to the undirected fluctuations 
  both of the environment 
  and in the effects of this large number of small fixating mutations.
  
  We can relate 
  the current coupling of $X$ and $N$
  to an approximation given 
  by the autonomous dynamics of a process $Y$ similar to $X$.
  For the approximation to be as valid as possible,
  the law of $Y$ should be biased by some extinction rate 
  (depending at time $t$ on the value $Y_t$)
  and its jump rate should be adjusted.
  By these means, we would take an implicit account of 
  what would be the fluctuations of $N$ 
  if $X$ would be around the value of $Y_t$.
  This approximation is particularly reasonable
  when the characteristic fluctuations of $N$
  around its quasi-equilibrium
  are much quicker 
  than the effect of the growth rate changing over time
  with the adaptation.
  Its validity is less clear 
  when the extinction has a strong effect 
  on the establishment of the quasi-equilibrium.
  
The exponential quasi-stationary convergence 
  is treated in Subsection 4.2 of \cite{AV_QSD}
  for a coupling $(X, N)$ that behaves as an elliptic diffusion,
while  Sections 4  and 5 of \cite{AV_disc}
deal with some cases of a biased autonomous process $Y$
 that behaves as a piecewise-deterministic process.
For such a process with jumps, 
it is manageable yet technical to deal with restrictions 
on the allowed directions or sizes of jumps,
while imposing $V_t$ to stay at zero actually makes the proof harder than choosing $V_t := v\ltm t$.
  While the proofs of $(A1)$ and $(A3)$ highly depends 
  on such local properties of the dynamics,
  the ones of $(A2)$ for these semi-groups
  rely on a common intuition. 
Although we allow $X$ to live in an unbounded domain,
the maladaptation of the process 
when it is far from the optimal position
constraints $X$ to stay confined 
conditionally upon survival.
This effect of the maladaptation has been modeled
either directly on the growth rate of the coupled process $N_t$
or with some averaged description in terms of extinction rate.
Such confinement property for the coupled process
is in fact the main novelty of \cite{AV_QSD}
and notably illustrated in Subsection 4.2.4.
For simplicity, 
we have dealt  there with a locally elliptic process,
for which the Harnack inequality is known 
to greatly simplify the proof,
as observed previously for instance in \cite{ChpLyap}.
The proof of this confinement
is actually simpler
with $Y$ behaving as an autonomous process
under the pressure of  a death rate
 going to infinity outside of compact sets.
The proof in this case
is naturally adapted from the proof of $(A2)$ 
given in Subsection 4.1.2 of \cite{AV_QSD}.
\\

Assume for now that 
the fluctuations of $N$ 
are much quicker
than the change of the growth rate 
in the domain where the population is well-adapted.
Then we conjecture that 
considering the autonomous process $Y$
(including the bias by the extinction rate)
instead of the coupled process $(X, N)$
would produce very similar results:
the extinction rates
and the rates of stabilization to equilibrium
should be close between these models,
while the QSD profile of $X$ should be similar to the one of $Y$.

The drop in the quality of the approximation
when extinction has a crucial contribution
 must have a quite limited effect
for our concern, 
which is to compare the extinction rate 
to the rate of stabilization to equilibrium,
see Subsection \ref{O_sec_ecoevo}.
Indeed, 
as long as the extinction rate is not way larger
than the rate of stabilization to equilibrium,
such domains of maladaptation are strongly avoided
when looking in the past of surviving populations.
On the other hand, 
the population is almost doomed 
when it enters these domains,
so that we should be able to neglect 
the contribution to the extinction rate
of the dynamics of the process there.


%


\subsection{The mathematical perspective on quasi-stationarity}
\label{O_MQS}

The subject of quasi-stationarity is now quite vast 
and a considerable literature is dedicated to it, 
as suggested by the bibliography collected by Pollett
\cite{QSDbibli}. 
Some insights into the subject 
can be found in general surveys 
like \cite{coll}, \cite{DP13} 
or more specifically for population dynamics in
\cite{MV12}.
However, it appears that
that much remains to be done 
for the study of  strong Markov processes 
both on a continuous space 
and in a continuous time,
without any property of reversibility.
For general recent results,
besides \cite{AV_QSD} and \cite{AV_disc} that we exploit,
 we refer to \cite{CVly2}, 
\cite{BCGM19}, \cite{CG20}, \cite{FRS20} or \cite{GNW20}.
The difficulty is increased
when the process is discontinuous 
(because of the jumps in $X$)
and multidimensional,
since the property of reversibility
becomes all the more stringent 
and new difficulties arise
(cf e.g. Appendix A of  \cite{CCM17}).

Thus, ensuring the existence and uniqueness 
of the QSD 
is already some breakthrough,
and we are even able to ensure 
an exponential rate of convergence in total variation
to the QSD and similar results on the Q-process.
This model is in fact a very interesting illustration 
of the new technique which we exploit.
Notably, we see how conveniently 
our conditions are suited for exploiting the Girsanov transform
as a way to disentangle couplings
(here between $X$ and $N$, that are respectively the evolutionary component 
and the demographic one).

Our approach relies on the general result 
presented in \cite{AV_disc}, which, 
as a continuation of \cite{AV_QSD},
has been originally motivated 
by this problem.
In \cite{AV_QSD},
the generalization of Harris recurrence property 
at the core of the results of \cite{ChQSD}
is extended 
to deal with exponential convergence 
which are not uniform with respect to the initial condition.
The fine control over the MCNE
has opened the way 
for the approach developed in \cite{AV_disc}
to deal with continuous-time and continuous-space 
strong Markov processes with discontinuous trajectories.

%
After their seminal article \cite{ChQSD}, 
these same authors have obtained 
quite a number of extensions,
for instance 
with multidimensional diffusions \cite{CCV17},
 inhomogeneous in time processes\cite{InhomChp},
and various examples of processes in a countable space 
notably with the use of Lyapunov functions, 
cf. \cite{ChpLyap} or \cite{CVly2}.
Exploiting the result 
of \cite{CVly2},
it may be possible to ensure
the properties 
of exponential quasi-ergodicity
for such a discontinuous process as the one of this article,
keeping a certain dependence on the initial condition.
At least, 
the conditions they provide 
as well as the ones from \cite{BCGM19} 
are necessarily implied by our convergence result
(cf  Theorem 2.3 of \cite{CV20} or Theorem 1.1 in \cite{BCGM19}).
Yet, in the approach of \cite{CVly2}
for continuous-time 
and continuous space Markov process,
the rather abstract assumption $(F3)$
appears tightly bound to
the Harnack inequality.
The similar Assumption $(A4)$ in \cite{BCGM19}
is also left without further guidance,
while the assumption of a strong Feller property
 in \cite{FRS20} and \cite{GNW20} appears too restrictive.
For discontinuous processes,
these two properties generally do not hold true,
which is what motivated us to look 
for an alternative statement in \cite{AV_disc}.
This technique
is very efficient here.

This dependence on the initial condition
is biologically expected, 
although its crucial importance becomes apparent 
when the population is already highly susceptible 
to extinction.
For a broader comparison of this approach with the general literature,
we refer to the introductions of \cite{CVly2}, 
and the comparison with the literature provided in \cite{AV_QSD} and  \cite{AV_disc}.

\section{Proof of Proposition \ref{O_Ex}}
\label{O_exU}
\setcounter{eq}{0}

\textbf{Uniqueness:}\\
\noindent\textbf{Step 1: A priori upper-bound on the jump rate.}

Assume that 
 we have a solution $(X_t, Y_t)_{t\le T}$ to $(S)$ 
until some (stopping) time $T$ 
(i.e. for any $t<T$)
satisfying $T\le t_\vee\wedge T_Y^m\wedge T_X^n$
for some $t_\vee>0$, $m, n \ge 1$
(see Equation \Req{O_TNO}).
We know from \RHyp{HrN}\, 
that the growth rate of the population 
remains necessarily upper-bounded 
by some $r^{\vee}>0$ until $T$. 
Thus, we deduce a stochastic upper-bound $(Y^{\vee}_t)_{t\ge 0}$ 
on $Y$:
\begin{align*}
Y^{\vee}_t 
&= y + \int_{0}^{t} \psi^{\vee}( Y_s)\, ds + B_t
\where
\hcm{0.5}
\psi^{\vee}(y) =  - \frac{1}{2\, y} + \frac{r^{\vee}\, y}{2} - \gamma\, y^3,
\EQn{O_Y+}{}
\end{align*} 
which is thus independent of $M$.
Since $\psi^{\vee}(y) \le r^{\vee}\, y/2$,
it is classical
 that $Y^{\vee}$ --and a fortiori $Y$--
  cannot explode before $T$,
  see for instance Lemma 3.3 in \cite{BM15}
  or \cite{L05} where such a process is described in detail. 

Under $\RHyp{HgB}$,
the jump rate of $X$ is uniformly bounded 
until $T$ by:
$$\textstyle
\nu(\bR^d)\ltm
\sup\Lbr  g(x', w)
\pv x'\in \obB(0, n), w\in \bR^d \Rbr
\ltm \sup\{f(y')\pv y'\le \sup\!_{s\le t_\vee} Y^{\vee}_s\} < \infty\;  a.s.$$

\noindent\textbf{Step 2: Identification until $T$.}

In any case, 
this means that the behavior of $X$ until $T$ 
is determined 
by the value of $M$ on a (random)  domain 
associated to an $a.s.$ finite intensity.
Thus, we need a priori to consider only 
a finite number $K$ of "potential" jump, 
that we can describe as the points
$(T\iJP^i, W^i, U_f^i, U_g^i)_{i\le K}$
in the increasing order of the times $T\iJP^i$.

From the a priori estimates, 
we know that for any $t< T\iJP^1\wedge T$:
$X_t = x - v\, t.$
By the improper notation $t<T\iJP^1\wedge T$,
we mean $t< T\iJP^1$ if $K\ge 1$ 
(since $T\iJP^1 < T$ by construction)
and $t< T$ if  $K = 0$, i.e.
when  there is no potential jump before $T$. 
We then consider the solution $\hat{Y}$ of:
$$\hat{Y}_t 
= y + \textstyle{\int_{0}^{t}} \psi\lp x-v\,s \mVg \hat{Y}_s\rp ds + B_t.$$

It is not difficult to adjust the proof of \cite{YW71}
to this time-inhomogeneous setting,
with $\RHyp{HrN}$,
so as to prove the existence and uniqueness 
of such a solution 
until any stopping time $T \le \hat{\tau}_\partial$,
where $\hat{\tau}_\partial 
:= \inf\{t\ge 0, \hat{Y}_t = 0\}$. 
Besides, $\hat{Y}$ 
is independent of $M$
and must coincide with $Y$ 
until $T_J^1\wedge T$.
Since $T \le T_Y^m$, 
the event $\{\hat{\ext} < T_J^1\wedge T\}$
is necessarily empty.
If there is no potential jump before $T$,
i.e. $K = 0$,
we have identified $(X_t, Y_t)$ for $t\le T$
as $X_t = x-v\, t$, $Y_t = \hat{Y_t}$.
Else, at time $T\iJP^1$,
we check whether $U_f^1\le f(\hat{Y}(T\iJP^1))$ 
and $U_g^1\le g(x - v\, T\iJP^1\mVg W^1)$.
If it holds, 
necessarily $X(T\iJP^1) = x - v\, T\iJP^1 + W^1$,
else $X(T\iJP^1) = x - v\, T\iJP^1$.
Doing the same inductively for the following time-intervals 
$[T\iJP^k, T\iJP^{k+1}]$,
we identify the solution $(X, Y)$ until $T$.
\\

\noindent\textbf{Step 3: Uniqueness of the global solution.}

Now, 
consider two solutions $(X, Y)$ and $(X', Y')$ 
of $(S)$ 
defined up to 
respectively $\ext$ and $\ext'$ as in Proposition \ref{O_Ex}
with in addition $X_t = Y_t = 0$ for $t\ge \ext$,
and $X'_t = Y'_t = 0$ for $t\ge \ext'$.

On the event 
$\{\sup_{m} T_y^m = \ext\wedge \ext'\}$,
we deduce by continuity of $Y'$
that $T_y^m = T_y^{'m}$
so that $\ext = \ext'$. 
On the event 
$\{\sup_{n} T_X^n = \ext \le \ext'< \infty\}$,
for any $n$ and $t_\vee>0$ 
there exists $m\ge 1$ and $n'\ge n$
such that $T_X^n\wedge t_\vee <  T_Y^m\wedge T_Y^{'m}$
 and $\|X(T_X^n\wedge t_\vee) \|
 \vee \|X'(T_X^n\wedge t_\vee) \| < n' < \infty$.
Thanks to Step 2, $(X, Y)$ and $(X', Y')$ 
must coincide until 
$T = (t_\vee+1)\wedge 
T_Y^{m}\wedge T_Y^{'m}
\wedge T_X^{n'}\wedge T_X^{'n'}$,
where the previous definitions 
ensure $T_X^n\wedge t_\vee < T$
(with the fact that $X$ and $X'$ are right-continuous). 
This proves that $T_X^n\wedge t_\vee  = T_X^{'n}\wedge t_\vee$,
and with $t_\vee, n\ifty$ that $\ext' = \ext$.

By symmetry between the two solutions,
we have a.s. $\ext = \ext',\quad$
$\frl{t<\ext} X_t = X'_t$
 and
$\frlq{t\ge \ext} X_t = X'_t = 0$.
It concludes the proof of the uniqueness.
\\

\noindent\textbf{Existence.}
We see that the identification obtained for the uniqueness
clearly defines the solution $(X, Y)$ until 
some $T = T(t_\vee, n)$ such that either $T = t_\vee$ or $Y_T = 0$ or $\|X_T\| \ge n$.
Thanks to the uniqueness property and the a priori estimates, 
this solution coincide 
with the ones for larger values of $t_\vee$ and $n$.
Thus, it indeed produces a solution up to time $\ext$.
\hfill $\square$

\section{Main properties leading to the proof of Theorem \ref{O_ECVGnl}}
\label{O_str}
\setcounter{eq}{0}

\subsection{General criteria for the proof of exponential quasi-stationary convergence}

The proof of Theorem \ref{O_ECVGnl}
relies on the set of Assumptions  $\mathbf{(AF)}$ 
presented in \cite{AV_disc},
and that we recall next.
$\mathbf{(AF)}$ is stated in the general context of a càdlàg process $Z$ on a Polish state $\cZ$,
with extinction at time still denoted $\ext$.
The notations are changed from \cite{AV_disc}
to avoid confusion with the current ones,
$Z$ corresponding now to the couple $(X, Y)$.
We introduce the following notations 
for the exit and first entry times of any set $\cD$:
\begin{align*}
	T_{\cD}:= \inf\Lbr  t \ge 0 \pv Z_t \notin \cD \Rbr
\mVg\quad
\tau_\cD:= \inf\Lbr t \ge 0\pv Z_t \in \cD \Rbr.
\EQn{TcD}{}
\end{align*}
The assumptions involved in $\mathbf{(AF)}$ are the following ones.

\begin{itemize}

	\item[$(A0_S)$]
	There exists a sequence $(\cD_\ell)_{\ell\ge 1}$ of closed subsets of $\cZ$ such that for any $\ell\ge 1$,\\ $ \cD_\ell \subset
	int(\cD_{\ell+1}) \quad$
	(with $int(\cD)$ the interior of $\cD$).
	
	\item[$(A1)$ ]
	There exists a probability measure $\zeta \in  \cM_1(\cZ)$ such that, 
	for any $\ell\ge 1$, 
	there exists $L>\ell$ and  $c, t>0$ such that:
	\begin{align*}
	&  \frl{z \in \cD_{\ell}}
	\hspace{.5cm}
	\PR_z \lc {Z}_{t}\in dx\pv
	t < \ext \wedge T_{\cD_L} \rc 
	\ge c\; \zeta(dz).	
	\end{align*}
	
	\item[$(A2)$] 
	$\Tsup{z\in \cZ} \;\E_{z} \lp 
	\exp\lc\rho\, (\ext\wedge \tau_E) \rc \rp < \infty.$
	
	\item[$(A3_F)$]
	for any $\eps\in (0,\, 1)$,
	there exist $\tZa, c >0$ 
	such that for any $z \in E$
	there exists a stopping time 
	$\Uza$ 
	such that:
	\begin{align*}
	&  \Lbr\ext \wedge \tZa \le \Uza \Rbr
	= \Lbr\Uza = \infty\Rbr
	\quad  \text{ and }\quad
	\PR_{z} (\Uza = \infty, \,  \tZa< \ext) 
	\le \eps\, \exp(-\rho\, \tZa),
	\EQn{D:FL}{}
	\end{align*} 
	
	while for some stopping time~
	$\UCa$:
	\begin{align*}
	&\PR_{z} \big(Z(\Uza) \in dz' \pv \Uza < \ext \big) 
	\le c \,\PR_{\zeta} \big(Z(\UCa) \in dz'
	\pv \UCa < \ext\big).
	\EQn{D:Abs}{}
	\end{align*} 
	We further require 
	that there exists a stopping time $U_{A}^\infty$ 
	extending $\Uza$
	in the following sense:
	
	$\star$   $U_{A}^\infty:= \Uza$ 
	on the event  $\Lbr \ext\wedge \Uza < \tau_E^1\Rbr$, 
	where $\tau_E^1:= \inf\{s\ge \tZa: 
	Z_s \in E\}$. 
	
	$\star$ On the event 
	$\Lbr \tau_E^1 \le \ext \wedge  \Uza \Rbr$
	and conditionally on $\cF_{\tau_E^1}$,
	the law of $U_{A}^\infty - \tau_E^1$ coincides
	with the one of $\wtd{U}_{A}^\infty$ 
	for a realization $\wtd{Z}$ 
	of the Markov process $(Z_t, t\ge 0)$
	with initial condition $\wtd{Z}_0:= Z(\tau_E^1)$ 
	and independent of $Z$ conditionally on $Z(\tau_E^1)$.
\end{itemize}
	$\rho$ as stated in Assumptions $(A2)$ and $(A3_F)$ 
is required to be strictly larger
than the following \textbf{"survival estimate"}: 
\begin{align*}
\rho_S
:= \sup\big\{\gamma \ge 0\pv 
\sup_{L\ge 1} \inf_{t>0} \;
e^{\gamma t}\,\PR_\zeta(t < \ext\wedge T_{\cD_L}) 
= 0
\big\}\vee 0.
\end{align*}
We are now in position to state $\mathbf{(AF)}$:

"$(A1)$ holds 
for some  \mbox{$\zeta \in \cM_1(\cZ)$}
and a sequence $(\cD_\ell)_\ell$ satisfying $(A0_S)$.
Moreover, there exist $\rho  > \rho_S$ 
and a closed set $E$ such that $E\subset \cD_\ell$ for some $\ell\ge 1$
and such that
$(A2)$  and $(A3_F)$  hold."
\\

As stated next
by gathering the results of Theorems 2.2, 2.3 and Corollary 2.2.3 of \cite{AV_disc},
 $\mathbf{(AF)}$ implies the convergence results that we aim,
noting that the sequence $(\cD_\ell)_\ell$ will cover the whole space.
Some additional properties of approximations
are also obtained, 
where the process is localized to large $\cD_L$ by extinction.

\begin{theo}
	\label{thECV}
Provided that $\mathbf{(AF)}$ holds,
the semi-group $P_t$ associated to the process $Z$ with extinction at time $\ext$ 
displays a uniform exponential quasi-stationary convergence
with some characte\-ristics
$(\alpha, h, \lambda) \in \cM_1(\cZ)\ltm B(\cZ)\ltm \bR$.

Moreover,
consider for any $L\ge 1$ the semi-group $P^L$
for which  the definition of $\ext$  	is replaced 
by $\ext^L:= \ext \wedge T_{\cD_L}$.
Then, 
for any $L\ge 1$ sufficiently large,
$P^L$ displays a uniform exponential quasi-stationary convergence
with some characteristics
$(\alpha^L, h^L, \lambda_L) \in \cM_1(\cD_L)\ltm B(\cD_L)\ltm \bR_+$.
The associated versions of \Req{PtCV}
hold true with constants that can be chosen uniformly in $L$.
As $L$ tends to  infinity, 
$\lambda_L$ converges to $\lambda$  and  $\alpha^L,h^L$ converge to $\alpha,h$  in total variation and pointwise respectively. 

If in addition, $\medcup_{\ell \ge 1} \cD_\ell = \cZ$,
then $h$ is positive
and the results of Theorem \ref{O_thQ} on the Q-process hold also true.
\end{theo}

\begin{rem}
Under $\mathbf{(AF)}$, the Q-process can generally be defined on $\mathcal{H}:=\{z\in \cZ;\; h(z)>0\}$ and the fact that $h$ is positive is not required 
or may be proven as a second step.
The proof of Theorem \ref{thECV}
however provides lower-bound of $h$ on any $\cD_\ell$,
so that $\cZ = \medcup_{\ell \ge 1} \cD_\ell$ is a practical assumption for the proof that $h$ is positive.
\end{rem}

\begin{rem}
The assumption $(A3_F)$ appears certainly technical 
and its usage is the main focus of \cite{AV_disc}.
It is referred to as the ``Absorption with failure'' property
and makes it possible to upper-bound 
the asymptotic survival probability
 from initial condition $z$
 as compared to the one from initial condition $\zeta$.
 To this purpose, 
 a coupling is introduced 
 where \Req{D:Abs} makes it possible to ``absorb'' most trajectories.
 Since failures where $U_A=\infty$ 
 while $\tZa < \ext$ are allowed,
 this step is to be iterated 
 and the probability of such failure 
is to be controlled through \Req{D:FL}.
\end{rem}

For the proof of Theorem \ref{O_ECVGnl},
the sequence $(\cD_\ell)_{\ell\ge 1}$  is defined as follows:
\begin{align*}
\cD_\ell:= \bar{B}(0, \ell) \times [1/\ell, \ell],
\EQn{O_Dn}{\cD_\ell}
\end{align*}
where $\bar{B}(0, \ell)$ denotes the closed ball of radius $\ell$
for the Euclidian norm.

Forbidding deleterious mutations 
in the case of unidimensional $\cX$
will make our proof a bit more complicated.
This case is thus treated later on.
The expression "with deleterious mutation"
 will be used a bit abusively 
to discuss the model under $\HgD$.
On the other hand, 
the expression "with only advantageous mutation" will refer to the case 
where $\HgA$ holds.

These criteria are proved to hold true
under the assumptions of Theorem \ref{O_ECVGnl}
in the following Theorems \ref{O_Mix_DF}-6.
We see in Subsection \ref{O_sec_mixP}
how these theorems together with Theorem \ref{thECV}
 imply Theorem \ref{O_ECVGnl}.
In the next subsections,
we then prove
  Theorems \ref{O_Mix_DF}-6.
  By mentioning first the mixing estimate, 
  we wish to highlight the constraint on the reachable domain
  under hypothesis \HgA.
 The order of the proofs is different and done for the clarity
 of their presentation.
The mixing estimates are handled similarly
 under the different sets of assumptions
and directly exploited 
in the proofs of the absorption estimates.
The escape estimates are very close 
to the ones of previously considered models,
so more easily dealt with.

\subsection{The whole space is accessible: with deleterious mutations or $d\ge 2$}

\subsubsection{Mixing property and accessibility}
\label{O_sec_mixP}

With deleterious mutations,
 the whole space becomes accessible. 
 It is in fact also the case with only advantageous mutations,
 provided $d\ge 2$:
\begin{theo}
\label{O_Mix_DF}
Suppose \HDef.
For $d = 1$, assume \HgD.
For $d\ge2$, assume either \HgD or \HgA.
Then, for any $\ell_I, \ell_M\ge 1$, 
   there exists $L>\ell_I\vee \ell_M$ and  $\cp, \tp>0$ such that:
\begin{align*}
&
  \frlq{(x_I,\, y_I) \in \cD_{\ell_I}}
 \PR_{(x_I,\, y_I)} \lc {(X,\, Y)}_{\tp}\in (dx,\,dy)\pv
  \tp < \ext \wedge T_{\cD_L} \rc 
  	\ge \cp\, \idg{\cD_{\ell_M}}(x,\, y) \,dx\, dy.
  	\EQn{O_MixDF}{}
  \end{align*}
\end{theo}

\begin{rem}
$\bullet$ \Req{TcD} is exploited when defining
$ T_{\cD_L} 
	:= \inf\Lbr \tp\ge 0\pv (X,\, Y)_\tp \notin \cD_{L} \Rbr$.
	
$\bullet$ Theorem \ref{O_Mix_DF} implies in particular that
the density w.r.t. Lebesgue's measure 
of any QSD 
is uniformly lower-bounded on any $\cD_\ell$.

$\bullet$ In the case where \HgD\, holds, 
$L:= \ell_I\vee \ell_M + \theta$ can be chosen.
The choice of $t$ cannot generally be made arbitrary,
at least for $d=1$,
since the lower-bound of the density of jump sizes
is only valid for jumps of size close to $\theta$.
Under \HgA with $d\ge 2$, 
the constraint that jumps 
must be advantageous 
makes the convenient choice of $L$ less clear.
\end{rem}

\subsubsection{Escape from the Transitory domain}

\begin{theo}
	\label{O_eT_D}
	Assume Assumptions $(H)$. Then,
	for any $\rho > 0$, there exists $\ell\iET\ge 1$
	such that $(A2)$ holds 	
	 with $E:= \cD_{\ell\iET}$.
\end{theo}

\begin{rem}
Heuristically, it means that the killing rate can be made arbitrarily large by adding killing 
when hitting some compact $\cD_\ell$ 
that sufficiently covers $\cZ = \bR\ltm \bR_+^*$.
\end{rem}
    
\subsubsection{Absorption with failures} 
We need some reference set
on which our reference measure has positive density. 
With the constants $\theta$ and $\eta$ involved in $[H4]$ let:
\begin{align*}
\varDelta 
:= \bar{B}(-\theta\, \mathbf{e_1}\mVg \dS) 
\times [1/2, 2].
\EQn{O_DAl}{\varDelta}
\end{align*}
This choice (rather arbitrary),
is made in such a way 
that the uniform distribution on $\varDelta$
can be taken as the lower-bound 
in the conclusions of Theorems \ref{O_Mix_AF}
and \ref{O_Mix_DF}.

Including deleterious mutations or with $d\ge 2$, 
we will exploit the following theorem 
for sets $E$ of the form  $E:= \cD_{\ell\iET}$,
where  $\ell\iET$ is determined thanks to Theorem \ref{O_eT_D}.
But the theorem holds generally for any  closed subsets $E$ 
of $\bR^d\times \bR_+^*$
for which there exists $\ell\ge 1$ such that $E\subset  \cD_\ell$,
property that we briefly denote as $E \in \mathbf{D}$. 

\begin{theo}
\label{O_AF_D}
Suppose \HDef.
For $d = 1$, assume \HgD.
For $d\ge2$, assume either \HgD or \HgA.
Then, 
for any $\rho > 0$, $\fl\in (0,\, 1)$  and $E \in \mathbf{D}$, 
there exist $\tZa, \cp >0$ 
which satisfy the following property
for any $(x, y) \in E$ 
and $(x_{\alc}, y_{\alc})\in \varDelta$.
There exists a stopping time $\Uza$ 
such that:
\begin{align*}
&  \Lbr\ext \wedge \tZa \le \Uza \Rbr
= \Lbr\Uza = \infty\Rbr
\quad  \text{ and }\quad
   \PR_{(x, y)} (\Uza = \infty, \,  \tZa< \ext) 
     \le \fl\, \exp(-\rho\, \tZa),
      \end{align*} 
and an additional stopping time~$\UCa$ such that:
\begin{align*}
&\PR_{(x, y)} \big[(X(\Uza),Y(\Uza))  \in (dx', dy')
 \pv \Uza < \ext \big]
 \\&\hcm{0.5}
       \le \cp \,\PR_{(x_{\alc}, y_{\alc})} \big[
       (X(\UCa),Y(\UCa))  \in (dx', dy')
       \pv \UCa < \ext\big].
\EQn{AbsZ}{}
   \end{align*} 
Moreover, there exists a stopping time $\Uza^{\infty}$
satisfying the following properties:
    
$\bullet$   $\Uza^{\infty}:= \Uza$ 
   on the event  $\Lbr \ext\wedge \Uza < \tau_{E}^1\Rbr$, 
   where $\tau_{E}^1:= \inf\{s\ge \tZa: 
   (X_s, Y_s) \in E\}$. 
     
$\bullet$ On the event $\Lbr \tau_{E}^1 < \ext\Rbr \cap 
\Lbr \Uza = \infty \Rbr$,
and conditionally on $\cF_{\tau_{E}^1}$,
the law of $\Uza^{\infty} - \tau_{E}^1$ coincides
with the one of $\wtd{U}_{A}^{\infty}$ 
for the solution $(\wtd{X}, \wtd{Y})$ of:
\begin{align}
&\left\{
\begin{aligned}
\wtd{X}_r 
&= X(\tau_{E}^1) - v\,r\, \mathbf{e_1} 
+ \int_{[0,r] \times \bR^d \times (\bR_+)^2 } w \; 
\varphi \left( \wtd{X}_{s^-},\,\wtd{X}_s,\,w,\,u_f,\, u_g \right)\; \MdT\\ 
\wtd{Y}_r
&= Y(\tau_{E}^1) + \int_{0}^{r} \psi\lp \wtd{X}_s,\, \wtd{Y}_s\rp ds + \wtd{B}_r,
\end{aligned}
\right.\EQn{O_St}{}
\end{align}
where $r\ge 0$, $\wtd{M}$ and $\wtd{B}$ are independent copies 
of respectively $M$ and $B$.
   \end{theo}

\subsubsection{Proof of Theorem \ref{O_ECVGnl} as a consequence of Theorems \ref{O_Mix_DF}-3}
\label{O_sec_TECV}


\begin{itemize}
	\item First, it is clear that the sequence 
	$(\cD_\ell)_\ell$ satisfies both $(A0_S)$ and
	$\medcup_{\ell \ge 1} \cD_\ell = \cZ$.

\item $(A1)$ 
holds true thanks to Theorem \ref{O_Mix_DF}, 
  where $\alc$ is the uniform distribution over $\varDelta$ --cf \Req{O_DAl}.
  
\item Theorem \ref{O_eT_D} implies $(A2)$ for any $\rho$,
and we also require that $\rho$ is chosen such that:
$$\rho > \rho_S:= \sup\big\{\gamma \ge 0\pv 
\sup_{L\ge 1} \inf_{t>0} \;
e^{\gamma t}\,\PR_\alc(t < \ext\wedge T_{\cD_L}) 
= 0
\big\}\vee 0.$$
Thanks to Lemma 3.0.2 in \cite{AV_QSD} and $(A1)$, 
we know that 
$\rho_S$ is upper-bounded 
by some value $\wtd{\rho}_S$.
In order to satisfy  $\rho>\rho_S$,
we set $\rho:= 2 \wtd{\rho}_S$. 
Thanks to Theorem \ref{O_eT_D}, 
we deduce $E = \cD_{\ell\iET}$ such that
assumption $(A2)$ holds for this value of $\rho$.

\item Finally, Theorem \ref{O_AF_D} 
implies that assumption $(A3_F)$ holds true,
for $E$ and $\rho$.
In the adaptation of \Req{AbsZ}
where $(x_{\alc}, y_{\alc})$ is replaced by $\zeta$,
$V$ is specified by the initial condition 
$(x_{\alc}, y_{\alc})\in \varDelta$
chosen uniformly according to $\zeta$.
\end{itemize}

This concludes the proof of the assumption $(AF)$ with 	$\medcup_{\ell \ge 1} \cD_\ell = \cZ$.
Exploiting Theorem \ref{thECV},
it implies Theorems \ref{O_ECVGnl} and \ref{O_thQ}
in the case where, besides \HDef, either \HgD\, holds
or $d\ge 2$ and \HgA\,  holds.
\hfill $\square$

%
%

\subsection{No deleterious mutations in the uni-dimensional case}

\subsubsection{Mixing property and accessibility}
When only advantageous mutations are allowed and $d =1$, 
as soon as the size of jumps is bounded, 
the process can't access some portion of space 
(there is a limit in the $X$ direction). 
We could prove that the limit is related to the quantity:
$\quad
L_A:= \sup\Lbr M;\; \nu [2\,M, +\infty) >0 \Rbr \in (\theta/2,\, \infty].
$\\
The accessible domains with maximal extension would then be rather of the form:
$
[-\ell\mVg L_A - 1/\ell] \times [1/\ell, \ell]$, for some $\ell\ge 1.$
To simplify the proof, the limit $L_A$ will however not 
appear in the next statements.
We just wanted to point out this potential constraint on the visited domain.
In fact, the $X$ component 
is assumed to be negative
in the following definition of the accessibility domains:
\begin{align*}
\Delta\iET := \{[-L, 0] \times [1/\ell, \ell];\;  L, \ell\ge 1\}.
\EQn{O_Dnac}{\Delta\iET}
\end{align*}

\begin{theo}
\label{O_Mix_AF}
Assume $d=1$, 
\HDef\ and \HgA. 
Then, for any $\ell_I\ge 1$ and $E\in \Delta\iET$, 
   there exists $\Lp>\ell_I$ and  $\cp, \tp>0$ such that:
\begin{align*}
&
  \frlq{(x_I,\, y_I) \in \cD_{\ell_I}}
 \PR_{(x_I,\, y_I)} \lc {(X_{\tp},\, Y_{\tp})}\in (dx,\,dy)\pv
  \tp < \ext \wedge T_{\cD_{\Lp}} \rc 
  	\ge \cp\, \idg{E}(x,\, y) \,dx\, dy.
  	\EQn{O_MixAF}{}
  \end{align*}
\end{theo}

\begin{rem}
Theorem \ref{O_Mix_AF} implies that
the density w.r.t. Lebesgue's measure of  any QSD 
is uniformly lower-bounded on any $E$
of the form given by \Req{O_Dnac}.
\end{rem}

\subsubsection{Escape from the Transitory domain}
\label{O_sec_escPD}

\begin{theo}
	\label{O_eT_A}
	Assume $d=1$, Assumptions $(H)$ and \HgA. 
	Then, for any $\rho > 0$, 
	there exists $E \in \Delta\iET$ such that $(A2)$ holds.
\end{theo}

\begin{rem}
	Heuristically, it means that the asymptotic killing rate can be made arbitrarily large by adding killing 
	when hitting some compact $E$ 
	that sufficiently covers $\bR_-\ltm \bR_+^*$.
\end{rem}

\subsubsection{Absorption with failures}
\label{O_sec_AFpD}
%
 
\begin{theo}
\label{O_AF_A}
Suppose \HDef\, and \HgA.
Then, 
for any $\rho > 0$, $\fl\in (0,\, 1)$  and $E \in \Delta\iET$, 
there exist $\tZa, c>0$ which satisfy the same property as in 
Theorem \ref{O_AF_D}.
\end{theo}

\begin{rem}
The definition of $\Delta$ is chosen to apply for both theorems
\end{rem}

\subsubsection{Proof of Theorem \ref{O_ECVGnl} as a consequence of Theorems \ref{O_Mix_AF}-6}
\label{O_sec_TECVD}

The argument being very similar to the one for the case $d\ge 2$ or with \HgD, we go briefly through it.

\begin{itemize}
\item $(A1)$ holds thanks to Theorem \ref{O_Mix_AF},
with again the choice of $\alc$ uniform on $\Delta$.

\item Thanks to Theorem \ref{O_eT_A},
and similarly as in the proof exploiting Theorem \ref{O_eT_D} in Subsection \ref{O_sec_TECV},
we deduce that there exists $E \in \Delta\iET$ 
such that $(A2)$ holds with some value $\rho > \rho_S$.

\item Finally, $(A3_F)$ holds  
for these choices of $\rho$ and $E$,
thanks to Theorem \ref{O_AF_D}.
\end{itemize}
This concludes the proof of the assumption $(AF)$ with 	$\medcup_{\ell \ge 1} \cD_\ell = \cZ$.
Exploiting Theorem \ref{thECV},
it implies Theorems \ref{O_ECVGnl} and \ref{O_thQ}
in the case where $d=1$, Assumptions $(H)$ and \HgA\, hold.
\hfill $\square$

\subsection{Structure of the proof}
To allow for fruitful comparison, 
the proofs are gathered according 
to the properties resp. $(A1)$, $(A2)$ and $(A3_F)$ they ensure.
We first prove Theorems \ref{O_eT_D} and \ref{O_eT_A}
in Section \ref{O_trans} since they are the simplest
and the closest to the proofs in \cite{AV_disc}
and the remaining theorems are more closely related.
We then prove Theorems \ref{O_Mix_DF} and \ref{O_Mix_AF}
in Section \ref{O_sec_mix},
and finally Theorems \ref{O_AF_D} and \ref{O_AF_A}
in Section \ref{O_sec_AF}.

\section{Escape from the transitory domain}
\label{O_trans}
\setcounter{eq}{0}

The most straightforward way 
to prove exponential integrability of first hitting times 
is certainly via Lyapunov methods. 
Yet, we highly doubt that this can be achieved 
as easily as we present next
given the interplay between the different domains 
on which the escape is to be justified.

\subsection
{With deleterious mutations or $d\ge 2$}
\label{O_sec_etD}

 Theorem \ref{O_eT_D}
 is a direct consequence of Proposition 4.2.2 in  \cite{AV_QSD}
and the proof is thus omitted. 
The process mainly considered in  \cite{AV_QSD}
is similar to this one 
in that there is also a coupling between a population size process $N$
and an adaptation process $X$.
Both population size processes are defined in the same way
in their relation to the process $X$ as 
$$N_t 
= n + \int_{0}^{t}  \left( r(X_s) \; N_s - \gamma_0 \ltm (N_s)^2 \right) ds 
+ \sigma \int_{0}^{t} \sqrt{N_s} \; dB_s.$$

Contrary to the current model,
$X$ is not evolving in \cite{AV_QSD} 
as piecewise deterministic with jumps,
but as a diffusion process.
If it changes significantly 
the proof of the other assumptions (A1) and (A3),
 this proof of $(A2)$ actually does not depend at all 
 on the dynamics of $X$,
 as expressed in Proposition 4.2.2 of \cite{AV_QSD}.
 The proof 
 developed in the next subsection 
 is an extension of this one
 and illustrates the technique.
 \\
 
 Let us give a few hints of how it works.
The proof relies
on uniform couplings
which ensure 
that with a probability close to 1,
the population size experience drastic decrease
sufficiently quickly,
be it when it starts at a very large value,
when the adaptation is very poor (large $\|X\|$)
or when the population is close to extinction.
In addition, 
we simply need to prove that
the probability of large increase 
is also very exceptional.

In practice,
 we distinguish between 3 different sets of initial conditions
depending on which of the above situation is to be considered,
like the sets $\cT_\infty$ , $\cT^X_\infty$  and $\cT_0$
from Figure \ref{O_eTA}.
The above-mentioned estimations
provide relations between the 3 exponential moments of return
starting from the different sets of initial conditions.
The decrease estimate proves 
that, prescribing a fixed time interval,
 the process exits during this time interval
 the set of conditions he starts in
with a probability sufficiently close to 1.
The increase estimate makes it possible to control 
the probability of trajectories rapidly navigating
 between the different sets of conditions.

 \subsection
 {Without deleterious mutations, $d=1$}
  In this section, we prove Theorem \ref{O_eT_A}, i.e.:
  
\textsl{ Suppose that $d=1$, \HDef\, and \HgA hold. 
 Then:}
 \begin{align*}
 &\frl{\rho >0}\quad
 \Ex{E\in \Delta\iET} \qquad
 \underset{(x, y)\in \bR\times \bR_+}{\sup} \;\E_{(x, \, y)} \lp 
       	\exp\lc\rho\, (\tau_{E}\wedge \ext)\rc \rp < \infty.
 \end{align*}
  
\begin{figure}[h]
\begin{center}
\includegraphics[width = 15cm, height = 9cm]{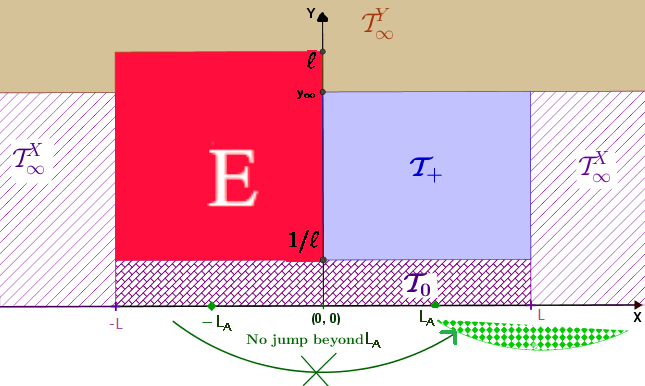}
\caption{subdomains for A2 }
\label{O_eTA}
\end{center}
\hrule
\vspace{0.2cm}
Until the process reaches $E$ or extinction, 
it is likely to escape any region either from below 
or from the side into $\cT_+$,
the reversed transitions being unlikely.
As long as $X_t > 0$, $\|X_t\|$ must decrease
(see Fact 5.2.5 in Subsection 5.2.4).
Once the process has escaped $\{x \ge L_A\}$,
there is no way (by allowed jumps and $v$)
that it reaches it afterwards.
\hrule
\end{figure}
   
\subsubsection{Decomposition of the transitory domain}  
The proof is very similar 
to the one of Subsection 4.2.4 of \cite{AV_disc} 
except that, due to Theorem \ref{O_AF_A}, 
the domain $E$ cannot be chosen as large.
We thus need to consider another subdomain of $\cT$, 
that will be treated specifically thanks to $\HgA$.

The complementary $\cT$ of $E$ is then made up of 4 subdomains:
 "$y=\infty$", "$y=0$", "$x>0$", and "$\|x\| = \infty"$, according to the figure \ref{O_eTA}.
Thus, we define:
 \begin{itemize}
	\item $\cT^Y_\infty:= 
	\Lbr (-\infty,\, -L) \cup(0, \infty)  \Rbr
	\times (y_{\infty}, \infty)$
	$\bigcup \; [-\ell, 0]\times [\ell, \infty)\hfill $("$y= \infty$"),
	\item $\cT_0:= (-L,\, L) \times [0, 1/\ell]\hfill$ ("$y = 0$"),
	\item $\cT_+:= (0, L)\times (1/\ell, y_{\infty}]\hfill$ ("$x > 0$"),
	 \item $\cT^X_\infty:= \Lbr \bR \setminus (-L,\, L) \Rbr \times (1/\ell, y_{\infty}]\hfill$
	 	 ("$|x| =\infty$").
\end{itemize}

With some threshold $t_\vee$ 
(meant to ensure finiteness
but whose effect shall vanish as it tends to $\infty$),
let us first introduce the exponential moments of each area 
(remember that $\tau_{E}$ 
is the hitting time of $E$):
  \begin{itemize}
  \begin{minipage}[t]{0.5\linewidth}
 	\item $\cE^Y_\infty 
 	:= 	\sup_{(x,\,y) \in \cT^Y_\infty}
		 	\E_{(x, y)}[\exp(\rho \; V_{E})]$,  	
 	\item  $\cE_0 
	:= \sup_{(x,\,y) \in \cT_0}
			\E_{(x, y)}[\exp(\rho \; V_{E})]$,
  \end{minipage}
  \begin{minipage}[t]{0.5\linewidth}
	 \item $\cE^X_\infty 
	:= \sup_{(x,\,y) \in \cT^X_\infty}
			\E_{(x, y)}[\exp(\rho \; V_{E})]$,
 	\item  $\cE_X
	:= \sup_{(x,\,y) \in \cT_+}
			\E_{(x, y)}[\exp(\rho \; V_{E})]$,
  \end{minipage}
 \end{itemize}
where $V_{E}:= \tau_{E} \wedge \ext\wedge t_\vee.$
  Implicitly, $\cE^Y_\infty,\, \cE^X_\infty,\, \cE_X $ and $\cE_0$ are functions of 
  $\rho,\, L,\, \ell, \,y_{\infty}$ that need to be specified.

\subsubsection{A set of inequalities}
Like in Subsection 4.2.4 in \cite{AV_QSD},
we first state some inequalities between these quantities, 
summarized in Propositions \ref{O_cEinfA}, \ref{O_cEXinfA}, \ref{O_cEx} and  \ref{O_cE0A}  that follow.
Thanks to these inequalities, 
we prove in Subsection \ref{O_CombA}
that those quantities are bounded.
This will end the proof of Theorem \ref{O_eT_A}. 

\begin{prop}
\label{O_cEinfA}
Suppose \HDef. Then, 
given any $\rho > 0$, there exist $y_{\infty} > 0$
and $C_{\infty}^Y \ge 1$
such that for any $\ell> y_{\infty}$ and $L>0$:
\begin{align*}
\cE^Y_\infty \le C_{\infty}^Y \, \lp 1+ \cE^X_\infty  + \cE_X \rp.
\EQn{O_EYi}{\cE^Y_\infty}
\end{align*}
\end{prop}

\begin{prop}
\label{O_cEXinfA}
Suppose  \HDef and \RHyp{HrN}. Then,  
given any $\rho>0$, 
there exists $C^X_\infty \ge 1$ 
which satisfies the following property
for any $\eps^X,\, y_{\infty} > 0$.
There exists $L>0$ and $\ell^X> y_{\infty}$ 
such that for any $\ell \ge \ell^X$:
\begin{align*}
\cE^X_\infty 
\le C^X_\infty \,\lp 1 + \cE_0 + \cE_X\rp
+ \eps^X \,\cE^Y_\infty.
\EQn{O_EXi}{\cE^X_\infty}
\end{align*}
\end{prop}

\begin{prop}
\label{O_cEx}
Suppose \HDef and \HgA. Then, 
given any $\rho,\, L>0$,
there exists $C_X \ge 1$ 
which satisfies the following property
 for any $\eps^{+},\,  y_{\infty} > 0$.
For any$\ell$ sufficiently large ($\ell \ge \ell^+> y_{\infty}$):
\begin{align*}
\cE_X \le C_X \, \lp 1+ \cE_0 \rp + \eps^{+} \, \cE^Y_\infty.
\EQn{O_EXL}{}
\end{align*}
\end{prop}

\begin{prop}
\label{O_cE0A}
Suppose  \HDef. Then, 
given any $\rho,\, \eps^0,\, y_\infty>0$,
there exists $C_0 \ge 1$
which satisfies the following property
for any $L$ and $\ell$ sufficiently large ($\ell \ge \ell^0> y_{\infty}$):
\begin{align*}
\cE_0 \le C_0  + \eps^0 \, \lp \cE^Y_\infty +\cE^X_\infty + \cE_X\rp.
\EQn{O_EY}{}
\end{align*}
\end{prop}

The proofs of Proposition  
\ref{O_cEinfA}, \ref{O_cEXinfA} and \ref{O_cE0A} 
can be taken mutatis mutandis
from the ones of Propositions 
respectively 4.2.1, 4.2.2 and 4.2.3
of \cite{AV_QSD}. 
The idea behind them is exactly the same 
as in the explanation provided for Subsection \ref{O_sec_etD},
made more precise with the statement of the propositions.
The justification relies on the estimate of what has happened 
in a time-interval $[0, t_0]$ for some arbitrary $t_0>0$:
by a proper definition of the domain boundaries,
observing a transition to some domain 
with lower-population sizes 
(or to extinction)
is justified to be a very likely event,
 with a probability 
 larger than $1-\exp(-\rho t_0)$.
Transitions to domain with a larger population size
shall be handled as very exceptional,
where an additional threshold  
is involved to specify transitions 
that are considered 
to happen during the time-interval $(0, t_0)$.
The only difference in the proofs 
is that transitions 
into $\cT_+$ are distinguished in the current paper,
which makes appear the term $\cE_X$ 
with factors resp. $C_{\infty}^Y$, $C^X_\infty$ and $\eps^0$.
The reason for $\eps^0$ to be as small as needed 
is that spontaneous extinction 
during some finite time-interval $[0, t_0]$
can be made as likely as needed,
while potential transitions (including those towards $\cT_+$)
are only considered at the end of the time-interval 
$[0,t_0]$.

We prove first how to deduce Theorem \ref{O_eT_A},
which naturally generalizes  
the similar argument in \cite{AV_QSD}.
Then, we will prove Proposition \ref{O_cEx},
that shall provide the main intuition 
for the proofs of the other propositions.
The core idea is that jumps are not allowed here 
to increase the maladaptation of the process.
Thus, the worst-case scenario 
for the exit time of $\cT_+$
is  that the process gets simply drifted 
by the environmental change at speed $v$
until $X$ gets negative.


\subsubsection{Proof that Propositions \ref{O_cEinfA}-4
	imply Theorem \ref{O_eT_A}}
\label{O_CombA}
We first prove that 
the inequalities \Req{O_EXi}, \Req{O_EXL} and \Req{O_EY}
given by Propositions \ref{O_cEinfA}-4
imply an upper-bound on
$\cE^Y_\infty \wedge \cE^X_\infty \wedge \cE_X\wedge \cE_0$
for sufficiently small 
$\eps^X$, $\eps^+$ and $\eps^0$.\\
Assuming first that $\eps^X \le  (2\,C_{\infty}^Y )^{-1}$,
we have:
\begin{align*}
	&\cE^X_\infty
	\le C^X_\infty \, \lp 3 + 3\, \cE_X + 2\, \cE_0\rp,
	\hcm{1}
	\cE^Y_\infty
	\le C_{\infty}^Y \,C^X_\infty \, 
	\lp 4 + 4\, \cE_X + 2\, \cE_0\rp.
\end{align*}
Assuming further that $\eps^+ \le (8\, C_{\infty}^Y \,C^X_\infty)^{-1}$:
\begin{align*}
	&\cE_X
	\le C_X\, \lp 2 + 3\, \cE_0\rp,
	\hcm{1}
	\cE^X_\infty
	\le C^X_\infty \,C_X \,
	\lp 9 + 11\, \cE_0\rp,
	\hcm{1}
	\cE^Y_\infty
	\le C_{\infty}^Y \,C^X_\infty \, 
	\lp 12 + 14\, \cE_0\rp.
\end{align*}
Assuming further that $\eps^0 \le (60\, C_{\infty}^Y \,C^X_\infty\, C_X)^{-1} $  (and exploiting $2\times[14+11+3]\le 60$):
\begin{align*}
	%
	&\cE_0
	\le 50\,C_0,
	\hcm{0.5}
	\cE_X
	\le 152\, C_X\,C_0,
	\hcm{1}
	\cE^X_\infty
	\le 559\, C^X_\infty \,C_X \,C_0,
	\hcm{1}
	\cE^Y_\infty
	\le 712\,C_{\infty}^Y \,C^X_\infty \, C_0,
	\\& 
	\text{In particular }
	\underset{(x, y)\in \bR\times \bR_+}{\sup} \;\E_{(x, \, y)} \lp 
	\exp\lc\rho\, (\tau_{E}\wedge \ext)\rc \rp 
	= \cE^Y_\infty \wedge \cE^X_\infty \wedge \cE_X\wedge \cE_0
	< \infty.
\end{align*}

Let us now specify the choice 
of the various parameters involved.
For any given $\rho$, we obtain from Proposition \ref{O_cEinfA}
the constant $y_{\infty}$, and $C_{\infty}^Y $ which gives us a value $\eps^X:= (2\,C_{\infty}^Y )^{-1}$.
We then deduce, thanks to Proposition \ref{O_cEXinfA}, some value for $C^X_\infty$, $\ell^X$ and $L$.
We can then be fix $\eps^+:= (8\, C_{\infty}^Y \,C^X_\infty)^{-1}$, 
and deduce, according to Proposition \ref{O_cEx}, 
some value $C_X$ and $\ell^+ > 0$. 
Now we fix $\eps^0:= (60\, C_{\infty}^Y \,C^X_\infty\, C_X)^{-1}$ and choose, according to Proposition \ref{O_cE0A}, 
some value $C_0$ and $\ell^0 > 0$. 
To make the inequalities 
\Req{O_EXi}, \Req{O_EXL} and \Req{O_EY}
hold, 
we can just take $\ell:= \ell^X \vee \ell^+ \vee \ell^0$.
With the calculations above, we then conclude Theorem \ref{O_eT_A}.
\hfill $\square$

\subsubsection{Proof of Proposition \ref{O_cEx}: 
	phenotypic lag pushed towards the negatives, }
\label{O_sec_cEx}

Since the norm of $X$ decreases at rate at least $v$ 
as long as the process stays
in $\widetilde{\cT}_+:= [0, L]\times \bR_+^*$, 
we know that the process cannot stay in this area
during a time-interval larger than $t_\vee:= \frac{L}{v}$.
This effect will give us the bound $C_X:= \exp\lp \rho\, L / v\rp$.

Moreover, 
we need to ensure 
that the transitions from $\cE_X$ to $\cE^Y_\infty$
are exceptional enough.
This is done exactly as for Proposition 4.2.2 in \cite{AV_QSD},
by taking $\ell^+$ sufficiently larger than $y_\infty$.
The event of having the process reaching
$\ell^+$ in the time-interval $[0, t_\vee]$
is then exceptional enough.

More precisely, 
given $L$ and $\ell > y_\infty\ge 1$ 
and initial condition $(x, y)\in \cT_+$, let:
\begin{align*}
	C_X:= \exp\lp \dfrac{\rho\, L}{v}\rp, \qquad
	\Tex:= \inf\Lbr t\ge 0\pv X_t \le 0 \Rbr \wedge V_{E}  \quad 
	\EQn{O_ToPos}{\Tex}
\end{align*}

\begin{fact}
	\label{APET}
	Assume that \HDef\, and \HgA\, hold. \\
	Then, for any initial condition $(x, y)\in \cT_+$,
	$(X, Y)_{\Tex} \notin \cT^X_\infty$ a.s. and:
\begin{equation*}
	\frl{t< \Tex}\qquad
	X_t \le x - v\, t \le L - v\, t
	\hcm{1}\text{so that }
	\Tex \le t_\vee:= L / v.
\end{equation*}
\end{fact}
Thanks to Assumption \RHyp{HdS},
an immediate induction 
on  the number of jumps previous to $\Tex\wedge t$
proves that the jumps of $X$ can only make its value decrease
(because it is positive 
while the absolute value must necessarily decrease).
It proves Fact \ref{APET}.
Thanks to it:
\begin{align*}
	&\E_{(x, y)}[\exp(\rho V_{E} )]
	= \E_{(x, y)}\Big[
	\exp(\Tex )\pv
	\Tex =  V_{E} \Big]
	+ \cE_0\; \E_{(x, y)}\Big[
	\exp(\Tex )\pv
	(X, Y)_{\Tex} \in \cT_0 \Big]
	\\&\hspace{3.5cm}
	+ \cE^Y_{\infty}\; \E_{(x, y)}\Big[
	\exp(\Tex )\pv
	(X, Y)_{\Tex} \in \cT^Y_\infty \Big]
	\\& \hspace{2cm}
	\le C_X\, \lp 1 + \cE_0\rp
	+ C_X\,  \cE^Y_{\infty}\;
	\PR_{y_\infty}(T_{\uparrow} \le t_\vee)	
	\\
	&\where 
	T_{\uparrow}:= \inf\Lbr t\ge 0
	\pv Y^{\uparrow}_t \ge \ell \Rbr,\quad
	\text{and } Y^{\uparrow} \text{is the solution of: } 
	\\&  Y^{\uparrow}_t
	:= y_{\infty} +  \int_{0}^{t} \psi_\vee \lp Y^{\uparrow}_s \rp \;ds  
	+ B_t 
	\qquad\text{ (again }
	\psi_\vee(y):= - \frac{1}{2\, y} +\frac{r_\vee\, y}{2} - \gamma\, y^3).
	\EQn{O_e:Yp}{Y^{\uparrow}}
\end{align*}
We conclude the proof of Proposition \ref{O_cEx} by noticing that:
$\quad  \PR_{y_\infty}(T_{\uparrow} \le t_\vee)	\cvifty{\ell} 0.$

\hfill $\square$

\section{Mixing properties and accessibility}
\label{O_sec_mix}
\setcounter{eq}{0}
Before we turn to the proofs 
of Theorems \ref{O_Mix_DF} and \ref{O_Mix_AF},
we describe the common elementary properties 
upon which they rely
in the three following subsections.
The first one gives the trick to disentangle 
the behavior of the processes $X$ and $N$
up to a factor on the densities.
Subsection \ref{O_Girsanov}
deals with the mixing property for the $Y$ process.
These results are exploited in Subsection \ref{sec_Xmix}
to obtain the elementary mixing properties
that allow to deduce $(A2)$.
The three next subsections starting from \ref{O_stabX}
deal respectively with the proofs of Theorem \ref{O_Mix_DF}
first under Assumption \HgD,
then under under Assumption \HgA\, and $d\ge 2$
and finally with the proof of  Theorem \ref{O_Mix_AF}.

\subsection*
{General mixing properties}

\subsection
{Construction of the change of probability under \RHyp{HdS}}

The idea of this subsection 
is to prove that we can think of $Y$ 
as a Brownian Motion 
up to some stopping time which will bound $\Uza$. 
If we get a lower bound 
for the probability of events in this simpler setup, 
this will prove that we also get a lower bound in the general setup.

\paragraph{The limits of our control}
\label{O_sec_Glim}

Let $t_G,\, x_\vee >0,$ $0< y_\wedge< y_\vee,$ $N_J \ge 1.$
Our aim is to simplify the law of $(Y_t)_{t\in [0, t_G]}$
as long as 
$Y$ stays in $[y_\wedge, y_\vee]$,
$\|X\|$ stays in $\obB(0, x_\vee)$, 
and at most $N_J$ jumps have occurred.
Thus, let: 
\begin{align*}
	&
	T_X:= \inf \Lbr t\ge 0\pv \|X_t\| \ge x_\vee \Rbr,\quad
	T_Y:= \inf \Lbr t\ge 0\pv Y_t \notin [y_\wedge, y_\vee] \Rbr.
	\EQn{O_TXL}{T_X, T_Y}
	\\& 
	g_{\vee}:= \sup\Lbr g(x, w)
	\pv \|x\|\le x_\vee, w\in \bR^d \Rbr,
	\quad
	f_{\vee} 
	:=\sup\Lbr f(y)\pv y\in [y_\wedge,\, y_\vee] \Rbr
	\EQn{O_fG}{f_G^\vee}
	\\& 
	\cJ:= \Lbr (w, u_g, u_f)\in \bR^d \times [0, f_{\vee}]\times [0, g_{\vee}]\Rbr,\quad
	\\&\hcm{2}
	\text{so that }
	\nu\otimes du_g\otimes du_f(\cJ) = \nu(\bR^d)\;g_{\vee}\, f_{\vee} <\infty.
\end{align*}

Our Girsanov's transform alters the law of $Y$
until the stopping time:
\begin{align*}
	&
	T_G:= t_G\wedge T_X \wedge T_Y \wedge U_{N_J},
	\EQn{O_TG}{T_G}
	\\&\where 
	U_{N_J}:= \inf\Lbr\, t\pv
	M([0, t] \times \cJ) \ge N_J + 1\; \Rbr.
	\EQn{O_UjNJ}{U_{N_J}}
\end{align*}
Note that the $(N_J+1)$-th jump of $X$ 
will then necessarily occur after $T_G$.

\paragraph{The change of probability}
\textcolor{white}{:}

We define:
\begin{align*}
	L_t:= - \int_{0}^{t\wedge T_G} \psi(X_s, Y_s) dB_s,
	\quad \AND D_{t}:= \exp\lc L_{t} -  \LAg L\RAg_t /2 \rc,
	\EQn{O_Lt}{L_t}
\end{align*}
the exponential local martingale associated with $(L_t)$.

\begin{theo}
	\label{O_G}
	Suppose \HDef. 
	Then, for any $t_G,\, x_\vee > 0,$ and $y_\vee > y_\wedge >0$, 
	there exists $C_G > c_G>0$
	such that a.s.
	and for any $t>0$,  $c_G \le D_{t} \le C_G$.
	In particular, 
	$D_t$ is a uniformly integrable martingale
	and $\beta_t = B_t - \LAg B, L\RAg_t$
	is a Brownian Motion under:
	$\qquad \PR^G_{(x, y)} 
	:= D_{\infty} \cdot \PR_{(x, y)}$.
	Moreover:
	\begin{align*}
		\frl{(x,\, y)\in \bR^d \times \bR_+}\hspace{1cm}   
		c_G \; \PR^G_{(x,\, y)} \le \PR_{(x,\, y)} \le C_G \; \PR^G_{(x,\, y)},
	\end{align*}
\end{theo}

On the event $\{t\le T_G\}$,
$Y_t = y + \beta_t$, i.e. 
$Y$ has the law of a Brownian Motion under $\PR^G_{(x, y)}$ 
up to time $T_G$.
This means 
that we can have bounds 
of the probability for events 
involving $Y$ as in our model 
by considering $Y$ as a simple Brownian Motion.
Meanwhile, 
the independence between its variations as a Brownian 
and the Poisson Process still hold 
due to Proposition \ref{O_Qind}.

\subsubsection
{Proof of Theorem \ref{O_G}}

The proof is achieved by ensuring uniform upper-bounds 
of $L_t$ and  $\LAg L\RAg_t$,
which corresponds to $L_\infty$ and  $\LAg L\RAg_\infty$
for $t_G$ replaced by $t\wedge t_G$.

\subsubsection*{Proof in the case where $r$ is $C^1$}
\begin{align*}
	\Let
	&\NGOM{r}:= \sup\Lbr\, |r(x)| \pv x\in \obB(0, x_\vee) \Rbr,
	\EQn{O_NLr}{\NGOM{r}}
	\\&
	\NGOM{r'}:= \sup\Lbr\, |r'(x)| \pv x\in \obB(0, x_\vee) \Rbr.
	\EQn{O_NLrp}{\NGOM{r'}}
\end{align*}
With $\psi_G^\vee$ an upper-bound of $\psi$ 
on $\obB(0,\, x_\vee) 
\times [y_\wedge,\, y_\vee]$
(deduced from \RHyp{HrN} ) 
and recalling that $(X, Y)$ belongs to this subset 
until $T_G$ (see \Req{O_TG}):
\begin{align*}
	\LAg L\RAg_{\infty}
	=  \int_{0}^{T_G} \psi(X_s, Y_s)^2 ds 
	\le t_G \times (\psi_G^\vee)^2. 
	\EQn{O_varL}{\LAg L\RAg}
\end{align*}

In the following,
we look for bounds on
$\int_{0}^{T_G} \psi(X_s, Y_s) dY_s$, 
noting that:
\begin{align*}
	&L_{T_G} + \int_{0}^{ T_G} \psi(X_s, Y_s) dY_s
	= \int_{0}^{T_G} \psi(X_s, Y_s)^2 ds
	\in [0, t_G \times (\psi_G^\vee)^2].
\end{align*}
\begin{align*}	
	&\int_{0}^{T_G} \psi(X_s, Y_s) dY_s 
	= \int_{0}^{T_G} \lp- \dfrac{1}{2 Y_s} + \frac{r(X_s)\; Y_s}{2} - \gamma\; (Y_s)^3 \rp dY_s.
	\EQn{O_PsiDY}{\psi\,dY}
\end{align*}
Now, thanks to Itô's formula:
\begin{align*}
	& \ln(Y_{T_G}) = \ln(y) 
	+ \int_{0}^{T_G} \dfrac{1}{ Y_s}  dY_s  
	- \frac{1}{2} 	\int_{0}^{T_G} \dfrac{1}{ (Y_s)^2} ds\\
	&
	\text{thus   }
	\left|\int_{0}^{T_G} \dfrac{1}{ Y_s}  dY_s \right|
	\le 2\;( |\ln(y_\wedge)|\vee|\ln(y_\vee)|) + \dfrac{t_G}{2\, (y_\wedge)^2} <\infty.
	\EQn{O_II}{I1}
\end{align*}
\begin{align*}
	&(Y_{T_G})^4 = y^4
	+ 4\, \int_{0}^{T_G} (Y_s)^3 dY_s
	+ 6\,  \int_{0}^{T_G} (Y_s)^2 ds\\
	&
	\text{thus   }\left|\int_{0}^{T_G} (Y_s)^3 dY_s \right|
	\le  (y_\vee)^4/4 + 3\, t_G\, (y_\vee)^2/2
	<\infty.
	\EQn{O_ID}{I2}
\end{align*}
\begin{align*}
	&r(X_{T_G-})\, (Y_{T_G})^2
	= r(x)\, y^2 
	+ 2\, \int_{0}^{T_G} r(X_s)\, Y_s\, dY_s
	+  \int_{0}^{T_G} r(X_s)\; ds 
	- v \int_{0}^{T_G} r'(X_s)\, (Y_s)^2 \, ds\\
	&\hspace{1 cm}
	+ \int_{[0,T_G) \times \bR^d \times \bR_+ } 
	\blp\;r(X_{s-} + w) - r(X_{s-})\;\brp \times (Y_s)^2
	\\&\hcm{5}
	\times \idc{ u_f\le f(Y_s) } \, \idc{ u_g \le g(X_{s^-},\,w)  }
	\Md.
	\EQn{O_RYY}{}
\end{align*}
Since $\forall s\le T_G,\,  Y_s \in [y_\wedge, y_\vee]$, 
we get from \RHyp{HgB} and \Req{O_fG}:
\begin{align*}
	\frl{s\le T_G}
	\frl{ w \in \bR^d} \quad 
	g(X_{s-}, w) \le g_{\vee},\quad f(Y_s) \le f_{\vee},\quad
	\AND T_G \le U_{N_J}.
\end{align*}
Since moreover $T_G \le T_X$:
\begin{align*}
	&\left|\int\!
	\underset{[0,T_G) \times \bR^d \times \bR_+ \hcm{2.5}}
	{\quad\blp\;r(X_{s-} + w) - r(X_{s-})\;\brp} \;(Y_s)^2\;
	\idc{ u_f\le f(Y_s) } \, \idc{ u_g \le g(X_{s^-},\,w)  }
	\Md\right|
	\\&\hcm{4}
	\le 2\, N_J\, \NGOM{r}\, (y_\vee)^2,
\end{align*}
so that \Req{O_RYY} leads to:
\begin{align*}
	2\; \left|\int_{0}^{T_G} r(X_s)\, Y_s\, dY_s \right|
	\le \blp 2\,(N_J +1)\, \NGOM{r} + \NGOM{r'}\; v\, t_G\brp \times (y_\vee)^2 
	+ \NGOM{r} \; t_G< \infty.
	\EQn{O_rYdY}{I3}
\end{align*}

Inequalities \Req{O_II}, \Req{O_ID}, \Req{O_rYdY}
combined with \Req{O_varL} 
conclude that $L_\infty$ and $\LAg L\RAg_\infty$ 
are uniformly bounded.
This proves the existence of $0<c_G<C_G$ such that a.s.
$c_G \le D_\infty \le C_G$.

This statement is a priori 
adapted for $t_G$ replaced by $t\wedge t_G$,
yet these bounds are actually the largest for $t=t_G$.
So it entails that $c_G \le D_t \le C_G$ hold uniformly in $t$.
The rest of the proof is only classical application
of Girsanov's transform theory.

\subsubsection*{Extension to the case 	where $r$ is only Lipschitz-continuous}
$\Req{O_II}$ and $\Req{O_ID}$ are still true, 
so we show that we can find the same bound on 
$\left|\int_{0}^{T_G} r(X_s)\, Y_s\, dY_s \right|$ 
where we replace $\NGOM{r'}$ 
by the Lipschitz-constant $\|r\|_{Lip}^G$ of $r$ on $\obB(0, x_\vee)$, 
by approximating $r$ by $C^1$ functions 
that are $\|r\|_{Lip}^G$-Lipschitz continuous.

\begin{lem}
	\label{APG}
	Suppose $r$ is Lipschitz continuous on $\obB(0,\, x_\vee)$
	for some $x_\vee>0$.
	Then there exists 
$r_n\in C^1\blp \obB(0,\, x_\vee), \bR\brp , n\ge 1$ such that:
$$\NGOM{r_n - r} \cvifty{n} 0 \quad \AND
\frl{n\ge 1} \NGOM{r'_n} \le \|r\|_{Lip}^G.$$
\end{lem}

\paragraph{Proof of Lemma \ref{APG}}

\textcolor{white}{:}

We begin by extending $r$ on $\bR^d$ 
with $r_G(x):= r\circ\Pi_G(x)$, 
where $\Pi_G$ is the projection on $\obB(0, x_\vee)$
(it is well-known that $r$ can be extended on 
$\obB(0, x_\vee)$ with the same Lischitz constant). 
Note that this extension $r_G$ is still $\|r\|_{Lip}^G$-Lipschitz continuous.
If we define now:
$\quad r_n:= r_G\ast \phi_n \in C^1,\quad$
where $(\phi_n)$ is an approximation of the identity 
of class $C^1$, then:
\begin{align*}
	&\frl{x, y} |r_n(x) - r_n(y)| 
	= \left| \int_{\bR^d} (r_G(x-z) - r_G(y-z)) \phi_n(z) dz\right|\\
	&\hspace{5cm}
	\le \|r\|_{Lip}^G \, \|x - y\|\, \int_{\bR^d}\phi_n(z) dz 
	= \|r\|_{Lip}^G \, \|x - y\|.
\end{align*}
\begin{align*}
	\text{Thus }\hcm{2}
	\frlq{n \ge 1} \NGOM{r'_n} \le \|r\|_{Lip}^G,\quad
	\NGOM{r_n - r_G} \cvifty{n} 0.
	\hcm{2}
\tag*{$\square$}
\end{align*}

\paragraph{Proof that Lemma \ref{APG} and the case $r\in C^1$ 
	proves Theorem \ref{O_G}.}

We just have to prove $\Req{O_rYdY}$ 
with $\|r\|_{Lip}^G$ instead of $\NGOM{r'}$.
If we apply this formula for $r_n$ and exploit Lemma \ref{APG}, we see that there will be some $C = C(t_G, y_\vee, N_J) > 0$  such that:
\begin{align*}
	2\; \left|\int_{0}^{T_G} r_n(X_s)\, Y_s\, dY_s \right|
	\le \blp 2\,(N_J +1)\, \NGOM{r} + \|r\|_{Lip}^G\; v\, t_G\brp  \, (y_\vee)^2
	+ r_\infty \; t_G + C \, \NGOM{r - r_n}.
\end{align*}

Thus, it remains to bound:
\begin{align*}
	&\left|\int_{0}^{T_G} (r_n(X_s) - r(X_s) )\, Y_s\, dY_s \right|
	\le t_G\, y_\vee\, \psi_G^\vee \, \NGOM{r - r_n} + \left|M_n \right|,\\
	&
	\text{where } M_n:= \int_{0}^{T_G} (r_n(X_s) - r(X_s) )\, Y_s\, dB_s 
	\text{ has mean 0  and variance:}\\
	&\E\lp (M_n)^2 \rp 
	= \E\lp \int_{0}^{T_G} (r_n(X_s) - r(X_s) )^2\, Y_s\, ^2\, ds \rp\\
	&\hspace{1.5cm}
	\le t_G\, (y_\vee) ^2\,  (\NGOM{r - r_n}) ^2 \cvifty{n} 0.
\end{align*}
Thus, we can extract some subsequence $M_{\phi(n)}$ which converges a.s. towards 0. So that  a.s.:
\begin{align*}
	&\left|\int_{0}^{T_G} r(X_s) \, Y_s\, dY_s \right|
	\le \underset{n\rightarrow \infty}{\liminf} \Lbr
	\left|\int_{0}^{T_G} r_{\phi(n)}(X_s)\, Y_s\, dY_s \right|
	+ t_G\, y_\vee\, \psi_G^\vee \, \NGOM{r - r_{\phi(n)}} +
	\left|M_{\phi(n)}\right|\Rbr\\
	&\hspace{3cm}
	\le \frac{1}{2}\,\blp 2\,(N_J +1)\, \NGOM{r} + \|r\|_{Lip}^G\; v\, t_G\brp \, (y_\vee)^2
	+ \frac{1}{2}\, \NGOM{r} \; t_G < \infty.
\end{align*}
The proof in the case $r\in C^1$ can then be exploited without difficulty.\hfill $\square$

\subsection{Mixing for $Y$ }
 \label{O_Girsanov}


The proof will rely on Theorem \ref{O_G}
and on the following classical property of Brownian Motion:
\begin{lem}
	\label{O_outoftube}
	Consider any constants $b_\vee>0$,
	$\eps >0$ and $0<t_0 \le t_1$.
	Then, 
	there exists $c_B>0$ such that for any
	$b_I \in [0, b_\vee]$ and $t\in [t_0, t_1]$:
	\begin{align*}
		\PR_{b_I}\lp B_t \in db\pv
		\min_{s \le t_1} B_s \ge -\eps\mVg
		\max_{s \le t_1} B_s \le b_\vee + \eps \rp
		\ge c_B\times \idg{[0,\, b_\vee]}(b) \, db,
	\end{align*}
where $B$ under $\PR_{b_I}$
has by definition the law of a Brownian Motion with initial condition $b_I$.
\end{lem}

Thanks to this lemma and Theorem \ref{O_G}, 
we will be able to control $Y$
to prove that it indeed diffuses and 
that it stays in some closed interval $I_Y$ away from $0$.
We can then control the behavior of $X$
independently of the trajectory of $Y$ 
by appropriate conditioning of $M$ --the PRaMe--
so as to ensure the jumps we need
(conditionally that it remains in $I_Y$).
\\

\paragraph{Proof:}
Consider the collection of marginal laws of $B_t$,
with initial condition $b\in (-\eps, b_\vee +\eps)$,
killed when it reaches $-\eps$ or $b_\vee + \eps$.
It is classical that these laws have a density $u(t; b, b')$, $t>0$, 
$b'\in [-\eps, b_\vee +\eps]$,
w.r.t. the Lebesgue measure 
(cf e.g. Section 2.4 in Bass \cite{Ba95} for more details). 
$u$ is a solution to the Cauchy problem 
with Dirichlet boundary conditions:
\begin{align*}
	& \partial_t u(t; b_I, b) = \Delta_{b} u(t; b_I, b), 
	\hcm{2} \text{for }t>0, b_I, b\in (-\eps, b_\vee +\eps)
	\mVg
	\\&
	u(t; b_I,-\eps) =  u(t; b_I, b_\vee +\eps) = 0, \hcm{1}\text{for }t>0.
\end{align*}
Thanks to the maximum principle (cf e.g. Theorem 4, Subs 2.3.3. in Evans \cite{Evans}),
$u>0$ on $\bR_+^*\ltm [0, b_\vee] \ltm (-\eps, b_\vee +\eps)$ 
and since $u$ is continuous in its three variables, 
it is lower-bounded by some $c_B$ on the compact subset 
$[t_0, t_1]\ltm [0, b_\vee] \ltm [0, b_\vee]$.
\hfill $\square$

\subsection{Mixing for $X$}
 \label{sec_Xmix}
 
 For clarity, we decompose the "migration" along $X$ 
 into different kinds of elementary steps,
 as already done in \cite{AV_disc}.
 Let:
 \begin{align*}
 \cA:= \obB(-\theta\, \mathbf{e_1},\,\dS/2), \qquad
 \tau_\cA:= \inf\Lbr t\ge 0\pv 
 X_t \in \cA\mVg
 Y_t \in [2, 3] \Rbr,
 \EQn{O_Atr}{\cA}
 \end{align*}
 where we assume w.l.o.g. that $\dS \le \theta/8$
 ($[2, 3]$ is chosen arbitrarily).

 Under any of the three sets of assumptions
 considered in the following,
 the proof is achieved in three steps.
The first step is to prove that, 
 with a lower-bounded probability 
 for any initial condition in $\cD_\ell$, 
 $\tau_\cA$ is upper-bounded by some constant $t_\cA$.
In the second step, we prove that
 the process is sufficiently diffuse 
 and that time-shifts are not a problem.
 In the third step, we specify 
 which sets we can reach from $\cA$.
 \\
 
 Recall that for any $\ell\ge 1$,
 $T_{\cD_{\ell}}:= \inf\Lbr t\ge 0\pv (X,\, Y)_t \notin \cD_{\ell} \Rbr< \ext$.
 For $\np\ge 3$, 
 let us define 
 $T_{(\np)}:=T_{\cD_{2 \np}}$. For  $\np\ge 3$ and $t, c>0$, let:
 \begin{align*}
 &\cR^{(\np)}(t, c):= \Lbr x_F\in \bR^d \pv
 \frl{(x_0,y_0) \in \cA \times [1/\np,\,\np]}
 \right.
 \EQn{O_cR}{\cR^{(\np)}}
 \\&\hcm{2}
 \left. 
 \PR_{(x_0,y_0)} \lc (X, Y)_t\in (dx, dy)
 \pv t <T_{(\np)}\rc
 \ge c\; \idg{B(x_F, \dS/2)}(x)\, \idg{[1/\np,\,\np]}(y) \, dx\, dy \Rbr.
 \end{align*}

 We will prove the mixing on a global scale 
 by translating local mixing properties
 into some induction properties of the sets 
 $(\cR^{(\np)}(t, c))_{\{t, c>0\}}$.
 
 Several local mixing properties require 
 local lower- and upper-bounds on $g$,
so that they can only be exploited in specific areas of $\bR^d$.
In order to provide a general framework 
for these 
through Proposition \ref{pr_step},
let us consider the following sequence of sets,
indexed by $\np\ge 1$:
\begin{align*}
\cG_\np 
:= \{x\in \obB(0, \np)\pv  
&\frl{z\in [0, \eta/4]}\frl{\delta\in \obB(0, \dS/2)} \frl{w\in  \obB(\theta\, \mathbf{e_1},\dS)}
&g(x - (\theta - z) \mathbf{e_1} + \delta, w) \ge 1/\np,&
\\
\text{and }&\frl{z\in [-\theta, \eta/4]}\frl{\delta\in \obB(0, \dS/2)} \frl{w\in  \bR^d}
&g(x+z \mathbf{e_1} + \delta, w) \le \np
\}&.
\end{align*}

 These steps are deduced from the following elementary properties: 
 \begin{lem}
 	\label{O_lem_v}
 	Assume that Assumption $(H)$ hold.
 	Then, for any $\np\ge 1$ there exists $c_D>0$ 
 	such that for any $(x_I, y_I)\in \cD_ \np$
 	and $u\in [0, u_\vee(x_I)]$, where 
 	$u_\vee(x):= \sup\{ u\ge 0
 	\pv (x-v\, u\, \mathbf{e_1}) \in \obB(0, \np)\}$:
 	\begin{align*}
 		\PR_{(x_I,y_I)}\lc (X_u, Y_u) \in (dx, dy)
 		\pv u <  T_{(\np)}\rc
 		\ge c_D\, \delta_{\{x_I -v\, u\, \mathbf{e_1}\}}(dx)\ltm \idg{[1/\np, \np]}(y)\,dy.
 	\end{align*}
 	In particular,
 	for any $t, c>0$, $\np\ge3$,
 	the fact that $x$ belongs to $\cR^{(\np)}(t, c)$
 	implies the following inclusion:
 	\begin{align*}
 \frlq{u\in [0, u_\vee(x)]}
 x -v\, u\, \mathbf{e_1} \in \cR^{(\np)}(t+u, c\ltm c_D).
 	\end{align*}
 \end{lem}

The proof of Lemma \ref{O_lem_v}
being easily adapted from the one of the next proposition, 
it is deferred after the proof of the latter.
 
 \begin{prop}
\label{pr_step}
For any $\np \ge 3$, there exists $t_P, c_P>0$
such that for any $x_I \in \cG_\np$,
for any $x_0\in B(x_I, \dS/4)$ and $y_0\in [1/\np, \np]$:
 \begin{align*}
	\PR_{(x_0, y_0)} \lc
	(X, Y)_{t_P} \in (dx, dy)\pv t_P <T_{(\np)}\rc
	\ge c_P\; \idg{B(x_I,\, 3\dS/4)}(x)\,
	\idg{[1/\np,\,\np]}(y) \; dx\, dy.
\end{align*}
 \end{prop}

 	A direct application of the Markov property implies 
the two following results.
\begin{cor}
	\label{lem_Rpro}
	For any $\np \ge 3$, there exists $t_P, c_P>0$
	such that for any $t, c>0$, the following inclusion holds:
	\begin{align*}
		\{x\in \bR^d\pv d(x, \cR^{(\np)}(t, c) \cap \cG_\np)\le \dS/4\}
		\subset \cR^{(\np)}(t+t_P, c\ltm c_P).
	\end{align*} 
\end{cor}

\begin{fact}
	\label{O_lem_PropA}
	For any $t, t', c, c'>0$ and $\np\ge 1$:
	\begin{align*}
		\cA \cap \cR^{(\np)}(t, c) \neq \emptyset
		\Rightarrow \cR^{(\np)}(t', c')
		\subset \cR^{(\np)}(t+t', c\times c').
	\end{align*}
\end{fact}

\paragraph{Corollary \ref{lem_Rpro} as a consequence of Proposition \ref{pr_step}}
 
 For $n\ge 3$, let $t_P, c_P>0$ 
 be prescribed by Proposition \ref{pr_step}. 
  We consider $x_I\in \cR^{(\np)}(t, c) \cap \cG_\np$
 $x_F$ such that $\|x_F - x_I\| \le \dS/4$.
Combining through the Markov property 
the fact that $x_I\in  \cR^{(\np)}(t, c)$ 
and Proposition \ref{pr_step},
we deduce that for any $(x_0,y_0) \in \cA \times [1/\np,\,\np]$:
\begin{align*}
	& \PR_{(x_0,y_0)} \lc (X, Y)_{t+t_P}\in (dx, dy)
	\pv t+t_P < T_{(\np)}\rc
	\\&\hcm{1}
	\ge c \int_{B(x_I, \dS/2)} dx_0' \int_{1/\np}^{\np} dy'_0
	 \PR_{(x_0',y_0')} \lc (X, Y)_{t_P}\in (dx, dy)
	\pv t_P < T_{(\np)} \rc
	\\&\hcm{1}
	\ge c\times Leb(B(x_I, \eta/4))\times (\np-1/\np)
\times c_P	 \idg{B(x_I,\, 3\dS/4)}(x)\,
	\idg{[1/\np,\,\np]}(y) \; dx\, dy
	\\&\hcm{1}
		\ge c \ltm c'_P \idg{B(x_F,\, \dS/2)}(x)\,
		\idg{[1/\np,\,\np]}(y) \; dx\, dy,
\end{align*}
where $c'_P := Leb(B(0, \eta/4))\times (\np-1/\np)\times c_P>0$.
This means that $x_F \in \cR^{(\np)}(t+t_P, c\times c'_P)$.
The proof of Corollary \ref{lem_Rpro}
is thus concluded with this choice of $t_P$ and $c'_P$,
indeed independent from $x_I$, $x_F$.
\hfill $\square$

 \paragraph{Proof of Proposition \ref{pr_step}}
 
 \paragraph*{Step 1: description of the random event.}
For $\np\ge 3$, we set $t_P:= \theta/v$,
$t_J := \eta/(4 v)$,
$y_\wedge:= 1/(2\,\np)$, $y_\vee:= 2\, \np$.
 Let also:
 \begin{align*}
& T^Y:= \inf\Lbr t\ge 0\pv Y_t \notin [y_\wedge,\, y_\vee]  \Rbr,
  \quad  \EQn{O_fwL}{T^Y}
 \\
& f_\wedge
 := \inf \Lbr f(y)
 \pv y \in [y_\wedge,\, y_\vee] \Rbr, \qquad
 f_\vee
 := \sup \Lbr f(y)
 \pv y \in [y_\wedge,\, y_\vee]  \Rbr.
 \EQn{O_fvL}{f_\wedge, f_\vee}
 \end{align*}
 $f_\vee$ is finite due to \RHyp{HfPos}. 
Thanks to \RHyp{HfPos}, we know that
 $f_\wedge$ is positive.

 On the event $\{t_P < T^Y\}$,
 we shall prove that
 the values of $X$ on $[0, t_P]$ 
 are prescribed  as functions of $M$ 
 restricted to the subset:
 \begin{align*}
 \cX^M:= [0, t_P]\times \bR^d \times [0, f_\vee]\times [0,\np].
 \EQn{O_DomM}{\cX^M}
 \end{align*}
 
 Let $x_0:= x_I + \delta_0$ 
 with $x_I\in \cG_\np$ and $\delta_0\in B(0, \dS/4)$,
 and $y_0 \in [1/\np,\, \np]$
 that we consider as the initial conditions 
 for the process $(X, Y)$.

 To ensure one jump of size around $\theta$, 
 at time nearly $t_P$,
 while ``deleting'' the contribution of $\delta_0$, 
 let:
 \begin{align*}
 \cJ 
 := [t_P - t_J, t_P]
 \times B(\theta\, \mathbf{e_1} - \delta_0,\, \dS/2)
 \times [0, f_\wedge]
 \times [0, 1/\np].
 \EQn{O_JI}{\cJ}
 \end{align*}
 We partition  $\cX^M = \cJ \cup  \cN$, where:
 $\quad \cN:= \cX^M \setminus \cJ.$
The event mostly under consideration is the following: 
 \begin{align*}
 \cW = \cW^{(x_0, y_0)}:= 
 \Lbr t_P < T^Y \Rbr 
 \cap \Lbr M(\cJ) = 1\Rbr
 \cap \Lbr M(\cN) = 0\Rbr.
 \EQn{O_SbX}{\cW}
 \end{align*}

 Thanks to Theorem \ref{O_G}, 
 (with $x_\vee  := \np + 2 \theta$, $t_G =t_P$, 
 and the same values for $y_\wedge$ and $y_\vee$),
 there exists $c_G >0$ such that:
 \begin{align*}
 &\PR_{(x_0, y_0)} \lp
 (X, Y)_{t_P} \in (dx, dy)\pv
 \cW\rp
 \\&\hcm{2}
 \ge c_G\; \PR^G_{(x_0, y_0)} \lp
 (X, Y)_{t_P} \in (dx, dy)\pv
 \cW \rp.
 \EQn{O_PQA}{c_G}
 \end{align*}
 
 Under the law $\PR^G_{(x_0, y_0)}$, 
 the condition $\Lbr M(\cJ) = 1\Rbr$ is independent 
 of $\Lbr M(\cN) = 0\Rbr$, 
 of $\Lbr t_P < T^Y \Rbr$
 and of $Y_{t_P}$,
 cf Proposition \ref{O_Qind}.
 Thus, 
 on the event $\cW$,
 the only "jump" coded in the restriction of $M$ on $\cJ$
 is given as $(T\iJP\mVg \theta\, \mathbf{e_1} - \delta_0 + W\mVg U_f\mVg U_g)$,
 where $T\iJP, U_f$ and $U_g$ 
 are chosen uniformly and independently 
 on respectively $[t_P - t_J, t_P]$, $[0,f_\wedge]$ and $[0, 1/\np]$,
 and $\theta\, \mathbf{e_1} - \delta_0 + W$ 
 independently according to the restriction 
 of $\nu$ on $B(\theta\, \mathbf{e_1}- \delta_0,\, 3\dS/4)$
 (see notably chapter 2.4 in \cite{DV08}).
 Thanks to \RHyp{HdS},
 $W$ has a lower-bounded density $d_W$ on  $B(0,\, 3\dS/4)$.
 
 The following fact 
 motivates this description:
 \begin{fact}
 	\label{O_Xt0}
 	Under $\PR^G_{(x_0, y_0)}$ 
 consider	on the event $\cW$
 the r.v. $W= W_J - \theta \mathbf{e_1} + \delta_0$ 
 	 	where $(T\iJP, W\iJP, U_f, U_g)$ 
 	 	is the only point encoded by $M$ on $\cJ$.
 	 	Then, 
 a.s.	$X_{t_P} = x_I + W$
 and $\cW$ is included in $\{t_P <T_{(\np)}\}$.
 \end{fact}
 
 \paragraph{Step 2: proof of Fact \ref{O_Xt0}.}
 \subparagraph{Step 2.1.} We prove 
 that on the event $\cW$
 defined by \Req{O_SbX}:
 \begin{align*}
 \frlq{t< T\iJP}
 X_t:= x_0 - v\, t\, \mathbf{e_1}.
 \EQn{O_bUj}{}
 \end{align*}
 Indeed, $t_P< T^Y$ implies 
 that for any $t\le T\iJP$,
 $Y_t \in [y_\wedge,\, y_\vee]$. 
 Thanks to \Req{O_fvL},
 any "potential jump" $(T\iJP', W', U_f', U_g')$ 
 such that $T\iJP'\le  T\iJP$ and either $U_f' > f_\vee$ or $U_g' > \np$
 will be rejected.
 Thanks to the definition of $T\iJP$, 
 with \Req{O_DomM}, \Req{O_JI} and \Req{O_SbX},
 no other jump can occur, thus \Req{O_bUj} holds.
 \\
 
 Note that, in order to prove this rejection very rigorously, 
 we would like to consider the first one of such jumps.
 This cannot be done however for $(X, Y)$ directly,
 but is easy to prove for any approximation of $M$
 where $u_f$ and $u_g$ are bounded. 
 Since the result does not depend on these bounds 
 and the approximations converge to $(X, Y)$ 
 (and even equal to it before $T\iJP$ for bounds larger than $(f_\vee, \np)$),
 \Req{O_bUj} indeed holds.

 \subparagraph{Step 2.2.} We then prove that the jump at time $T\iJP$ is surely accepted.\\
Since $x_I\in \cG_\np$, by \Req{O_fwL}
 and the definition of $(T\iJP, W, U_f, U_g)$:
 \begin{align*}
 &U_f 
 \le f_\wedge
 \le f(Y_{T\iJP}),\quad
 U_g 
 \le 1/\np
 \le g(x_0 - v\, T\iJP\, \mathbf{e_1}\mVg \theta\, \mathbf{e_1} - \delta_0 + W)
 \\&\hcm{5.3}
 =  g(X_{T\iJP-}\mVg \theta\, \mathbf{e_1} - \delta_0 + W).
 \\&	\text{Thus } 
 X_{T\iJP} 
 = x_I + \delta_0 - v\, T\iJP\, \mathbf{e_1} + \theta\, \mathbf{e_1} - \delta_0 + W
 = x_I + (\theta - v\, T\iJP)\, \mathbf{e_1} + W.
 \end{align*}
 
 \subparagraph{Step 2.3.} We say that no jump can be accepted after $T\iJP$,
 which is proved as in Step 1.\\
 This means:
 $\frlq{T\iJP\le t \le t_P} 
 X_t = X_{T\iJP} - v\, (t-T\iJP)\,\mathbf{e_1} 
 = x_I + W.$\\
 
 This concludes in particular the proof of Fact \ref{O_Xt0} 
 with $t=t_P = \theta/v$.\hfill $\square$

 \paragraph{Step 3: concluding the proof of Proposition \ref{pr_step}.}
 Note that under $\PR^G$, 
 $\Lbr M(\cN) = 0\Rbr$
 is also independent
 of $\Lbr t_P < T^Y \Rbr$
 and  of $Y_{t_P}$,
 so that:
 \begin{align*}
 &\PR^G_{(x_0, y_0)} \lc
 (X, Y)_{t_P} \in (dx, dy)\pv
 \cW \rc
\\
 &= \PR[M(\cN) = 0]
 \times
 \PR(M(\cJ) = 1)\times
 \PR^G_{y_0} \lp Y_{t_P} \in dy\pv
 t_P < T^Y\rp
 d_W\, \idg{B(x_I,\, 3\dS/4)}(x)\, dx.
 \EQn{O_difX}{}
 \end{align*}
 From \Req{O_DomM} and \Req{O_JI}: 
 \begin{align*}
 &\PR(M(\cN) = 0)\; \PR(M(\cJ) = 1)
 = (t_J\, f_\wedge/\np)\ltm \nu\{ B(\theta\, \mathbf{e_1}- \delta_0,\, 3\dS/4)\}
 \times \exp[ - t_P\, f_\vee\,\np\, \nu(\bR^d) ] 
 \\&\hcm{2}
 \ge (t_J\, f_\wedge\,d_W/\np)\ltm Leb\{B(0, 3\dS/4)\} 
 \times \exp[ - t_P\, f_\vee\, \nu(\bR^d) ] 
 := c_X,
 \EQn{O_cX}{c_X}
 \end{align*}
 where the lower-bound $c_X$ is independent of $x_0$ and $y_0$.
 
 Thanks to Lemma \ref{O_outoftube}
 (recall the definitions of $y_\wedge$ and $y_\vee$ at the beginning of this subsection),
 \begin{align*}
 \PR^G_{y_0} \lp Y_{t_P} \in dy\pv
 t_P < T^Y\rp
 \ge c_B\; \idg{[1/\np,\,\np]}(y) \, dy
 \EQn{O_cB}{c_B}.
 \end{align*}
 Again, $c_B$ is independent of $x_0$ and $y_0$.
 
 With \Req{O_PQA},	\Req{O_difX},	\Req{O_cX}, \Req{O_cB},
Fact \ref{O_Xt0}
 and setting the value: $ c_P:= c_G\, c_X\, c_B\;d_W >0$:
 \begin{align*}
 &\frl{x_0 \in B(x_I, \dS/4)}
 \frl{y_0\in [1/\np,\,\np]}
 \\&\hcm{.5}
 \PR_{(x_0, y_0)} \lc
 (X, Y)_{t_P} \in (dx, dy)\pv t_P <T_{(\np)}\rc
 \ge c_P\; \idg{B(x_I,\, 3\dS/4)}(x)\,
 \idg{[1/\np,\,\np]}(y) \; dx\, dy.
 \end{align*}
 This ends the proof of Proposition \ref{pr_step}. \hfill $\square$
\paragraph{Proof of Lemma \ref{O_lem_v}}
 The proof of Lemma \ref{O_lem_v}
 relies on similar principles 
 as the one of Proposition \ref{pr_step}.
 In this case, $t_P$ is to be replaced by $u\in [0, u_\vee(x_I)]$
 and the event under consideration  is simply the following:
 \begin{align*}
 		\cW' := 
 		\Lbr u < T^Y \Rbr 
 		\cap \Lbr M([0, u]\times \bR^d \times [0, f_\vee]\times [0,\np]) = 0\Rbr.
 \end{align*}
 The reasoning given for Step 2.1. can be applied 
 to prove that  for any $t\le u$, 
$X_t:= x_0 - v\, t\, \mathbf{e_1}.$
We also exploit  Theorem \ref{O_G}
for the independence property between $X$ and $Y$ under $\PR^G_{(x_I, y_I)}$
and Lemma \ref{O_outoftube} to control the diffusion along the $Y$ coordinate.
Note that $c_B$ can be taken independently of $x_I, y_I$ and $t$
(noting that $t$ is uniformly upper-bounded by $2\np$).
These arguments conclude the proof of the lower-bound on the marginal of $(X, Y)$ on the event $\{t<T_{(\np)}\}$.
 
 The implication in term of the sets $\cR^{(\np)}(t, c)$
 is simply exploiting the Markov property,
 similarly as Corollary \ref{lem_Rpro} is deduced as a consequence of Proposition \ref{pr_step}.
 $\hfill \square$
 
 \subsection*{Application to the various sets of assumptions}
 
\subsection
{Proof of  Theorem \ref{O_Mix_DF} under Assumption \HgD}
\label{O_stabX}
We treat in this subsection the mixing of $X$ 
when both advantageous and deleterious 
mutations are occurring.
More precisely, 
each step corresponds to each of the following Lemmas:
\begin{lem}
	\label{O_lem_cRD}
	Assume that \HDef\, and
	\HgD\, hold.
	Then, for any $m \ge 3$, we can find $\np\ge m$ 
	and $c, t>0$
	such that
	$\obB(0, m)$ is included in $\cR^{(\np)}(t, c)$.
\end{lem}

%
\begin{lem}
\label{O_lem_Sb}
Assume that \HDef\, and
\HgD\, hold.
Then,
there exists $\np \ge 3$
which satisfies the following property
for any $t_1, t_2>0$.
There exists $t_R > t_1$ and $c_R>0$ 
such that for any $t\in [t_R, t_R+t_2]$ 
and $(x_0,y_0) \in \cA \times [2, 3]$:
\begin{align*}
\PR_{(x_0,y_0)} \lc (X, Y)_t\in (dx, dy)
\pv t < T_{(\np)}\rc
\ge c_R\; \idg{\cA}(x)\, \idg{[2, 3]}(y) \, dx\, dy.
\end{align*}
\end{lem}

\begin{lem}
	\label{O_lem_A}
	Assume that \HDef\, and
	\HgD\, hold.
	Then, for any $\ell_I>0$, there exists $c_I, t_I>0$
	and $\np \ge \ell_I$ 
	such that:
	\begin{align*}
		\frlq{(x, y)\in \cD_{\ell_I} }
		\PR_{(x, y)}( \tau_\cA \le t_I
		\wedge T_{(\np)}) \ge c_I.
	\end{align*}
\end{lem}

In the following Subsections, we prove these three lemmas then how  Theorem \ref{O_Mix_DF} is deduced as a consequence of these.

\subsubsection{Step 1: proof of Lemma \ref{O_lem_cRD}.}
\label{sec_lem_cRD}
Let $x_I = -\theta\, \mathbf{e_1}$.
Since $g$ is positive and continuous under Assumption \HgD,
there exists $\np_0$ such that $\obB(x_I, \dS/2)$ 
is included in $\cG_{\np_0}$.
With $t_0, c_0$ the values associated to $\np_0$
through Proposition \ref{pr_step},
we deduce that $x_I\in \cR^{(\np_0)}(t_0, c_0)$.

For $\mpr\ge 3$,
let $K := \Lfl 4\,\|\mpr+ \theta\| /\dS\Rfl+1$.
Similarly, we can choose $\np_1$ such that $B(0, \mpr)$
is a subset of $\cG_{\np_1}$.
Consider any $x_F\in \obB(0, \mpr)$ and
for $0\le k\le K$, let 
\mbox{$x_k:= -\theta\, \mathbf{e_1} + k/K\,(x_F+ \theta\, \mathbf{e_1})$.} 
This choice is made to ensure that 
$d(x_k, x_{k+1})\le \dS/4$ and $\frl{k\le K} x_k\in \cG_{\np_1}$
Thanks to  Corollary \ref{lem_Rpro}, 
we deduce by immediate induction over $k\le K$
that there exist $\np_2, t_k, c_k >0$ independent of $x_F$
such that:
$x_k \in \cR^{(\np_2)}(t_k, c_k)$.
$t_k$ and  $c_k$  are of the form $t_k := t_0 + k\, t_P$ and $c_k :c_0\ltm (c_P)^k$.
In particular with $k=K$, and $n:=n_2$,
Lemma \ref{O_lem_cRD} is proved.
\hfill $\square$

\subsubsection{Step 2: proof of Lemma \ref{O_lem_Sb}.}
\label{sec_Sb}
 We keep $x_I:= -\theta\, \mathbf{e_1}$ 
 and $x_1:= (-\theta+\dS/2)\, \mathbf{e_1}$.
Thanks to Lemma \ref{O_lem_cRD},
there exists $\np, t_1, c_1>0$ such that:
\begin{align*}
	\Lbr\, x_I + u\, \mathbf{e_1}\pv u\in [\dS/6\mVg 5\,\dS/6]\Rbr 
	\subset \cR^{(\np)}(t_1, c_1).
\end{align*}
There exists $t_2, c_2 >0$
thanks to  Lemma \ref{O_lem_v} such that:                         
$\frlq{t\in [t_2, t_2 + 2\,\dS/(3\,v)]}$
$x_I \in \cR^{(\np)}(t, c_2).$
Applying twice Corollary \ref{lem_Rpro}
with the knowledge that $B(x_I, \dS/2)$ is a subset of $\cG_\np$,
we deduce that there exists $t_3, c_3>0$
such that:
$$\frlq{t\in [t_3,  t_3 + 2\,\dS/(3\,v)]}
\cA \subset \cR^{(\np)}(t, c_3).$$
Applying inductively Fact \ref{O_lem_PropA}, 
we deduce the following for any $k\ge 1$:
$$\frlq{t\in [k\,t_3,  k\,t_3 + 2\,k\,\dS/(3\,v)]}
\cA  \subset \cR^{(\np)}(t, [c_3] ^k).$$
Let $t_1, t_2>0$ and consider $k\ge 1$ sufficiently large 
for $k\,t_3>t_1$ and $2\,k\,\dS/(3\,v)>t_2$ to hold.
Thus, Lemma \ref{O_lem_Sb} is proved with this value of $\np$,
$t_R:= k\,t_3$ and $c_R :=  [c_3] ^k$.
\hfill $\square$

\subsubsection{Step 3: proof of Lemma \ref{O_lem_A}.}
\label{sec_lemA}

As before, we can find $\np\ge \ell_I$ be such that $ \cD_{\ell_I}\subset \cG_\np$.
We go backwards in time from $\cA$ by defining 
for $t\ge 0$, $c>0$:
\begin{align*}
	&\cR'(t, c):= \Lbr (x, y)\in \cG_\np\pv
	\PR_{(x,y)} \lc \tau_\cA \le t \wedge T_{(\np)}\rc \ge c \Rbr.
\end{align*}
It is clear that $\cA\subset \cR'(0, 1)$.
Thanks to Proposition \ref{pr_step}
and the Markov property,
there exists $t_P, c_P>0$
such that for any $t, c>0$:
\begin{align*}
	\Lbr x\in \cG_\np\pv d(x, \cR'(t, c))\le \dS/4\Rbr
	\subset \cR'(t+t_P, c\ltm c_P).
\end{align*}
Since $\cD_{\ell_I}\subset \cG_\np$ is bounded,
an immediate induction ensures that there exists $t_I, c_I>0$ such that $\cD_{\ell_I}\subset \cR'(t_I, c_I)$.
This concludes the proof of Lemma \ref{O_lem_A}.
\hfill $\square$

\subsubsection
{Theorem \ref{O_Mix_DF}
	as a consequence of Lemmas \ref{O_lem_cRD}-3 }
\label{O_pr_MixD}

The proof is quite naturally adapted from the one of Lemma 3.2.1 in \cite{AV_disc}.
Note that 
 for any $\np_1\le \np_2$, $T_{(\np_1)}\le T_{(\np_2)}\le \ext$ 
holds a.s.

Let $\ell_I, \ell_M\ge 0$.
According to Lemma \ref{O_lem_A}, we can find $c_I, t_I>0$ 
and $\np_1\ge \ell_I\wedge \ell_M$ 
such that 	for any $(x_I, y_I)\in \cD_{\ell_I}$:
\begin{align*}
	\PR_{(x_I, y_I)}( \tau_\cA \le t_I\wedge T_{(\np_1)} ) \ge c_I.
			\EQn{O_tA}{}
\end{align*}
Let also $\np_2 \ge \np_1$, $c_R,t_R>0$ 
chosen thanks to Lemma \ref{O_lem_Sb}
to satisfy that for any $t\in [t_R, t_R+t_I]$ 
and $(x_0,y_0) \in \cA \times [2, 3]$:
\begin{align*}
	\PR_{(x_0,y_0)} \lc (X, Y)_t\in (dx, dy)
	\pv t < T_{(\np_2)} \rc
	\ge c_R\; \idg{\cA}(x)\, \idg{[2, 3]}(y) \, dx\, dy.
			\EQn{O_tR}{}
\end{align*}

Thanks to Lemma \ref{O_lem_cRD}, 
since $\cD_{\ell_M}$ is a bounded set,
we know that there exists $\np\ge \np_2$,
$c_F$ and $t_F> 0$ such that for any $(x_0,y_0) \in \cA \times [2,\,3]$:
\begin{align*}
	\PR_{(x_0,y_0)} \lc (X, Y)_{t_k}\in (dx, dy)
	\pv t_k < T_{(\np)} \rc
	\ge c_F\; \idg{\cD_{\ell_M}}(x)\, \idg{[1/\np,\,\np]}(y) \, dx\, dy.
	\EQn{O_tF}{c_F,t_F}
\end{align*}
The fact that $\np$ is larger than $\np_1$ and $\np_2$ 
implies without difficulty 
that \Req{O_tA} and \Req{O_tR} hold with $\np_1$ and $\np_2$
 replaced by $\np$, which is how these statements are exploited in the following reasoning.

Let $t_M:= t_I + t_R + t_F$ 
and $c_M:= c_I \ltm  c_R \ltm Leb(\cA) \ltm c_F$.
For any $(x_I, y_I)\in \cD_{\ell_I}$,
by combining \Req{O_tR}, \Req{O_tF} and the Markov property,
we deduce that a.s. on the event 
$\{\tau_\cA \le t_I\wedge T_{(\np)}\}$:
\begin{align*}
&		\PR_{(X, Y)[\tau_\cA]} \lc
	(\wtd{X}, \wtd{Y})[t_M - \tau_\cA]
	\in (dx, dy)
	\pv  t_M- \tau_\cA < \wtd{T}_{(\np)} \rc
	\\&\hcm{1}\ge 
	c_F\ltm \PR_{(X, Y)[\tau_\cA]} \lc
	(\wtd{X}, \wtd{Y})[t_M - t_F- \tau_\cA]
	\in \cA\times [2, 3]
	\pv  t_M - t_F - \tau_\cA < \wtd{T}_{(\np)} \rc
	\\&\hcm{1.5}
	\ltm \idg{\cD_{\ell_M}}(x)\, \idg{[1/\np,\,\np]}(y) \, dx\, dy
	\\&\hcm{1}\ge 
		c_R \ltm Leb(\cA) \ltm c_F
			\ltm \idg{\cD_{\ell_M}}(x)\, \idg{[1/\np,\,\np]}(y) \, dx\, dy,
\end{align*}
where we exploited the knowledge that $\tau_\cA\le t_I$ 
to deduce that $t_M - t_F - \tau_\cA \in [t_R, t_R+t_I]$.
By combining this estimate with \Req{O_tA}
and again the Markov property, 
we conclude:
\begin{align*}
	&\PR_{(x_I, y_I)} \lc
	(X_{t_M}, Y_{t_M})	\in (dx, dy)
	\pv  t_M< T_{(\np)} \rc
		\\&\hcm{1}
	\ge \PR_{(x_I, y_I)}( \tau_\cA \le t_I\wedge T_{(\np)} ) 
	\ltm c_R \ltm Leb(\cA) \ltm c_F
	\ltm \idg{\cD_{\ell_M}}(x)\, \idg{[1/\np,\,\np]}(y) \, dx\, dy
	\\&\hcm{1}
	\ge c_M\, \idg{\cD_{\ell_M}}(x)\, \idg{[1/\np,\,\np]}(y) \, dx\, dy.
\end{align*}
This ends the proof of Theorem \ref{O_Mix_DF}
with $\Lp=2\np$, $c:= c_M$ and $t:=t_M$
under Assumption \HgD.

\hfill $\square$

\subsection
{Proof of  Theorem \ref{O_Mix_DF} under Assumption \HgA\, and $d\ge 2$}
\label{sec_Mix_DF}

The proof of Theorem \ref{O_Mix_DF} 
is handled under Assumption \HgA\, and $d\ge 2$
in the same way as in Subsection \ref{O_pr_MixD}.
Notably, the lemmas that replace  Lemmas \ref{O_lem_Sb}-3
have identical implications:

\begin{lem}
	\label{lem_cRD}
	Assume that $d\ge 2$ and that  \HDef\, and
	\HgA\, hold.
	Then, for any $m \ge 3$, we can find $\np\ge m,$
	$t, c>0$ such that
	$\obB(0, m)$ is included in $\cR^{(\Lp)}(t, c)$.
\end{lem}

\begin{lem}
	\label{lem_Sb}
	Assume that $d\ge 2$ and that  \HDef\, and
	\HgA\, hold.
Then,
there exists $\np \ge 3$
which satisfies the following property
 for any $t_1, t_2>0$.
There exists $t_R > t_1$ and $c_R>0$ such that, 
for any $t\in [t_R, t_R+t_2]$ 
and $(x_0,y_0) \in \cA \times [2, 3]$:
\begin{align*}
	\PR_{(x_0,y_0)} \lc (X, Y)_t\in (dx, dy)
	\pv t < T_{(\np)}\rc
	\ge c_R\; \idg{\cA}(x)\, \idg{[2, 3]}(y) \, dx\, dy.
\end{align*}
\end{lem}

\begin{lem}
	\label{lem_A}
	Assume that $d\ge 2$ and that \HDef\, and
	\HgA\, hold.
	Then, for any $\ell_I>0$, there exists $c_I, t_I>0$
	and $\np \ge \ell_I$ 
	such that:
	\begin{align*}
		\frlq{(x_0, y_0)\in \cD_{\ell_I} }
		\PR_{(x_0, )}( \tau_\cA \le t_\cA 
		\wedge T_{(\np)} ) \ge c_\cA.
		\EQn{O_tA1}{}
	\end{align*}
\end{lem}

Since the implications are the same,
the proof of Theorem \ref{O_Mix_DF}
under Assumption \HgA\, with $d\ge 2$
as a consequence of Lemmas \ref{lem_cRD}-3
is mutatis mutandis the same as the one
given in Subsection \ref{O_pr_MixD}.
Since deleterious mutations 
are now forbidden,
the proof of Lemma \ref{lem_cRD}
is much trickier than the one of Lemma \ref{O_lem_cRD}.
The first step is given by the two following lemmas.
To this purpose, 
given any direction $\bU$ on the sphere  $S^d$ of radius $1$,
we denote its orthogonal component by:
\begin{align*}
	x^{(\perp \bU)}:= x - \LAg x, \bU\RAg \bU\mVg
	\quad
	\text{ and specifically for $\mathbf{e_1}$:}\quad
	x^{(\perp 1)}:= x - \LAg x, \mathbf{e_1}\RAg \mathbf{e_1}.
	\EQn{O_xu}{x^{(\perp \bU)}}
\end{align*}

\begin{lem}
	\label{O_lem_RSA}
	Assume that $d\ge2$, \HDef\, and
	\HgA\, hold.
	Then, for any $x_\vee> 0$,  
	there exists $\eps \le \dS/8$
	 which satisfies the following property
	for any $\np \ge 3\vee (2\, \theta)$,  
	$x\in B(0, \np)$ and $\bU\in S^d$
	such that both $\langle x,\bU\rangle \ge \theta$ 
	and $\|x^{(\perp \bU)}\|\le x_\vee$.
There exists $t_P, c_P>0$ such that for any $t, c>0$:
	\begin{align*}
		x \in \cR^{(\np)}(t, c)
		\Rightarrow \obB(x - \theta\, \bU, \eps)\subset \cR^{(\np)}(t+t_P, c\times c_P).
	\end{align*}
\end{lem}

\begin{lem}
	\label{O_lem_RlocA}
	Assume that $d\ge2$, \HDef\, and
	\HgA\, hold.
	Then, for any $m \ge 3\vee (2\, \theta)$, there exists $\eps \le \dS/8$ which satisfies the following property
		for any $x\in B(0, m)$ with $\LAg x, \mathbf{e_1} \RAg\le 0$.
There exists $t_P, c_P>0$ such that:
	\begin{align*}
		\frlq{t, c>0}
		x \in \cR^{(\Lp)}(t, c)
		\Rightarrow \obB(x, \eps)\subset \cR^{(\Lp)}(t+t_P, c\times c_P).
	\end{align*}
\end{lem}

Lemma \ref{O_lem_RlocA} is actually directly 
implied from Lemma \ref{O_lem_v} 
(first applied for a time-interval $[0, \theta/v]$)
then Lemma \ref{O_lem_RSA} with $\bU :=\mathbf{e_1}$,
combined with the Markov property.
Subsection \ref{sec_RSA} is dedicated to the proof of Lemma \ref{O_lem_RSA}.


\subsubsection{Step 1: proof of Lemma \ref{O_lem_RSA}.}
\label{sec_RSA}

Fix $x_\vee>0$. Consider $\eps>0$ 
that is to be fixed later,
but assume already that $\eps \le \theta/8$.
We recall that $\eta\le \theta/8$ si assumed w.l.o.g.
Let $\np \ge 3\vee(2\theta)$,
$x_0\in B(0, \np)$ and $\bU \in S^d$ such that 
both $\langle x_0,\bU\rangle \ge \theta$ 
and $\|x_0^{(\perp \bU)}\|\le x_\vee$ hold.

Compared to Proposition \ref{pr_step},
the first main difference is that the jump is now 
almost instantaneous.
The second is that, in order that $g_\wedge >0$, 
we have way less choice in the value of $w$ 
when $\|x^{(\perp \bU)}\|$ is large.
In particular, 
the variability of any particular jump 
will not be sufficient to wipe out 
the initial diffusion around $x$ 
deduced from $x \in \cR^{(\np)}(t, c)$,
but will rather make it even more diffuse.

To fix $\eps>0$, let us first compute, 
for $\delta\in B(0, \dS)$, $\wP\in B(-\theta\, \bU, \eps)$:
\begin{align*}
	&\|x_0+\delta\|^2 - \|x_0+\delta + \wP\|^2 
	= 2\, \langle x_0 + \delta\mVg \wP\rangle
	- \|\wP\|^2 
	\\&\hcm{1}
	\ge \lp\frac{7}{4} - \frac{9}{8}
	\ltm(\frac{1}{4}+\frac{9}{8})\rp\, \theta^2- 2\, \eps \; x_\vee,
\end{align*}
where we exploited that 
$\langle \bU\mVg \wP\rangle \ge 7\, \theta/8$.
We note that:
\begin{align*}
	c:= \frac{7}{4} - \frac{9}{8}
	\ltm(\frac{1}{4}+\frac{9}{8})
	= \frac{13}{64}>0.
\end{align*}
By taking $\eps:= \{c\,\theta^2/(4\,x_\vee)\}\wedge \{\theta/8\}$, we thus ensure that
$\|x_0+\delta\|^2 
>  \|x_0+\delta+\wP\|^2 $.
Note that $\eps$ does not depend on the specific choice of $x_0$.

Let $t_P:=\eps/(2\,v)$. 
The initial condition for $X, Y$ is taken as $x_I\in B(x, \dS/2)$ and $y_I\in [1/\np,\,\np]$.
\begin{align*}
	&g_\wedge
	:= \inf\Lbr\, g(x, w)\pv 
	x \in \obB(x_0, \dS)\mVg
	w\in \obB(-\theta\,\bU, \eps)\Rbr >0,
	\\&
	\cX^M:= [0, t_P]\times \bR^d \times [0, f_\vee]\times [0,\np],
	\\&
	\cJ 
	:= [0, t_P]
	\times B(-\theta\, \bU +(\eps/2)\, \mathbf{e_1} \mVg 	\eps/2)
	\times [0, f_\wedge]
	\times [0, g_\wedge].
\end{align*}
With the same reasoning 
as in  the proof of Proposition \ref{pr_step},
we obtain a change of probability 
$\PR^G_{(x_I, y_I)}$ and an event $\cW$ 
on which the r.v. $W$ is uniquely defined from $M$
under $\PR^G_{(x_I, y_I)}$
and such that it satisfies a.s.: 
\begin{align*}
	X_{t_P} = x_I - (\eps/2)\, \mathbf{e_1}  - \theta\, \bU 
	+ (\eps/2)\, \mathbf{e_1} + W
= x_I - \theta\, \bU  + W,
\end{align*}
where the density of $W$ is lower-bounded by $d_W$ on $B(0, \eps/2)$,
uniformly over $x_I$ (given $x$), and $y_I$.
We thus similarly obtain some constants $c_P, c'_P >0$ 
independent of $x_0$
such that
for any such $x_0$:
\begin{align*}
& \int_{B(x_0, \dS/2)}\,dx_I \int_{[1/\np,\,\np]}\, dy_I\,
 	\PR_{(x_I, y_I)} \lc
 (X, Y)_{t_P} \in (dx, dy)\pv t_P <T_{(\np)}\rc
\\&\hcm{1}
\ge c_P\,  \int_{B(x_0, \dS/2)}\,dx_I
\idg{B(x_I - \theta\, \bU,\, \eps/2)}(x)
\ltm
\idg{[1/\np,\,\np]}(y) \; dx\, dy
\\&\hcm{1}
\ge c_P'\,  \idg{B(x_0 - \theta\, \bU,\,\dS/2 + \eps/3)}(x)
\idg{[1/\np,\,\np]}(y) \; dx\, dy.
\end{align*}
We then reason similarly as in the proof of Corollary \ref{lem_Rpro} 
as a consequence of Proposition \ref{pr_step}.
Assuming further that $x_0 \in \cR^{(\np)}(t, c)$
for some $t,c >0$,
we can deduce:
$$B(x- \theta\,\bU\mVg \eps /3) \in \cR^{(\np)}(t+t_P, c\times c_P).$$
This is exactly the implication of Lemma \ref{O_lem_RSA},
stated in terms of  $\eps/3$ instead of $\eps$. 
\hfill $\square$

\subsubsection
{Step 2: Lemma \ref{lem_cRD}
as a consequence of Lemmas \ref{O_lem_RlocA} and \ref{O_lem_RSA}.}
\label{sec_cRD}

\paragraph{Step 2.1: $x_I\in \cR^{(\np_0)}(t_0, c_0)$.}
Let $x_I:= -\theta \mathbf{e_1}$. 
We check that there exists $\np_1\ge 1$
 such that  $B(x_I, \dS/2)$ is a subset of  $\cG_{\np_1}$.
Since $g$ is continuous and thanks to \HgA, 
it is sufficient to prove that 
$\|x_I-z \mathbf{e_1} + \delta\| > \|x_I-z \mathbf{e_1} + \delta +w\| $
holds for any $z\in [0, \theta]$, $\delta\in \obB(0, \dS)$,
and $w\in  \obB(\theta\, \mathbf{e_1},\dS)$:
	\begin{align*}
	&\|x_I-z \mathbf{e_1}  + \delta\|^2
	- \|x_I-z \mathbf{e_1}  + \delta + w\|^2
	=  2 \langle (\theta +z) \mathbf{e_1} - \delta\mVg w\rangle
	- \|w\|^2
	\\&\hcm{1}
	\ge 2\, [\theta\ltm (\theta - \dS) - \dS\ltm(\theta+\dS)] - (\theta+\dS)^2 
		\\&\hcm{1}
	= \theta^2 -6\, \theta\,\dS - 3\, \dS^2
	\ge \dfrac{13 \theta^2}{64}  > 0,
\end{align*}
	since $\dS \le \theta/8$, as assumed above, just after \Req{O_Atr}.
Applying twice Proposition \ref{pr_step},
	it concludes that there exists $\np_0\ge 1$, $t_0, c_0>0$ 
	such that $x_I\in \cR^{(\np_0)}(t_0, c_0)$.

\paragraph{Step 2.2: under the condition that $\LAg x_F, \mathbf{e_1}\RAg:= -\theta$.}
The purpose of this step is the following lemma,
in which we employ the notation $\pi_1:x\mapsto \LAg x, \mathbf{e_1}\RAg$.
\begin{lem}
	\label{lem_e1OT}
For any $\np\ge 1$ sufficiently large,
there exists $t, c>0$ such that
$\pi_1^{-1}(-\theta)\cap B(0,\np)$
is a subset of $\cR^{(\np)}(t, c)$.
\end{lem}

Let $x_F \in \pi_1^{-1}(-\theta)\cap B(0,\np)$,
where we assume that $\np$ is larger than $\np_0$, 3 and $2\theta$.
First, we define $\bU$ as $\mathbf{e_1}$ if $x_F^{(\perp 1)}=0$ and else as $\bU:= x_F^{(\perp 1)}/\|x_F^{(\perp 1)}\|.$
Note that $\|x_F^{(\perp 1)}\|\le \np$.
We consider the value of $\eps$ 
given by Lemma \ref{O_lem_RlocA}
for $x_\vee := \np$
and define: 
\begin{align*}
K:= \Lfl \np \eps\Rfl +1 \mVg
\text{ for } 0\le k\le K,\quad
x_k:= -\theta\, \mathbf{e_1} + \frac{k\, \|x_F^{(\perp 1)}\|}{K}\; \bU.
\end{align*}
This choice ensures that for any $k\in \II{0, K-1}$, 
$x_{k+1}\in B(x_k, \eps)$, 
while $x_k\in B(0, \np)$, $\LAg x_{k}\bv \mathbf{e_1}\RAg \le 0$
and $x_K = x_F$.
Thanks to Step 2.1, $x_0\in \cR^{(\np)}(t_0, c_0)$.
Thus, by induction over $k\le K$ with Lemma \ref{O_lem_RlocA}, 
$x_k \in \cR^{(\np)}(t_0 + k\, t_P,\, c_0\, [c_P]^k)$.
In particular, 
there exists $t, c>0$ such that $x\in \cR^{(\np)}(t, c)$,
which concludes Step 2.2.

\paragraph{Step 2.3: the general case}

Assume solely that $x \in \cB(0, m)$.
We consider the value of $\eps$ 
given by Lemma \ref{O_lem_RSA}
for $x_\vee:= \mpr$.
The choice of $\bU$ is as in Step 2.2.

Let:
\begin{align*}
&\hcm{3}
K:= \Lfl\frac{m + \theta}{\eps}\Rfl +1,\quad
\text{so that }
\frac{\LAg x, \mathbf{e_1}\RAg + \theta}{K} \le \eps,
\EQn{O_Lst}{K}
\\&
\text{and for } 0\le \kp \le K,\quad
x_\kp:= (-\theta+ (\kp/K)\ltm (\LAg x, \mathbf{e_1}\RAg + \theta)) \, \mathbf{e_1} 
 +(K- \kp)\,\theta\, \bU+ x^{(\perp 1)}.
\end{align*}
In particular $\LAg x_0,\mathbf{e_1}\RAg = -\theta$,
$x_K = x_F$
while for any $k\le K-1$, $x_{k+1}\in B(x_k, \eps)$,
$x_k \in B(0, \mpr+K\, \theta)$
and $\LAg x_k, \bU\RAg \le \theta\vee \LAg x, \mathbf{e_1}\RAg\le \mpr=x_\vee$.

Since $\LAg x_0,\mathbf{e_1}\RAg = -\theta$,
we can exploit Lemma \ref{lem_e1OT} to prove that
 there exists $\np\ge 1$ and $t_0, c_0>0$ independent of $x_F$
 such that $x_0\in \cR^{(\np)}(t_0, c_0)$.
Thanks to Lemma \ref{O_lem_RSA} and induction on $\kp$,
we deduce that there exist $t_P, c_P>0$ such that:
$x_\kp 
\in \cR^{(\np)}(t_0 + \kp\, t_P,\, c_0\, [c_P]^\kp)$.
In particular, there exists $t,c>0$ such that 
$x_F \in \cR^{(\np)}(t, c)$.
\hfill $\square$

\subsubsection
{Step 3: proof of Lemma \ref{lem_Sb}.}
The proof can be taken mutatis mutandis from the one given 
in Subsection \ref{sec_Sb}.
The fact that $B(x_I, \dS/2)$ is a subset of  $\cG_{\np_1}$
is already proved in Step 2.1 (cf Subsection \ref{sec_cRD}),
while Lemma \ref{lem_cRD} 
replaces Lemma \ref{O_lem_cRD} with identical implication.
\hfill $\square$

\subsubsection
{Step 4: proof of Lemma \ref{lem_A}.}
\label{sec_lem_A}
\begin{rem}
The presented proof 
efficiently exploits the already known lemmas 
but is probably very far from optimal in its estimations.
\end{rem}

\paragraph{Step 4.1: study of $\cG_\np$.}
We look for conditions on $x\in \bR^d$
that ensures  that it belongs to  $\cG_\np$ for some $\np$.
Let $x_\theta:= x -(\theta -\dS/2)  \mathbf{e_1}$.
By definition of $\cG_\np$,
it is necessary that for any $z\in [0, \eta/4]$, $\delta\in \obB(0, \dS/2)$ and $w\in  \obB(\theta\, \mathbf{e_1},\dS)$
$g(x_\theta - z \mathbf{e_1} + \delta, w) >0$ which under \HgA\, is equivalent to 
$\|x_\theta - z  \mathbf{e_1} + \delta\|> \|x_\theta + z \mathbf{e_1} + \delta+w\|$.
We first restrict ourselves to the values of $x$ such that $\pi_1(x) \le 0$,
and we compute:
\begin{align*}
&	\|x_\theta - z \mathbf{e_1} + \delta\|^2 
-  \|x_\theta +z \mathbf{e_1} + \delta+w\|^2
=  - 2 \langle x_\theta +z \mathbf{e_1} + \delta, w\rangle
	- \|w\|^2
	\\&\hcm{1}
	\ge 2\, (-\pi_1(x_\theta)-\dS/2)\ltm (\theta-\dS)
	- 2(\|x^{(\perp 1)}\|+ \dS/2)\ltm \dS
	- (\theta+\dS)^2 
	\\&\hcm{1}
\ge (-7\pi_1(x_\theta)/32 - \|x^{(\perp 1)}\|/4)\ltm \theta
+ (7/4)\ltm  (\theta-\dS/2) \ltm  (\theta-\dS)
- \dS\ltm  (\theta-\dS)- \dS^2 - (\theta+\dS)^2.
	\\&\hcm{1}
\ge	(-7\pi_1(x_\theta)/32 - \|x^{(\perp 1)}\|/4)\ltm \theta
	+ (7\ltm 15\ltm 7- 8\ltm 7- 8 -8\ltm 81) \theta^2 /2^9
		\\&\hcm{1}
	\ge	(-7\pi_1(x_\theta)/32 - \|x^{(\perp 1)}\|/4)\ltm \theta
	+ 23\theta^2/2^9.
\end{align*}
From these computations, we see that $g(x_\theta - z \mathbf{e_1} + \delta, w) >0$ 
holds true provided $\pi_1(x)\le 0$ 
and $|\pi_1(x_\theta)|\ge 8 \|x^{(\perp 1)}\|/7$
thus a fortiori if $|\pi_1(x)|\ge 8 \|x^{(\perp 1)}\|/7$.
Since $g$ is continuous, we deduce that for any $\mpr\ge 1$,
there exists $\np\ge 1$ such that $\cG_\np$
contains the following set:
\begin{align*}
	\{x\in B(0, \mpr)\pv -\pi_1(x) \ge 8 \|x^{(\perp 1)}\|/7\}.
\end{align*}

\paragraph{Step 4.2.}
Let $\ell_I\ge 1$.
Thanks to Step 4.1, we can find $\np\ge \ell_I\vee 3$ 
such that $\cG_\np$
contains the following set:
\begin{align*}
\cA_1:=	\{x\in B(0, 2 \ell_I)\pv -\pi_1(x)\ge 8 \|x^{(\perp 1)}\|/7\}.
\end{align*}
We go backwards in time from $\cA$ by defining 
for $t\ge 0$, $c>0$:
\begin{align*}
	&\cR'(t, c):= \Lbr (x_I, y_I)\in \cG_\np\pv
	\PR_{(x_I,y_I)} \lc \tau_\cA \le t \wedge T_{(\np)}\rc \ge c \Rbr.
\end{align*}
Similarly as for the proof of Lemma \ref{O_lem_A},
by exploiting inductively Proposition \ref{pr_step},
we deduce that $\cA_1$ is a subset of 
$\cR'(t_1, c_1)$ for some $t_1, c_1>0$.

Consider now any $x_I\in \obB(0, \ell_I)$.
If $x_I \notin \cA_1$,
let $u_*:= 8 \|x^{(\perp 1)}\|/(7v) + \pi_1(x)/v$
and $x_1 := x- v \,u_*\mathbf{e_1}\in \cA_1$.
If $x_I \in \cA_1$, we simply define $x_1:= x_I$ and $u_*:=0$.
Since $\|x^{(\perp 1)}\|\le \ell_I$,
this choice necessarily satisfies
$0\le -\pi_1(x_1)=8 \|x^{(\perp 1)}\|/7\le 8\np/7$.
In any case, $x_1\in  B(0, 2 \ell_I)$ thus $x_1\in \cA_1$.
Since $\cA_1\subset \cR'(t_1, c_1)$ and thanks to Lemma  \ref{O_lem_v},
there exists a value $c_D>0$ uniform over $x$
such that $x_I\in \cR'(t_1+u_*, c_1\ltm c_D)$.
Since $u_*$ is upper-bounded by $2\ell_I$ 
and the set $\cR'(t, c)$
are increasing with $t$,
it concludes that $\obB(0, \ell_I)$
is a subset of $\cR'(t_2, c_2)$ with $t_2 := t_1 + 2\ell_I$ 
and $c_2 := c_1\ltm c_D$.
This ends the proof of Lemma \ref{lem_A}.
\hfill $\square$
\\

As mentioned at the beginning of Subsection \ref{sec_Mix_DF},
the last step of the proof of Theorem \ref{O_Mix_DF}
can be taken mutatis mutandis from Subsection \ref{O_pr_MixD}.
With this, the proof of the theorem is complete.

\subsection
{Proof of  Theorem \ref{O_Mix_AF}}
\label{O_stabXA}

We treat in this subsection the mixing for $X$ 
when simply advantageous mutations are occurring
and the phentotype is unidimensional.
The proof of Theorem \ref{O_Mix_DF} 
is handled under Assumption \HgA\, and $d\ge 2$
in the same way as in Subsection \ref{O_pr_MixD},
except that Lemmas \ref{O_lem_Sb}-3
are replaced by the following ones,
in respective order.
Note that only the first one has a different implication.

\begin{lem}
	\label{O_lem_cRd1}
	Assume that $d = 1$, \HDef\, and
	\HgA\, hold.
	Then, for any $\mpr\ge 3$, there exists $n\ge\mpr$, $t, c>0$ such that
	$[-\mpr, 0]$ is included in $\cR^{(\np)}(t,c).$
\end{lem}

\begin{lem}
	\label{lem_SbA}
	Assume that $d= 1$ and that  \HDef\, and
	\HgA\, hold.
Then,
there exists $\np \ge 3$
which satisfies the following property
 for any $t_1, t_2>0$.
There exists $t_R > t_1$ and $c_R>0$ such that, 
for any $t\in [t_R, t_R+t_2]$ 
and $(x_0,y_0) \in \cA \times [2, 3]$:
\begin{align*}
	\PR_{(x_0,y_0)} \lc (X, Y)_t\in (dx, dy)
	\pv t < T_{(\np)}\rc
	\ge c_R\; \idg{\cA}(x)\, \idg{[2, 3]}(y) \, dx\, dy.
\end{align*}
\end{lem}

\begin{lem}
	\label{lem_AA}
	Assume that $d= 1$ and that  \HDef\, and
\HgA\, hold.
	Then, for any $\ell_I>0$, there exists $c, t>0$
	and $\np \ge \ell_I$ 
	such that:
	\begin{align*}
	\frlq{(x_I, y_I)\in \cD_{\ell_I} }
	\PR_{(x_I, y_I)}( \tau_\cA \le t 
	\wedge T_{(\np)} ) \ge c.
	\EQn{O_tA2}{}
	\end{align*}
\end{lem}

\paragraph{Step 1: proof of Lemmas \ref{O_lem_cRd1}
and \ref{lem_SbA}}
Considering the calculations given in Step 4.1,
Subsection \ref{sec_lem_A},
in this case where there is no contribution from $x^{(\perp 1)}$,
we can conclude that for any $\mpr$,
there is $\np\ge \mpr$ such that $[-\mpr, 0]$ is included in $\cG_\np$. 
Adapting the reasoning given resp. in Subsections \ref{sec_lem_cRD} and \ref{sec_Sb},
we can directly conclude the proof of Lemmas \ref{O_lem_cRd1}
and \ref{lem_SbA}.

Note that the set first introduced in the proof of Lemma \ref{O_lem_Sb} here takes the form 
$[-\theta + \eta/6, -\theta + 5\eta/6]$. It is included in
$[-\mpr, 0]$ for any choice of $\mpr \ge \theta$,
so that Lemma \ref{O_lem_cRd1} can indeed replace Lemma \ref{O_lem_cRD}.
\hfill $\square$

\paragraph{Step 2: proof of Lemma \ref{lem_AA}}
Let $(x_I, y_I)\in \cD_{\ell_I}$.
\subparagraph{Case 1: $x_I\ge -\theta$}
Thanks to Lemma \ref{O_lem_v} with $u := x_I+\theta$,
there exits $t_+, c_+>0$ 
which satisfies the following property 
for any $(x_I, y_I)\in \cD_{\ell_I}$ 
such that $x_I\ge -\theta$:
\begin{align*}
	\PR_{(x_I, y_I)}( \tau_\cA \le t_+ \wedge T_{(\np)} ) \ge c_+.
\end{align*}

\subparagraph{Case 2: $x_I< -\theta$}
We recall from the proof of Lemmas \ref{O_lem_cRd1}
that there exists $\np\ge 1$ such that $[-\ell_I, 0]$
is included in $\cG_\np$.
In this set, the proof of Lemma  \ref{O_lem_A} given in Subsection \ref{sec_lemA} 
can be directly exploited to prove that there exists $t_-, c_->0$
which satisfy the following property
for any $(x_I, y_I)\in \cD_{\ell_I}$ such that $x_I\le 0$:
\begin{align*}
	\PR_{(x_I, y_I)}( \tau_\cA \le t_- \wedge T_{(\np)} ) \ge c_-.
\end{align*}

The combination of these two cases with $t:= t_+\vee t_-$ and $c:= c_+\wedge c_-$ concludes the proof of Lemma \ref{lem_AA}.
\hfill $\square$

\paragraph{Step 3: concluding the proof of  Theorem \ref{O_Mix_AF}}
By replacing Lemmas \ref{O_lem_cRD}, \ref{O_lem_Sb} and \ref{O_lem_A}
by Lemmas \ref{O_lem_cRd1}, \ref{lem_SbA} and \ref{lem_AA}
in the proof given in Subsection \ref{O_pr_MixD},
it is clear that the conclusion of Theorem \ref{O_Mix_AF}
is reached.
\hfill $\square$

\section{Absorption with failures}
\label{O_sec_AF}
\setcounter{eq}{0}

\subsection
{Proof of Theorem \ref{O_AF_D} in the case $d=1$}
\label{O_sec_AFD}

\subsubsection
{Definition of the stopping time and its elementary properties}
We consider a first process $(X, Y)$ 
with some initial condition $(x\iET,\, y\iET) \in E$.

We will prove that 
considering $\Uza = \tZa$ 
is sufficient, 
except for exceptional behavior of the process.
Given $\fl, \rho>0$, $\tZa$ shall be chosen sufficiently small to ensure 
that, with probability close to 1
 (the thresholds depending on $\fl$ and $\rho$),
 no jump has occurred before time $\tZa$,
and that the population size has not changed too much.
We define:
\begin{align*}
\dy:= \blp 3\, \ell\iET (\ell\iET+ 1)\brp^{-1},\quad
y_\wedge &:= 1/(\ell\iET +1) = 1/\ell\iET - 3\,\dy,\quad
y_\vee:= \ell\iET +1 > \ell\iET + 3\,\dy,\\
T_{\dy} &:= \inf\Lbr t\ge 0\pv |Y_t - y\iET| \ge 2\,\dy \Rbr 
< \ext.
\EQn{O_EDY}{\dy;y_\wedge;y_\vee; T_{\dy}}
\end{align*}
We recall that we can upper-bound the first jump time of $X$ by:
\begin{align*}
T\iJP:= \inf\Lbr t\ge 0\pv 
M([0, t] \times \cJ) \ge 1 \Rbr,
\EQn{O_UJ1}{T\iJP}
\end{align*}
where $\cJ$ is defined as in Subsection \ref{O_sec_Glim}.

\begin{itemize}
\item On the event $\Lbr \tZa < T_{\dy} \wedge T\iJP \wedge \ext\Rbr$, 
we set $\Uza:= \tZa$. 

\item On the event 
$\Lbr T_{\dy} \wedge T\iJP\wedge \ext  \le \tZa \Rbr$, 
we set  $\Uza:= \infty$.
\end{itemize}

Before we turn to the details 
of the proof of Theorem \ref{O_AF_D},
we first give the main scheme for proving the following lemma, 
noting that we will not go too deeply in the details of this proof.
\begin{lem}
\label{O_UAB}
We can define a stopping time $\Uza^\infty$ 
extending the above definition of $\Uza$ 
as described in Theorem \ref{O_AF_D}.
\end{lem}

\subsubsection{Step 1: main argument for the proof of Lemma \ref{O_UAB}}
\label{ApUAB}

Recall (with simplified notations)
that considering the process $(X, Y)$ 
with initial condition $(x, y)$, we define
for some $t>0$:
$\quad \Uza:= t$ on the event $\{t < T_{\delta y}\wedge T\iJP\}$,
$\quad \Uza:= \infty$ otherwise,
\begin{align*}
	\where
	& T_{\dy}:= \inf\Lbr s\ge 0\pv |Y_s - y| \ge 2\,\dy \Rbr 
	< \ext,
	&\text{ for some  $\delta y>0$, }
	\\&
	T\iJP:= \inf\Lbr s\ge 0\pv 
	M([0, s] \times \cJ) \ge 1 \Rbr,
	\\&\cJ:= \bR^d \times [0, f_{\vee}]\times [0, g_{\vee}]
	&\text{ for some } f_{\vee}, g_{\vee} >0.
\end{align*} 

Recursively, we also define:
\begin{align*}
	\tau_{E}^{i+1}:= \inf\{s\ge \tau_{E}^i +t: X_s \in E\}
	\wedge \ext,
	\AND
	\tau_{E}^0 = 0,
\end{align*}
and on the event $\Lbr \tau_{E}^{i}  < \ext\Rbr$,
for any $i$, we set:
\begin{align*}
	&T^i_{\dy} 
	:= \inf\Lbr s\ge \tau_{E}^{i} \pv |Y_s - Y(\tau_{E}^{i})| \ge 2\,\dy \Rbr,
	\\&
	U^i_j:= \inf\Lbr s\ge 0¸\pv 
	M([\tau_{E}^{i}, \tau_{E}^{i} + s] \times \cJ) \ge 1 \Rbr,
	\\&
	\Uza^\infty 
	:= \inf\{\tau_{E}^{i} + t
	\pv t\ge 0\mVg
	\tau_{E}^{i} < \infty \mVg 
	\tau_{E}^{i} + t < T^i_{\dy}  \wedge U^i_j\},
\end{align*}
where in this notation, 
the infimum equals $\infty$ if the set is empty,
$T^i_{\dy}:= \infty$ and $U^i_j =\infty$ 
on the event 
$\Lbr \ext \le \tau_{E}^{i}\Rbr$.

The proof that all these random times 
define stopping times is classical although very technical
and the reader is spared the details. 
The main point is that there is a.s. a positive gap between 
any of these iterated stopping times.
We can thus ensure recursively in $I$
that there exists a sequence of stopping times
with discrete values
$
	(\tau_{E}^{i, (n)}, T^{i, (n)}_{\dy}, U^{i, (n)}_j)
	_{\{i\le I, n\ge 1\}},$
such that
a.s. for $n$ sufficiently large, 
and $1\le i\le I$:
\begin{align*}
	&\hcm{2}
	\tau_{E}^{i}
	\le \tau_{E}^{i, (n)} 
	\le \tau_{E}^{i} + 1/n
	< \tau_{E}^{i} + t,
	\\&
	T^{i}_{\dy}
	\le T^{i, (n)}_{\dy}
	\le T^{i}_{\dy} + 1/n,\,
	\hcm{1}
	U^{i}_j
	\le U^{i, (n)}_j
 \le	U^{i}_j + 1/n.
\end{align*}

It is obvious that $\Uza^\infty$ coincide with $\Uza$ 
on the event $\Lbr \Uza \wedge \ext \le \tau_{E}^1 \Rbr$,
while the Markov property at time $\tau_{E}^1$
and the way $\Uza^\infty$ is defined
entails that on the event 
$\Lbr  \tau_{E}^1 < \Uza \wedge \ext \Rbr$,
$\Uza^\infty - \tau_{E}^1$ has indeed the same law 
as $\wtd{U}_{A}^{\infty}$ associated the process 
$(\wtd{X}, \wtd{Y})$ solution of the system \Req{O_St}
with initial condition $(X(\tau_{E}^1),\, Y(\tau_{E}^1)\,)$.
\hfill $\square$

\subsubsection
{Step 2: end of the proof of Theorem \ref{O_AF_D} when $d=1$}
Let $\ell\iET \ge 1$,
$\epsilon, \rho>0$ be prescribed.
We first require $\tZa\le 1$ to be sufficiently small.

Note that our definitions  ensure
that for any $t< \tZa \wedge T_{\dy} \wedge T\iJP$,
we have a.s.:
\begin{align*}
(X_t, Y_t) 
\in [-\ell\iET -1,\, \ell\iET]\times[y_\wedge,\, y_\vee].
\end{align*}
Thanks to Theorem \ref{O_G}, with some constant $C_G$ uniform over any $(x\iET, y\iET )\in E$:
\begin{align*}
\PR_{(x\iET,\, y\iET)}\lp  T_{\dy}  < \tZa\wedge T\iJP  \rp
&\le C_G\; 
\PR^G_{(x\iET,\, y\iET)}\lp  T_{\dy} < \tZa\wedge T\iJP  \rp 
\\&\le C_G\, 
\PR^G_{0}\lp  T_{\dy} < \tZa \rp  
\cvz{\tZa}{0},
\hcm{3}
\end{align*}
where $T_{\dy}$ under $\PR^G_{0}$ 
denotes the first time the process $|B|$ reaches $\dy$,
with $B$ a standard Brownian Motion.
Moreover:
\begin{equation*}
\PR_{(x\iET,\, y\iET)}\lp  T\iJP < \tZa\wedge T_{\dy}  \rp 
\le \PR \lp  M([0,\tZa]\times \cJ) \ge 1 \rp
\le \nu(\bR)\; f_{\vee}
\; \tZa 
\cvz{\tZa}{0}.
\end{equation*}
By choosing $\tZa$ sufficiently small, we can thus ensure the following property for any $(x\iET, y\iET )\in E$:
\begin{align*}
\PR_{(x\iET,\, y\iET)}(\Uza = \infty,  \tZa < \ext)
&\le \PR_{(x\iET,\, y\iET)}\lp  T_{\dy} < \tZa\wedge T\iJP \rp 
+ \PR_{(x\iET,\, y\iET)}\lp T\iJP < \tZa \wedge T_{\dy} \rp 
\\&\le \fl\;e^{-\rho}
\le \fl\;\exp(-\rho\, \tZa). 
\EQn{tZaDef}{}
\end{align*}

On the event $ \Lbr \tZa < T_{\dy} \wedge T\iJP \Rbr$, we have: $X_{\Uza} = x\iET- v\; \tZa$ 
and $Y_{\Uza} \in [y\iET-\dy, y\iET+\dy]$.
Indeed, 
as in the proof of Lemma \ref{O_lem_v},
we have chosen our stopping times to ensure 
that no jump for $X$ can occur before time $T\iJP \wedge \tZa\wedge T_{\dy}$.
We also rely on the Girsanov transform and Theorem \ref{O_G}
to prove that,
during the time-interval $[0, \tZa]$, 
$Y$ is indeed sufficiently diffused
(since we care now for an upper-bound, 
we can neglect the effect of assuming $\tZa < T_{\dy}$).
It makes us conclude that 
there exists $D^X > 0$ 
such that for any $x\iET \in [-\ell\iET, \ell\iET]$ and $y\iET \in [1/\ell\iET,\, \ell\iET]$:
\begin{align*}
\PR_{(x\iET,\, y\iET)}\lc
(X,\, Y) (\Uza) \in (dx,\, dy)\pv
\Uza < \ext\rc
\le D^X \,\idg{[ y\iET-2\,\dy,\,  y\iET+2\,\dy]}(y) \;
\delta_{x\iET- v\; \tZa}(dx)\;  dy.
\EQn{Umin}{}
\end{align*}

With $\alc$ the uniform distribution over $\cD_1$,  
thanks to Theorem \ref{O_Mix_DF},
there exists $c_M, t_M >0$
such that:
\begin{equation*}
	\PR_{\alc}\lc (X, Y)_{t_M} \in (dx',\, dy') \rc
	\ge c_M\; \idc{(x', y') \in \cD_{L\iET}}\; dx'\, dy'.
\end{equation*}

The idea is then to let $X$ decrease 
until it reaches $x\iET- v\, \tZa$ 
by ensuring that no jump occurs.
We then identify $u$ as the time needed for this to happen. 
Then, thanks to Theorem \ref{O_G} and Lemma \ref{O_outoftube}, 
we deduce a lower-bound on the density of $Y$ on 
$[y\iET-2\,\dy, y\iET+2\,\dy]$.
We have already proved a stronger result 
for Lemma \ref{O_lem_v}, that we let the reader adapt to obtain the following property.
For any $\tZa>0$, 
there exists $d^X_2$
which satisfies the following property
for any $x\iET \in [-\ell\iET, \ell\iET]$ and $y\iET \in [1/\ell\iET,\, \ell\iET]$.
There exists a stopping time $\UCa$  such that:
\begin{equation*}
	\PR_{\alc}\lc
	(X,\, Y) (\UCa) \in (dx,\, dy)\rc 
	\ge d^X_2 \, c_M\; \idg{[ y\iET-2\,\dy,\,  y\iET+2\,\dy]}(y) \;
	\delta_{x\iET- v\; \tZa}(dx)\;  dy.
	\EQn{VMaj}{}
\end{equation*}
The proper definition of $\UCa$
is given by $\UCa:= t_M + \tZa + (X_{t_M} - x\iET)/v\ge t_M$
on the event 
$\{X_{t_M}  \in [x\iET,\; x\iET + v\, \tZa]\}
\cap \{Y_{t_M}  \in [ y\iET-\dy/2,\,  y\iET+2\,\dy/2]\}$
(and can be made arbitrary as $t_M$ otherwise).

Thanks to Lemma \ref{O_UAB}, \Req{tZaDef}, \Req{Umin}
and \Req{VMaj},
we conclude the proof of Theorem \ref{O_AF_D},
with 
$\cp:= D^X / (d^X_2\; c_M)$.
\hfill $\square$

\subsection
{Proof of Theorem \ref{O_AF_A}}
Except that we exploit Theorem \ref{O_Mix_AF} instead of \ref{O_Mix_DF}, 
which constrains the shape of $E$,
the proof is immediately adapted from the previous Subsection \ref{O_sec_AFD}.
\hfill $\square$

\subsection
{Proof of Theorem \ref{O_AF_D} in the case $d\ge 2$}
The difficulty in this case is that, 
as long as no jump has occurred, 
$X_t$ stays confined in the line $x + \bR_+.\mathbf{e_1}$.
The "absorption" thus cannot occur before a jump.
Thus, we first wait for a jump to diffuse on $\bR^d$ 
and then let $Y$ diffuse independently 
in the same way as in Subsection \ref{O_sec_AFD}.
These two steps are summarized in the following:

\begin{prop}
\label{O_AFDX}
Given any $\rho>0$, 
$E\in \mathbf{D}$ and $\fl_X \in (0,\, 1)$,
there exists $t^X,\, c^X,\, x_\vee^X>0$ and $0< y_\wedge^X < y_\vee^X$ which satisfies the following property
for any $(x\iET,y\iET) \in E$.
There exists a stopping time $U^X$ 
such that:
\begin{align*}
&  \Lbr\ext \wedge t^X \le U^X \Rbr
= \Lbr U^X = \infty\Rbr\mVg
\quad  
   \PR_{(x\iET,y\iET)} (U^X = \infty\mVg t^X< \ext) 
     \le \fl_X\, \exp(-\rho\, t^X),
\\&\AND
    \PR_{(x\iET,y\iET)} \big(  X(U^X) \in dx
    \pv Y(U^X) \in [y_\wedge^X \mVg y_\vee^X]
     \mVg U^X <  \ext \big) 
    \le c^X \,\idg{B(0, x_\vee^X)}(x) \, dx.
      \end{align*} 
\end{prop}
We defer the proof in Subsection \ref{O_sec_AFX}.

\begin{prop}
\label{O_AFDY}
Given any $\rho,\, x_\vee^X>0$, 
$0< y_\wedge^X < y_\vee^X$ and $\fl_Y \in (0,\, 1)$,
there exists $t^Y,\, c^Y>0$ and 
$0<y_\wedge^Y  <   y_\vee^Y$
 which satisfies the following property
for any $(x,y) \in B(0, x_\vee^X) \times [y_\wedge^X \mVg y_\vee^X]$.
There exists a stopping time $T^Y$
such that:
\begin{align*}
&\hcm{2}
   \PR_{(x, y)} (T^Y\le t^Y \wedge  \ext) 
     \le \fl_Y\, \exp(-\rho\, t^Y),
\\&
\AND   \PR_{(x, y)} \big(  (X, Y)\,(t^Y) \in (dx, dy) 
    \pv t^Y < T^Y\wedge \ext \big) 
   \le c^Y \,\delta_{\Lbr x -  v\, t^Y\,\mathbf{e_1}\Rbr}(dx)\; 
   \idg{[y_\wedge^Y \mVg y_\vee^Y]}(y) \, dy.
\end{align*} 
\end{prop}

The proof of Lemma \ref{O_AFDY}
is taken mutatis mutandis from the one in Subsection \ref{O_sec_AFD}.
It leads to define $\Uza$ as below.

   
$\bullet$ $\Uza:= U^X + t^Y$ 
on the event  $\Lbr U^X < t^X\wedge \ext\Rbr\cap \{t^Y < \wtd{\ext}\wedge \wtd{T}^Y\}$, 
   where $\wtd{\ext}$ and $\wtd{T}^Y$
   are defined as respectively $\ext$ and $T^Y$ 
   for the solution $(\wtd{X}_t, \wtd{Y}_t)$,
   defined on the event $\Lbr U^X < t^X\wedge \ext\Rbr$, of:
   \begin{align*}
   &\left\{
   \begin{aligned}
   \wtd{X}_t
   &= X(U^X) - v\,t\, \mathbf{e_1} 
   + \int_{[U^X,U^X+t] \times \bR^d \times (\bR_+)^2 } 
   w \;    \varphi \left( \wtd{X}_{s^-},\,\wtd{Y}_s,\,w,\,u_f,\, u_g \right)\; \Md\\ 
   \wtd{Y}_t
   &= Y(U^X) 
   + \int_{0}^{t} \psi\lp \wtd{X}_s,\, \wtd{Y}_s\rp ds 
   + \int_{U^X}^{U^X + t}dB_r.
   \end{aligned}
   \right.
   \end{align*}
   
$\bullet$ Else $\Uza:= \infty$.


\begin{fact}
\label{O_AFDst}
There exists a stopping time $\Uza^\infty$ 
extending the above definition of $\Uza$ 
as described in Theorem \ref{O_AF_D}
(with $t = t^X + t^Y$ here).
\end{fact}

The proof of Fact \ref{O_AFDst} is technical 
but classical from the way
 we define $U^X$ and $T^Y$
 and similar to the proof of Lemma \ref{O_UAB}.
 The reader is spared this proof.

\subsubsection
{Proof of Theorem \ref{O_AF_D}
as a consequence of Propositions \ref{O_AFDX}-2 and Fact  \ref{O_AFDst}}

Given  $\rho>0,$ $\fl\in (0,1)$ and some $E\in \mathbf{D}$, we define 
$\fl_X:= \fl/4$ and deduce from Proposition \ref{O_AFDX} 
the values $t^X$, $c^X$, $x_\vee^X$, $y_\wedge^X$, $y_\vee^X$ 
and the definition for the stopping times 
$U^X$ with the associated properties.

With $\fl_Y:= \fl \, \exp(-\rho\, t^X)/2$, 
we then deduce from Proposition \ref{O_AFDY} 
the values $t^Y$, $c^Y$, $y_\wedge^Y$, $y_\vee^Y$ 
and the stopping time $T^Y$ with the associated properties.
Defining, for some $(x, y) \in E$, $\Uza$ as in Fact \ref{O_AFDst} 
 and combining these results:
\begin{align*}
&
 \Lbr\ext \wedge (t^X+t^Y) \le \Uza \Rbr
= \Lbr \Uza = \infty\Rbr\mVg
\EQn{O_Uza}{}
\\&
   \PR_{(x, y)} \lc
    (X, Y)\,(\Uza) \in (dx, dy) 
    \pv \Uza <  \ext \rc
 \\&\hcm{3}
   \le c^X\, c^Y \, \idg{B(0, x_\vee^X + v\, t^Y)}(x) \, 
      \idg{[y_\wedge^Y \mVg y_\vee^Y]}(y)\; dx \, dy,
      \EQn{O_UY}{}
\end{align*}
\begin{equation*}
   \PR_{(x, y)} \big(\Uza = \infty, t^X+t^Y < \ext \big) 
   \le \fl_X\, \exp(-\rho\, t^X) + \fl_Y\, \exp(-\rho\, t^Y) 
   \le \fl\, \exp(-\rho\, [t^X+ t^Y]),
   \EQn{O_flX}{}
\end{equation*}
where we exploited the definitions of $\fl_X$, $\fl_Y$ and 
that $t^Y \le \ln(2)/\rho$ (i.e. $1/2 \le \exp(-\rho\, t^Y)$)
  in the last inequality.

For the opposite upper-bound, 
we recall first that $\alc$ is chosen to be uniform 
over the compact space $\Delta$,
that is included in some $\cD_\ell$.
Exploiting Theorem \ref{O_Mix_AF} on this set $\cD_\ell$, 
we deduce that
there exists $t, c>0$ such that:
\begin{align*}
\PR_{\alc} \big[  (X, Y)\,(t) \in (dx, dy) 
\pv t < \ext \big]
\ge c\; \idg{B(0, x_\vee^X + v\, t^Y)}(x) \, 
   \idg{[y_\wedge^Y \mVg y_\vee^Y]}(y)\; dx \, dy.
\EQn{O_tal}{}
\end{align*}

Combining \Req{O_Uza}--\Req{O_tal}\,
 ends the proof of Theorem \ref{O_AF_D} in the case $d \ge 2$.
\hfill $\square$

\subsubsection
{Proof of Proposition \ref{O_AFDX}}
\label{O_sec_AFX}

For readability,
note that most of the subscripts $"X"$ 
(except for $t^X$)
from  Proposition \ref{O_AFDX}
are removed in this proof.
\\

First, remark that without any jump, 
$\|X\|$ tends to  infinity, 
which makes the population almost doomed to extinction.
We can thus find some time-limit $t_\vee$ such that,
even with an amplification 
of order $\exp(\rho\, t_\vee)$,
the event that the 
population survived 
without any mutation occurring in the time-interval $[0, t_\vee]$
is exceptional enough.
With this time-scale, we can find an upper-bound $y_\vee$ on $Y$:
that the population reaches such size before $t_\vee$ 
is an exceptional enough event.
For the lower-bound, 
we exploit the fact that extinction is very strong 
when the population size is too small.
Thus, that the population has survived --at least for a bit--
after declining below this lower-bound $y_\wedge$
is also an exceptional enough event.

The last part is to ensure that this first jump is indeed diffuse in $X$
(which is why we need $\nu(dw)$ to have a density w.r.t. Lebesgue
 with the bound of \RHyp{Hac}).
 \\

For $y_\vee > \ell\iET > 1/\ell\iET > y_\wedge > 0$, 
$t_{\vee}, w_\vee >0$ 
 and initial condition $(x, y) \in E$, 
 let:
\begin{align*}
&\hcm{3}
T\iJP:= \inf\Lbr t\ge 0\pv \Delta X_t \neq 0 \Rbr,
\EQn{O_UJF}{T\iJP}
\\&
T_Y^\vee 
:= \inf\Lbr t\ge 0\pv Y_t = y_\vee \Rbr,\quad
T_Y^\wedge
:= \inf\Lbr t\ge 0\pv Y_t = y_\wedge \Rbr < \ext.
\EQn{O_TYF}{T_Y^\vee;T_Y^\wedge}
\end{align*}
On the event $\{ T\iJP < t_{\vee}  \wedge T_Y^\vee \wedge T_Y^\wedge\}
\cap\{\|\Delta X_{T\iJP}\| < w_\vee\}$,
we define $ U:= T\iJP$.
Else $U:= \infty$.

To choose $y_\wedge, y_\vee$, $t_\vee$ and $w_\vee$, we refer to the following lemmas,
which are treated as the four first steps of the proof:

\begin{fact}
\label{O_tAF}
For any $\rho, \fl_{1}>0$, there exists $t_{\vee}>0$
such that:
\begin{align*}
 \frlq{(x, y) \in E}
\PR_{(x, y)} (t_{\vee} < T\iJP \wedge \ext ) \le \fl_{1}\, \exp(-\rho\, t_{\vee}).
\end{align*}
\end{fact}

\begin{fact}
\label{O_yvA}
For any $t_{\vee}, \fl_{2}>0$, there exists $y_\vee>0$
such that:
\begin{align*}
\frlq{(x, y) \in E}
\PR_{(x, y)} (T_Y^\vee  < t_{\vee}\wedge \ext ) 
\le \fl_{2}.
\end{align*}
\end{fact}

\begin{fact}
\label{O_ywA}
For any  $t_S, \fl_{3}>0$, there exists $y_\wedge>0$
such that:
\begin{align*}
\frlq{x \in \bR^d}
\PR_{(x, y_\wedge)} (t_S < \ext ) 
\le \fl_{3}.
\end{align*}
\end{fact}

\begin{fact}
	\label{O_wv}
For any  $t_\vee, \fl_{4}>0$, there exists $w_\vee>0$
	such that:
	\begin{align*}
	\frlq{(x, y) \in E}
	\PR_{(x, y)} (\|\Delta X_{T\iJP}\| \ge w_\vee\mVg T\iJP < t_\vee\wedge \ext) 
	\le \fl_{4}.
	\end{align*}
\end{fact}

\begin{fact}
\label{O_UjDif}
For any $t_{\vee}>0$,
and any $y_\vee > \ell\iET > 1/\ell\iET > y_\wedge > 0$,
 there exists $c, x_\vee>0$
such that:
\begin{align*}
\frlq{(x, y) \in E}
\PR_{(x, y)} \big( 
 X(U) \in dx \pv
U < \ext \big) 
\le c \,\idg{B(0, x_\vee)}(x) \, dx.
\end{align*}
\end{fact}

\paragraph
{Step 1: proof of Fact \ref{O_tAF}}
Exploiting assumption \RHyp{HrN},
as long as $\|X\|$ is sufficiently large,
we can ensure that the growth rate of $Y$ is largely negative,
leading to a quick extinction.
The proof is similar to the one of Lemma 3.2.2
in \cite{AV_QSD}, where more details can be found.
We consider the autonomous process $Y^D$ 
as an upper-bound of $Y$ where the growth rate is replaced by $r_D$.
For any $t_D$ and $\rho$,
there exists $r_D$ (a priori negative)
such that whatever $y_D$ the initial condition of $Y^D$,
survival of $Y^D$ until $t_D$ (i.e. $t_D < \ext^D$)
happens with a probability smaller than $\exp(-2\, \rho\, t_D)$.
Thanks to Assumption \RHyp{HrN}, 
we define $x_\vee$ such that for any $x$, 
$\|x\|\ge x_\vee$ implies $r(x)\le r_D$.
We then deduce:
\begin{align*}
\frl{(x, y)}\quad
\PR_{(x,y)} (\frl{t\le t_D} \|X_t\| \ge x_\vee \pv t_D < \ext)
\le \sup_{ y_{D}>0} \,\PR_{y_{D}} \lp  t_D < \ext^D \rp
\le \exp(-2\, \rho\, t_D).
\end{align*}

Let $t\iET:= (x_\vee + \ell\iET)/v$
and assume $t_{\vee} \ge t\iET$.
A.s. on $\Lbr t_{\vee} < T\iJP \wedge \ext\Rbr$ for any $(x, y)\in E$:
\begin{align*}
\frlq{t\iET \le t\le t_{\vee}}
 \|X(t)\| = \|x - v\, t\, \mathbf{e_1}\|
\ge x_\vee.
\end{align*}
Exploiting inductively the Markov property 
at times $t_{\vee}:= t\iET + k\, t_D$ for $k\ge 1$, 
we obtain:
\begin{align*}
\frl{(x, y)}\quad
\exp[\rho\, t_{\vee}]\; \PR_{(x,y)} (t_{\vee} < T\iJP \wedge \ext)
\le \exp(\rho\, [t\iET-k\,  t_D]) 
\cvifty{k} 0.
\end{align*}

\paragraph
{Step 2: proof of Fact \ref{O_yvA}} 
This is an immediate consequence of the fact that $Y$ is
 upper-bounded by the process $Y^{\vee}$ 
 given in \Req{O_Y+} with initial condition $\ell_M$.
 This bound is uniform in the dynamics of $X_t$ and $M$
 and uniform for any $(x, y) \in E$.
It is classical that a.s. $\sup_{t\le t_{\vee}} Y^{\vee}_t < \infty$, 
which proves the Lemma,
see e.g. Lemma 3.3 in \cite{BM15}.

\paragraph
{Step 3: proof of Fact \ref{O_ywA}} 
Like in the proof of Proposition 4.2.3 in \cite{AV_QSD},
cf Appendix D,
we exploit $r_\vee$ as the upper-bound 
of the growth rate of the individuals
to relate to the formulas 
for Continuous State Branching Processes.
Referring for instance to \cite{P16} Subsection 4.2, notably Lemma 5, 
it is classical that $0$ 
is an absorbing boundary for these processes 
(we even have explicit formulas for the probability of extinction).
This directly entails the result 
of the present lemma
that the probability of extinction 
tends uniformly to zero 
as the initial population size tends to  zero.

\paragraph
{Step 4: proof of Fact \ref{O_wv}} 
On the event $\{T_J< t_\vee \wedge \ext\}$,
for any initial condition $(x, y)\in E$, 
there exists a compact $K$ of $\bR^d$
that contains $X_t= x-v\,t$ for any $t\in [0, T_J)$.
Thanks to Assumption $[H2]$, 
there exists an upper-bound $g_\vee$
of $g$ valid on $K\ltm \bR^d$.

Let $\eps_4>0$ and $\rho_W:= (-1/t_\vee) .  \log(1-\eps_4)$.
 We define $w_\vee$ such that $\nu(B(0, w_\vee)^c)\le \rho_W/g_\vee$.
 Then we can couple the process $X$ 
 to an exponential r.v. $T_W$ of mean $1/\rho_W$
 such that 
on the event $\{T_J< t_\vee \wedge \ext\}\cap \{\|\Delta X_{T\iJP}\| \ge w_\vee\}$,
$T_J \le T_W$ holds a.s.
We can conclude with  the following upper-bound:
\begin{align*}
	\frlq{(x, y) \in E}
\PR_{(x, y)} (\|\Delta X_{T\iJP}\| \ge w_\vee\mVg
 T\iJP < t_\vee\wedge \ext) 
\le \PR(T_W < t_\vee)= 1 - \exp(-\rho_W\, t_\vee) \le \fl_4.
\end{align*}

Note that under Assumption $(A)$,
the jump at time $T_J$ cannot make the process escape $K$.
This provides a deterministic upper-bound $w_\vee$
such that $\|\Delta X_{T\iJP}\| \ge w_\vee$ 
a.s. on $\{T\iJP < t_\vee\wedge \ext\}$.

\paragraph
{Step 5: proof of Fact \ref{O_UjDif}} 
For $x_\vee:= \ell\iET + v\, t_{\vee}$, let:
\begin{align*}
\cp
:= \sup \Lbr 
\dfrac{g(x, w)\, \nu(w)}{
\int_{\bR^d} g(x, w')\, \nu(w')\, dw'}	\pv
\|x\| \le x_\vee\mVg
 w\in \bR^d \Rbr < \infty.
 \EQn{O_dX}{\cp}
\end{align*}
\AP{APfuj}
We exploit a sigma-field $\cF^*_{T\iJP}$
that includes the whole knowledge of the process
until time $T\iJP$,
except for the size of the jump at this time.
It is rigorously defined and studied in Appendix \Alph{APfuj}.
Conditionally on $\cF^*_{T\iJP}$
on the event 
$\Lbr U < \ext \Rbr
\in \cF^*_{T\iJP}$,
the law of $X(T\iJP)$ is given by:
\begin{align*}
\dfrac{g(X[T\iJP-]\mVg x - X[T\iJP-])\; \nu(x - X[T\iJP-])}{
\int_{\bR^d} g(X[T\iJP-], w')\; \nu(w')\, dw'}
\; dx.
\end{align*}
Note also that a.s. $\|X[T\iJP-]\| \le \ell\iET + v\, t_{\vee} = x_\vee$ 
(since no jump has occurred yet).

Since $\|\Delta X_{T\iJP} \| \le w_\vee$ 
on the event $\Lbr U < \ext \Rbr$,
with $\bar{x}_\vee:= x_\vee + w_\vee$, 
we get the following upper-bound 
of the law of $X(T\iJP)$:
%
%
%
%
\begin{align*}
&\PR_{(x, y)} \lp
 X(U) \in dx \pv
U < \ext \rp 
= \PR_{(x, y)} \lp
\E\lc X(U) \in dx \bv 
 \cF^*_{T\iJP} \rc  \pv
U <  \ext \rp
\\&\hcm{4}
\le \cp \,\idg{B(0, \bar{x}_\vee)}(x) \, dx.
\end{align*}

\paragraph
{Step 6: concluding the proof of Proposition \ref{O_AFDX}}
\textcolor{white}{:}

Let $\ell\iET, \rho, \fl >0$.
We first deduce $t_{\vee}$ thanks to  Fact \ref{O_tAF}  such that:
\begin{align*}
\frlq{(x, y) \in E}
\PR_{(x, y)} (t_{\vee} < T\iJP \wedge \ext ) 
\le \fl\, \exp(-\rho\, t_{\vee})/8.
\EQn{O_tv}{t_{\vee}}
\end{align*}

Thanks to Fact \ref{O_yvA},
we deduce some $y_\vee>0$ such that:
\begin{equation}
\frlq{(x, y) \in E}
\PR_{(x, y)} (T_Y^\vee  < t_{\vee}\wedge \ext ) 
\le \fl\, \exp(-\rho\, t_{\vee})/8.
\EQn{O_EyvX}{y_\vee}
\end{equation}
We could take any value for $t_S$ 
(so possibly $1$),
yet $\, t_S = \log(2) / \rho\,$ seems somewhat more practical.
We then deduce $y_\wedge$
thanks to Fact \ref{O_ywA} such that:
\begin{equation*}
\Tsup{x\in \bR^d} 
\PR_{(x, y_\wedge)}
(t_S < \ext)
\le \fl\, \exp(-\rho\, t_{\vee})/8.
\end{equation*}
This implies that for any $(x, y) \in E$:
\begin{align*}
&\PR_{(x, y)} \lp 
t_\vee + t_S < \ext \mVg
T_Y^\wedge < t_\vee\wedge \ext 
\wedge T_Y^\vee \wedge T_J\rp 
\\&\hcm{1}
\le\E_{(x, y)} \lp
\PR_{(X_{T_Y^\wedge}, y_\wedge)}
(t_S < \ext)
\pv T_Y^\wedge < t_\vee\wedge \ext 
\wedge T_Y^\vee \wedge T_J\rp 
\\&\hcm{1}
\le \fl\, \exp(-\rho\, t_{\vee})/8.
\EQn{O_EywX}{y_\wedge}
\end{align*}
$w_\vee$ is chosen thanks to Fact \ref{O_wv} such that:
\begin{align*}
\frlq{(x, y) \in E} 
\PR_{(x, y)} \lp \|\Delta X_{T\iJP}\| \ge w_\vee
\mVg T_J < t_\vee\wedge \ext\rp
\EQn{O_Ewv}{}
\end{align*}

Thanks to Fact \ref{O_UjDif},
there exist $\cp, x_\vee>0$  such that:
\begin{align*}
\frlq{(x, y) \in E}
\PR_{(x, y)} \big( 
X(U) \in dx \pv
U < \ext \big) 
\le \cp \,\idg{B(0, x_\vee)}(x) \, dx.
\end{align*}
Thanks to theconstruction of $U$,
and noting that $t^X:= t_\vee + t_S$,
it is clear that $U 
\ge \ext\wedge t^X$ 
is equivalent to $U = \infty$.
Combining \Req{O_tv}, \Req{O_EyvX}, \Req{O_EywX}
and \Req{O_Ewv}:
\begin{align*}
&
\PR_{(x, y)} (U = \infty
\mVg t_\vee + t_S< \ext) 
\le \PR_{(x, y)} (t_{\vee} < T\iJP \wedge \ext ) 
+ \PR_{(x, y)} (T_Y^\vee  < t_{\vee}\wedge \ext ) 
\\&\hcm{1.5}
+\PR_{(x, y)} \lp \|\Delta X_{T\iJP}\| \ge w_\vee
\mVg T_J < t_\vee\wedge \ext\rp
+ \PR_{(x, y)} \lp 
t_\vee + t_S < \ext \mVg
T_Y^\wedge < t_\vee\wedge \ext 
\wedge T_Y^\vee \wedge T_J\rp 
\\&\hcm{5}
\le \fl\, \exp(-\rho\, t_{\vee})/2
= \fl\, \exp(-\rho\, t^X).
\end{align*}
This ends the proof of Proposition \ref{O_AFDX}.\hfill $\square$
\\

The proof of Theorem \ref{O_AF_D} in the case $d\ge 2$ is now completed.
All the theorems have been proved at this point.
\AP{APSim} 
There are two appendix,
the first one being devoted 
to the filtration $\cF^*_{T\iJP}$ up to the jumping time.
We finish in Appendix \Alph{APSim}
with first results of simulations 
that shall help illustrate the discussion given in Subsection \ref{O_sec_ecoevo}.


%


\section*
{Appendix \Alph{APfuj}: A specific filtration for jumps}
\setcounter{eq}{0}
\stepcounter{section}

This appendix extends to our case the result already presented 
in \cite{AV_disc}: 
there exists a sigma-field $\cF^*_{T\iJP}$ 
which informally 
``includes the information carried by $M$ and $B$''
 until the jump time $T_J$
except the realization of the jump itself.

Denote $W$ as the additive effect
on $X$ of the first jump of $X$,
occurring at time $T\iJP$. We then define:
$$\cF^*_{T\iJP}:= \sigma\lp A_s \cap\Lbr s< T\iJP\Rbr\pv s >0,\, 
A_s \in \cF_s\rp.$$
\textbf{Properties of $\cF^*_{T\iJP}$:}
If $Z_s$ is $\cF_s$-measurable and $s<t\in (0, \infty]$, 
$Z_s\, \idc{s < T\iJP \le t}$ is  $\cF^*_{T\iJP}$-measurable.

\textbf{Lemma (\Alph{APfuj}1).}
\textsl{
For any left-continuous and adapted process $Z$, $Z_{T\iJP}$ is $\cF^*_{T\iJP}$-measurable.
Reciprocally, 
$\cF^*_{T\iJP}$ is in fact the smallest $\sigma$-algebra 
generated by these random variables.\\
In particular, for any stopping time $T$, $\Lbr T\iJP \le T\Rbr \in \cF^*_{T\iJP}$.}
\\


\textbf{Lemma (\Alph{APfuj}2).}
\textsl{
 For any $h: \bR \rightarrow \bR_+$ measurable, $(x,y) \in (-L, L) \times \bR_+$:}
 $$ \E_{(x,y)} \lc h(W) \Bv \cF^*_{T\iJP}\rc 
 = \dfrac{\intR h(w)\, f(Y_{T\iJP}) g(X_{T\iJP-}, w) \; \nu(dw)}
  {\intR f(Y_{T\iJP}) g(X_{T\iJP-}, w') \; \nu(dw')}.$$

\subsection*
{Proof of Lemma (\Alph{APfuj}1):} 

For any left-continuous and adapted process $Z$,
$\quad
Z_{T_J} = \limInf{n} \sum_{k \le n^2} Z_{\frac{k-1}{n}} \;
\idc{\frac{k-1}{n}< T_J \le \frac{k}{n}},
$
\\
where by previous property and the fact that $Z$ is adapted:
$Z_{\frac{k-1}{n}} \;
\idc{\frac{k-1}{n}< T_J \le \frac{k}{n}}$
is $\cF^*_{T_J}$-measurable for any $k, n$.
Reciprocally, for any $s>0$ and $A_s \in \cF_s$:
\begin{align*}
\idg{A_s\cup\{s< T_J\} } 
= \lim_{n\ge 1} Z^n_{T_J},
\where Z^n_t:= \{1\wedge [n\, (t-s)_+]\}\ltm \idg{A_s}.
\end{align*}
Now, for any stopping time $T$,
and any $t \ge 0$,  $\Lbr t \le T \Rbr \in \cF_t$
and $ \Lbr t \le T \Rbr
= \underset{s<t}{\cap} \Lbr s \le T  \Rbr$,
thus  $\Lbr T_J \le T \Rbr\cap\Lbr T_J < \infty\Rbr \in \cF^*_{T_J}$.
Similarly:
\begin{align*}
\Lbr T_J = T = \infty \Rbr
= \underset{s>0}{\cap} \Lbr s < T \Rbr \cap \Lbr s< T_J \le \infty\Rbr \in \cF^*_{T_J}.
\end{align*}

\subsection*
{Proof of Lemma (\Alph{APfuj}2):} 

Let:
$$Z_t:= \dfrac{\intR h(w')\, f(Y_t) g(X_{t-}, w') \; \nu(dw')}
  {\intR f(Y_{t}) g(X_{t-}, w'') \; \nu(dw'')},$$
which is a left-continuous and adapted process.
Thanks to Lemma (\Alph{APfuj}1),
$Z_{T\iJP}$ 
    is $\cF^*_{T\iJP}$-measurable.

We note the two following identities:
\begin{align*}
&h(W) = \intT h(w) \;\idc{t = T\iJP} \; M(dt,\, dw,\, du_f,\, du_g)
\\&\dfrac{\intR h(w)\, f(Y_{T\iJP}) g(X_{T\iJP-}, w) \; \nu(dw)}
  {\intR f(Y_{T\iJP}) g(X_{T\iJP-}, w') \; \nu(dw')}
\\ &\hspace{0.5cm}
= \intT \dfrac{\intR h(w')\, f(Y_t) g(X_{t-}, w') \; \nu(dw')}
  {\intR f(Y_{t}) g(X_{t-}, w'') \; \nu(dw'')}\;
\idc{t = T\iJP} \; M(dt,\, dw,\, du_f,\, du_g),
\end{align*}

Then, we exploit Palm's formula
 to prove that their product with any
$Z_s\; \idc{s < T\iJP \le r}$
has the same average
for any $s<r$ and $Z_s$ $\cF_s$-measurable:
\begin{align*}
&\E_{(x, y)} \lc\, h(W) \, Z_s \pv s < T\iJP \le r\,\rc\\
 &\hspace{0.5cm}
= \E_{(x, y)} \lc Z_s\, 
\intT h(w) \;
\idc{t = T\iJP} \; M(dt,\, dw,\, du_f,\, du_g)
\pv s < T\iJP \le r \rc\\
 &\hspace{0.5cm}
= \E_{(x, y)} \lc 
\intT Z_s\, h(w) \;\idg{(s,\, r]}(t)\, 
\idc{t = T\iJP} \; M(dt,\, dw,\, du_f,\, du_g) \rc\\
 &\hspace{0.5cm}
= \E_{(x, y)} \lc 
\intT \idg{(s,\, r]}(t)\, Z_s\, h(w) \; 
\idc{t = \widehat{T\iJP}} \; dt\, \nu(dw)\, du_f\, du_g \rc,
\end{align*}
where, according to Palm's formula, $\widehat{T\iJP}$ is the first jump of the process $(\widehat{X},\, \widehat{Y})$ encoded by $M + \delta_{(t,w,u)}$ and $B$
(cf e.g. \cite{DV08} Proposition 13.1.VII).
 Since $(\widehat{X},\, \widehat{Y})$ coincide with $(X,\,Y)$ at least up to time $t>s$, $Z_s$ was not affected by this change. Moreover:
 $$\Lbr t = \widehat{T\iJP} \Rbr 
 = \Lbr t\le T\iJP \Rbr \cap \Lbr u \le f(Y_t)\, g(X_{t-}, w)\Rbr. $$
Thus:
\begin{align*}
&\E_{(x, y)} \lc\, h(W) \, Z_s \pv s < T\iJP \le r\,\rc
\\&
= \E_{(x, y)} \lc 
\intT \idg{(s,\, r]}(t)\, Z_s\, h(w) \; 
\idc{u_f \le f(Y_t)}\,
\idc{u_g \le g(X_{t-}, w)}\;
\idc{t \le T\iJP} \; dt\, \nu(dw)\, du_f\, du_g \rc,
\\&\hspace{0.5cm}
= \E_{(x, y)} \lc Z_s\,
\int_s^r \int_{\bR} 
\idc{t \le T\iJP}\; h(w) \; f(Y_t)\, g(X_{t-}, w)\;
 \nu(dw)\, dt \rc.
\end{align*}

On the other hand, and with the same spirit:
\begin{align*}
&\E_{(x, y)} \lc \dfrac{\intR h(w')\, f(Y_{T\iJP}) g(X_{T\iJP-}, w') \; \nu(dw')}
  {\intR f(Y_{T\iJP}) g(X_{T\iJP-}, w'') \; \nu(dw'')}
 \, Z_s \pv s < T\iJP \le r\,\rc\\
 &\hspace{0.5cm}
= \E_{(x, y)} \Big[ Z_s\, 
\intT \dfrac{\intR h(w')\, f(Y_t) g(X_{t-}, w') \; \nu(dw')}
  {\intR f(Y_{t}) g(X_{t-}, w'') \; \nu(dw'')}\\
  &\hspace{7.5cm} 
\times \idc{t = T\iJP} \; M(dt,\, dw,\,  du_f,\, du_g) 
\pv s< T\iJP \le r\,\Big]\\
 &\hspace{0.5cm}
= \E_{(x, y)} \Big[  
\intT Z_s\, \idg{(s, r]}(t)\, 
\dfrac{\intR h(w')\, f(Y_t) g(X_{t-}, w') \; \nu(dw')}
  {\intR f(Y_{t}) g(X_{t-}, w'') \; \nu(dw'')}\\
   &\hspace{7.5cm}  
\times \idc{t = T\iJP} \; M(dt,\, dw,\,  du_f,\, du_g) \Big]\\
 &\hspace{0.5cm}
= \E_{(x, y)} \Big[  
\intT Z_s\, \idg{(s, r]}(t)\, 
\dfrac{\intR h(w')\, f(Y_t) g(X_{t-}, w') \; \nu(dw')}
  {\intR f(Y_{t}) g(X_{t-}, w'') \; \nu(dw'')}\\
   &\hspace{6.5cm}  
\times \idc{t \le T\iJP} \;
\idc{u_f \le f(Y_t)}\,
\idc{u_g \le g(X_{t-}, w)}\;
dt\, \nu(dw)\,  du_f\, du_g \Big]\\
&\hspace{0.5cm}
= \E_{(x, y)} \lc Z_s\,
\int_s^r \int_{\bR} 
\idc{t \le T\iJP}\; h(w') \; f(Y_t)\, g(X_{t-}, w')\;
 \nu(dw')\, dt \rc,
\end{align*}
which is indeed the same integral as for $h(W)$. 
\hfill $\square$

\section*
{Appendix \Alph{APSim}: Brief overview of  characteristic profiles of the quasi-stationary regime
obtained by simulations}

We provide in this Appendix \Alph{APSim} some results of a particular choice 
of three parameters regime whose comparison shall shed light on the discussion 
given in Subsection \ref{O_sec_ecoevo}.
We present the profiles of the characteristic distributions and functions 
of the quasi-stationary regime, namely the QSD, the quasi-ergodic distribution (QED) and the survival capacity
(the limiting properties are recalled just beside the figures).

The details of the exploited parameters are as follows. 
For population size dynamics, the growth rate as a function of $x$ is here chosen to be of the form $r(x) = 4 - 30\ltm |x|$.
A parabolic profile would give very similar results.
The competition rate is $c = 0.1$, 
which leads to population sizes at quasi-equilibrium (carrying capacity) close to 40
(in arbitrary units).
$2$ and $6$ are respectively the values
for the 
the diffusion coefficient $\sigma$ and the speed of the environment $v$.
Thus, there are rapid fluctuations in population size in the time-scale where adaptation changes.

The profile of additive effects of mutations 
is given by $\nu(dw) = \frac{1}{2 w_0} \exp(-|w|/w_0)$.
It is therefore symmetrical exponential, with $w_0 = 0.03$, so with many small mutations.
The effect of population size on the fixation rate is simply proportional $f_N(n) = m\ltm n$.
The mutation rate $m$ is the only parameter modified here between the 3 simulation sets:
it takes the values $m = 0.85$, $m = 0.55$ and $m = 0.25$.
The choice of these values is done so that 
the adaptation is critical at $m = 0.55$:
for larger $m$ like $m = 0.85$, 
 extinction is kept almost negligible, so that we say that adaptation is spontaneous;
 whereas for smaller values of $m$, extinction plays a consequent role 
 and produces differences in shape between the QSD and the QED.

We exploited the following expression 
for the probability of invasion:
$$g(x, w):= \dfrac{ N_H(x)\ltm \Delta r / \sigma}{1- \exp[-N_H(x)\ltm \Delta r / \sigma]},$$
where $\Delta r:= r(x+w) - r(x)$ is the variation of the growth rate between the mutant and the resident,
and $N_H(x)$ is the harmonic mean of a resident population  with fixed trait $x$
 (averaged against its associated QSD).
  Deleterious mutations are allowed, 
  but their probability of fixation is very reduced 
  if they are strong in relation to population fluctuations.
 The values of $N_H$ are estimated numerically, with the following profile:
 \begin{figure}[h]
 	\begin{center}
 		\begin{minipage}[c]{.46\linewidth}
 		\includegraphics[width=70mm, height = 40mm]{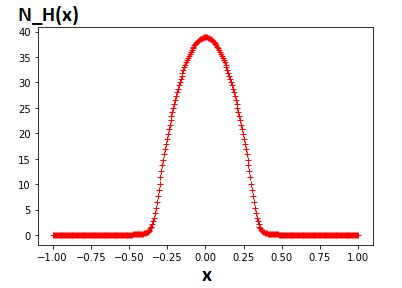}
 		\end{minipage}
 		 		\begin{minipage}[c]{.46\linewidth}
 			\includegraphics[width=70mm, height = 40mm]{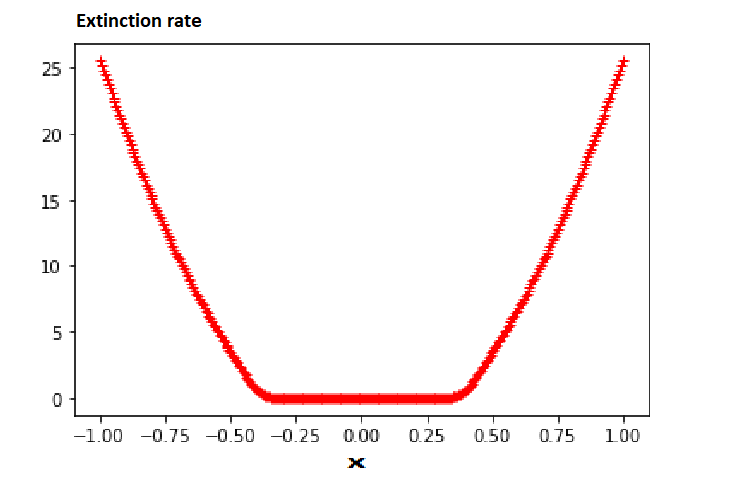}
 		\end{minipage}
 		\caption{$N_H(x)$ on the left side: the harmonic means of the population size fluctuations
 			of the process $(\wtd{N}^x_t)_{t\ge 0}$ with fixed traits $x$
 			 (given by the associated QSD);
 			 the extinction rate of these QSD is on the right.}
 		\label{MHarm}
 	\end{center}
 \end{figure}

The formula relies on the Kimura diffusion approximation
that has been derived 
in the case of fixed population size.
Assuming rapid size fluctuations, 
we choose the harmonic mean as the reference
by referring to classical approximations obtained in the case 
of periodically fluctuating population sizes
(cf notably \cite{OW97}).
More details are given (in French) in my PhD manuscript,
and a subsequent paper is planned to discuss these results 
and the relevance of this estimation.
The comparison of such a two-component stochastic model to the individual-based model
through the QSD and QED
shall be a good test for the relevance of such formula.
The kind of dependence in the difference in growth rate
seems to play a crucial role for having a QED as much conserved.
\\

These simulations were obtained by calculating the evolution of the densities themselves.
This method is related to those of finite volumes, with an explicit numerical scheme 
and a renormalization of density estimates at each time step.
The transitions to $X$ and $N$ are performed successively to reduce the calculation time.

\begin{figure}
	\textbf{Profiles of the "QSD"}
	
	\begin{minipage}[c]{.46\linewidth}
		\includegraphics[width=70mm, height = 36mm]{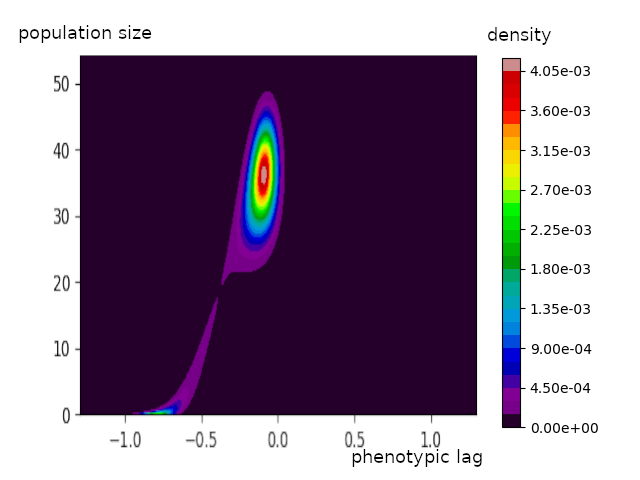}
		\label{PZ60}
	Case of a spontaneous adaptation
	\end{minipage}
	\begin{minipage}[c]{.46\linewidth}
		$\PR_{(x,\,n)} \lc (X, N)_t \in (dx, dn)\Bv  t < \tau_\partial\rc$
		
		\hspace*{1cm}$\cvifty{t}\alpha(dx, dn)$
	\end{minipage}
	
	\begin{minipage}[c]{.46\linewidth}
		\includegraphics[width=70mm, height = 36mm]{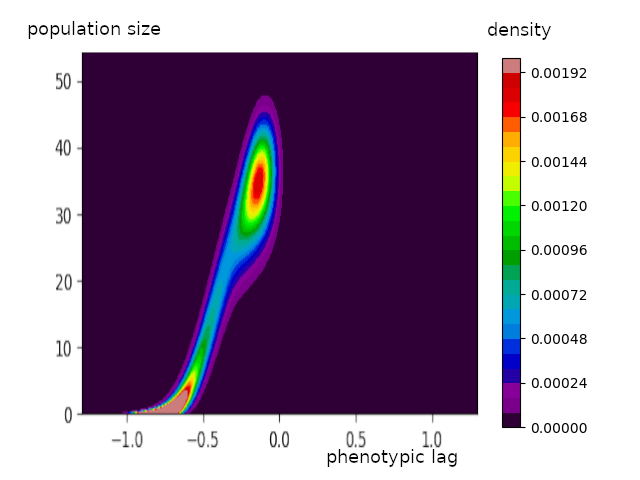}
		\label{PZ61}
	Critical regime of adaptation
	\end{minipage}
	\begin{minipage}[c]{.46\linewidth}
		\includegraphics[width=70mm, height = 36mm]{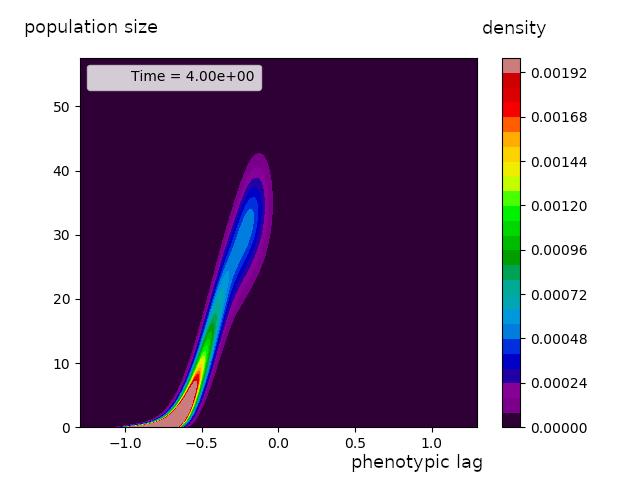}
		\label{PZ62}
	\mbox{Selection through the extinction of populations}
	\end{minipage}

	\textbf{Profiles of the QED: the invariant measure of the Q-process}
	
	\begin{minipage}[c]{.46\linewidth}
		\includegraphics[width=70mm, height = 40mm]{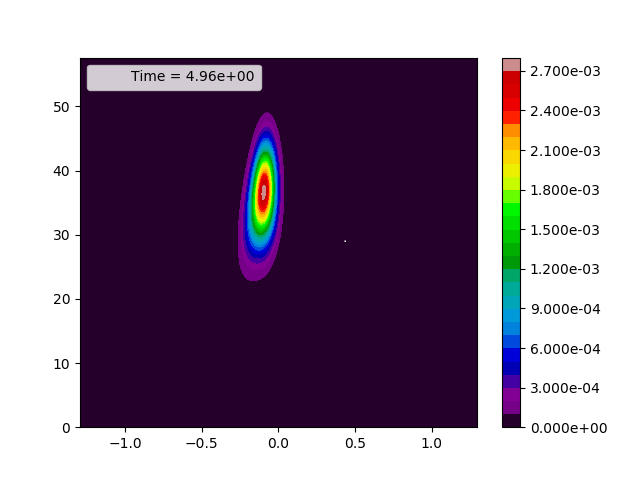}
		
		Case of a spontaneous adaptation
	\end{minipage}
	\begin{minipage}[c]{.46\linewidth}
		$\lim_{t\ifty}\lim_{T\ifty}\PR_{(x,\,n)} \lc (X, N)_t \in (dx, dn)\Bv  T < \tau_\partial\rc$
		
		\hspace*{1cm}$=\beta(dx, dn)$
	\end{minipage}
	
	\begin{minipage}[c]{.46\linewidth}
		\includegraphics[width=70mm, height = 40mm]{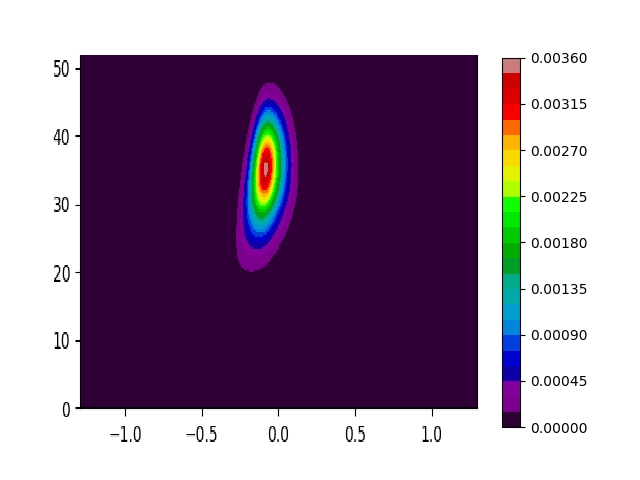}
		
		Critical regime of adaptation
	\end{minipage}
	\begin{minipage}[c]{.46\linewidth}
		\includegraphics[width=70mm, height = 40mm]{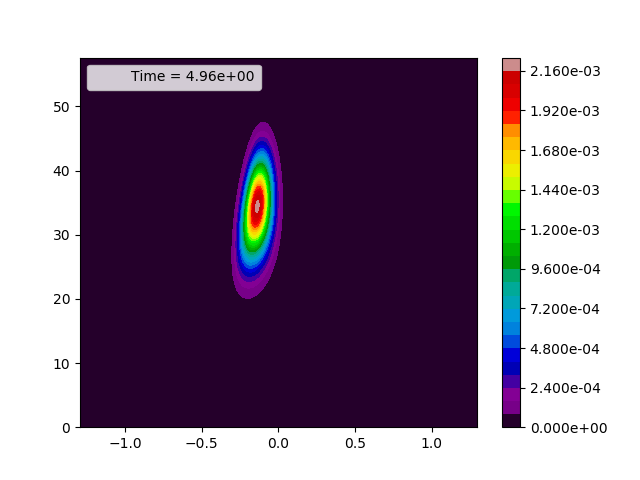}
			\mbox{Selection through the extinction of populations}
	\end{minipage}
	
	\vspace{0.3cm}
	Given the very different profiles obtained for the QSD, it is quite remarkable that the quasi-ergodic measures are similar as much.
	In particular, we can see that the histories of the surviving populations are still shaped by the maintenance of these populations at large sizes with almost optimal traits, even when such traits are very rare according to the QSD.
\end{figure}
\newpage
\begin{center}
\begin{figure}
	\textbf{Profiles of the survival capacity}
	
	\begin{minipage}[c]{.46\linewidth}
		\includegraphics[width=70mm, height = 36mm]{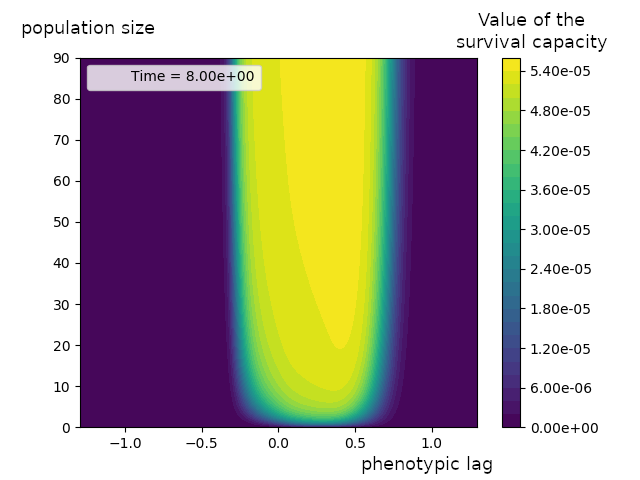}
		\label{Eta60}
		Case of a spontaneous adaptation
	\end{minipage}
	\begin{minipage}[c]{.46\linewidth}
		$\heig(x, n):=
		\underset{t\rightarrow\infty}{\lim}
		\dfrac{\PR_{(x,\, n)} (t < \tau_\partial)  }
		{\PR_{\alpha} (t< \tau_\partial) }$
	\end{minipage}
	
	\begin{minipage}[c]{.46\linewidth}
		\includegraphics[width=70mm, height = 36mm]{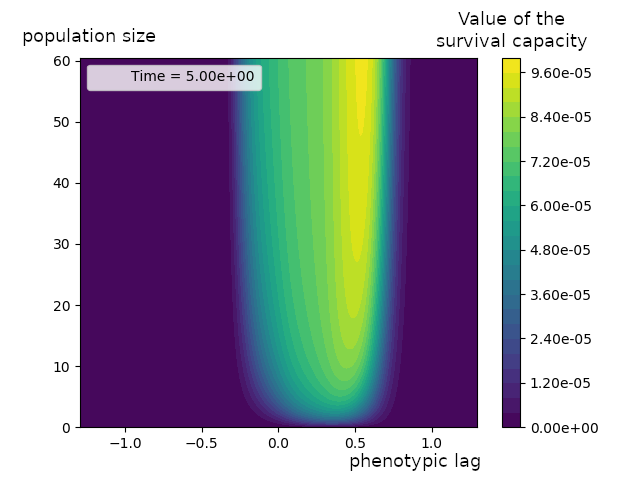}
		\label{Eta61}
		Critical regime of adaptation
	\end{minipage}
	\begin{minipage}[c]{.46\linewidth}
		\includegraphics[width=70mm, height = 36mm]{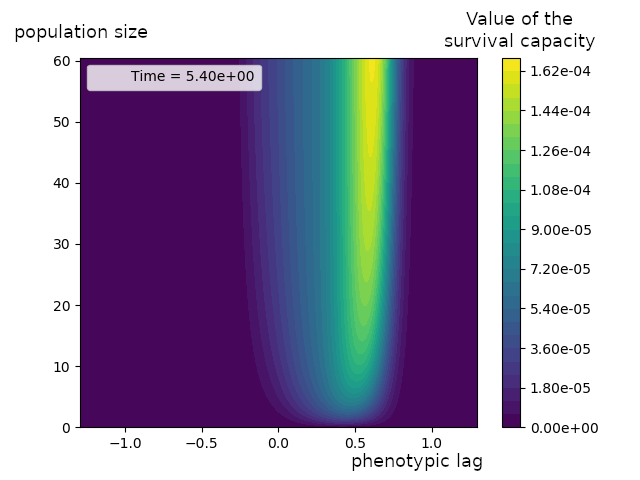}
		\label{Eta62}
		\mbox{Selection through the extinction of populations}
	\end{minipage}
\vspace{0.5cm}

		\textbf{Profiles 3D of the survival capacity}
		
		\begin{minipage}[c]{.30\linewidth}
			\includegraphics[width=40mm, height = 25mm]{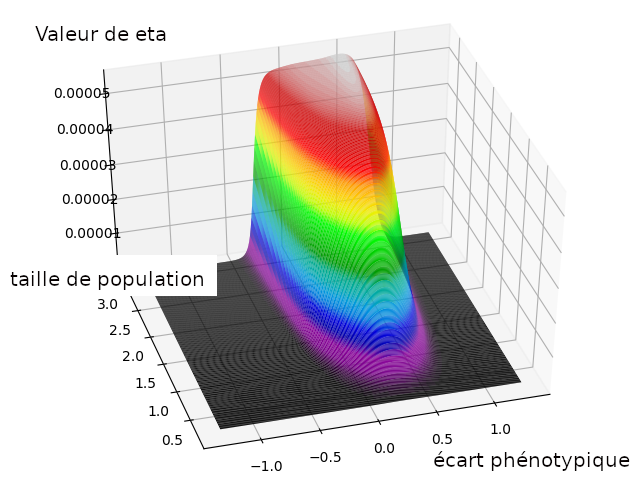}
			
			Spontaneous adaptation
		\end{minipage}
		\begin{minipage}[c]{.30\linewidth}
			\includegraphics[width=40mm, height = 25mm]{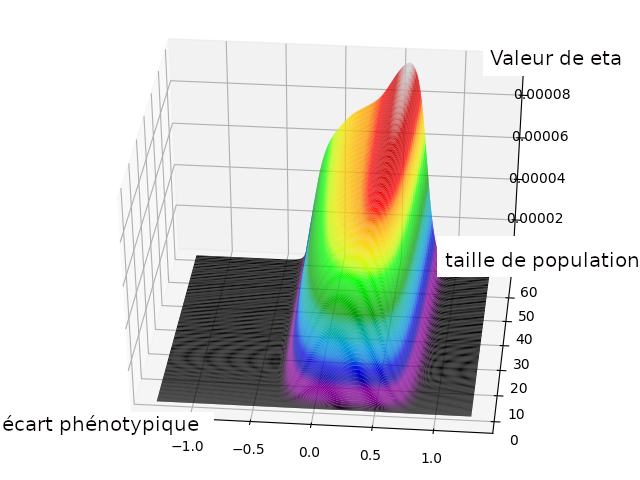}
			
			Critical adaptation
		\end{minipage}
		\begin{minipage}[c]{.30\linewidth}
			\includegraphics[width=40mm, height = 25mm]{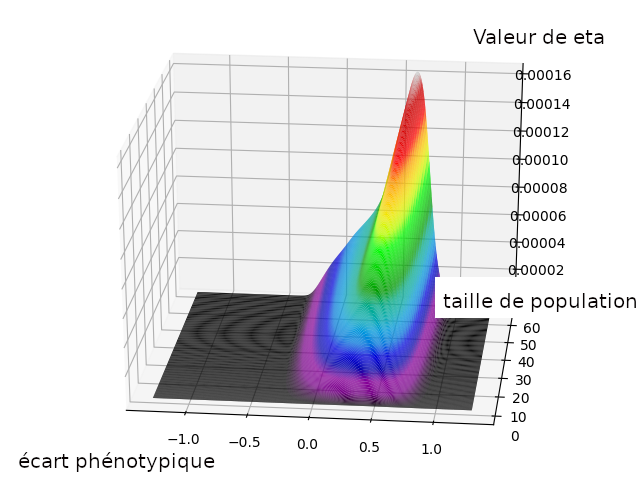}
			
			\mbox{Adaptation through extinction}
		\end{minipage}
	\end{figure}
\end{center}

\end{document}
